\documentclass[12pt]{amsart}
\usepackage{amssymb,latexsym}
\usepackage{fullpage}
\usepackage{xypic}
\usepackage[mathcal]{euscript}

\newcommand{\Mdef}[2]{\newcommand{#1}{\relax \ifmmode #2 \else $#2$\fi}}

\newcommand{\cok}{\mathrm{cok}}

\newcommand{\sm }{\wedge}

\newcommand{\tensor}{\otimes}

\newcommand{\map}{\mathrm{map}}
\newcommand{\Hom}{\mathrm{Hom}}

\newcommand{\Ext}{\mathrm{Ext}}
\Mdef{\bhom}{\mathbf{\hat{H}om}}

\Mdef{\Mod}{\mathrm{mod}}

\newcommand{\st}{\; | \;}

\newcommand{\lr}[1]{\langle #1\rangle}

\newtheorem{thm}{Theorem}[section]
\newtheorem{lemma}[thm]{Lemma}
\newtheorem{prop}[thm]{Proposition}
\newtheorem{cor}[thm]{Corollary}
\theoremstyle{definition}

\newtheorem{defn}[thm]{Definition}

\newtheorem{warning}[thm]{Warning}

\newtheorem{example}[thm]{Example}
\newtheorem{remark}[thm]{Remark}

\newtheorem{convention}[thm]{Convention}

\newcommand{\ppf}{{\flushleft {\bf Proof:}}}
\newcommand{\qqed}{\qed \\[1ex]}
\renewenvironment{proof}[1][\hspace*{-.8ex}]{\noindent {\bf Proof #1:\;}}{\qqed}

\Mdef{\PH} {\Phi^H}
\Mdef{\PK} {\Phi^K}
\Mdef{\PL} {\Phi^L}
\Mdef{\PT} {\Phi^{\T}}

\Mdef{\ef}{E{\cF}_+}
\Mdef{\etf}{\tilde{E}{\cF}}
\Mdef{\eg}{E{G}_+}
\Mdef{\etg}{\tilde{E}{G}}
\Mdef{\tf}{\T / \! \!  / {\cF}_+}

\newcommand{\piT}{\pi^\mathbb{T}}
\newcommand{\piA}{\pi^{\cA}}

\Mdef{\infl}{\mathrm{inf}}
\Mdef{\defl}{\mathrm{def}}
\Mdef{\res}{\mathrm{res}}
\Mdef{\ind}{\mathrm{inf}}

\Mdef{\univ}{\mathcal{U}}

\Mdef{\Fp}{\mathbb{F}_p}
\Mdef{\Zpinfty}{\Z /p^{\infty}}
\Mdef{\Zpadic}{\Z_p^{\wedge}}

\newcommand{\Ga}{\G_a}
\newcommand{\Gm}{\G_m}

\newcommand{\bi}{\begin{itemize}}
\newcommand{\be}{\begin{enumerate}}
\newcommand{\bc}{\begin{center}}
\newcommand{\bd}{\begin{description}}
\newcommand{\ei}{\end{itemize}}
\newcommand{\ee}{\end{enumerate}}
\newcommand{\ec}{\end{center}}
\newcommand{\ed}{\end{description}}
\newcommand{\dichotomy}[2]{\left\{ \begin{array}{ll}#1\\#2 \end{array}\right.}
\newcommand{\trichotomy}[3]{\left\{ \begin{array}{ll}#1\\#2\\#3 \end{array}\right.}

\newcommand{\adjunction}[4]{
\diagram
#1:#2 \rrto<0.7ex> &&
#3  \llto<0.7ex> :#4 
\enddiagram}
\newcommand{\lra}{\longrightarrow}
\newcommand{\lla}{\longleftarrow}

\newcommand{\spec}{\mathrm{spec}}

\Mdef{\we}{\mathbf{we}}
\Mdef{\fib}{\mathbf{fib}}
\Mdef{\cof}{\mathbf{cof}}
\Mdef{\BI}{\mathcal{BI}}

\newcommand{\ann}{\mathrm{ann}}

\newcommand{\cofibre}{\mathrm{cofibre}}

\newcommand{\ilim}{\mathop{ \mathop{\mathrm{lim}} \limits_\leftarrow} \nolimits}
\newcommand{\colim}{\mathop{  \mathop{\mathrm {lim}} \limits_\rightarrow} \nolimits}

\Mdef{\A}{\mathbb{A}}
\Mdef{\B}{\mathbb{B}}
\Mdef{\C}{\mathbb{C}}
\Mdef{\D}{\mathbb{D}}
\Mdef{\E}{\mathbb{E}}
\Mdef{\T}{\mathbb{T}}
\Mdef{\F}{\mathbb{F}}
\Mdef{\G}{\mathbb{G}}
\Mdef{\I}{\mathbb{I}}
\Mdef{\N}{\mathbb{N}}
\Mdef{\bbP}{\mathbb{P}}
\Mdef{\Q}{\mathbb{Q}}
\Mdef{\R}{\mathbb{R}}
\Mdef{\bbS}{\mathbb{S}}
\Mdef{\Z}{\mathbb{Z}}

\Mdef{\bA}{\mathbb{A}}
\Mdef{\bB}{\mathbb{B}}
\Mdef{\bC}{\mathbb{C}}
\Mdef{\bD}{\mathbb{D}}
\Mdef{\bE}{\mathbb{E}}
\Mdef{\bF}{\mathbb{F}}
\Mdef{\bG}{\mathbb{G}}
\Mdef{\bH}{\mathbb{H}}
\Mdef{\bI}{\mathbb{I}}
\Mdef{\bJ}{\mathbb{J}}
\Mdef{\bK}{\mathbb{K}}
\Mdef{\bL}{\mathbb{L}}
\Mdef{\bM}{\mathbb{M}}
\Mdef{\bN}{\mathbb{N}}
\Mdef{\bO}{\mathbb{O}}
\Mdef{\bP}{\mathbb{P}}
\Mdef{\bQ}{\mathbb{Q}}
\Mdef{\bR}{\mathbb{R}}
\Mdef{\bS}{\mathbb{S}}
\Mdef{\bT}{\mathbb{T}}
\Mdef{\bU}{\mathbb{U}}
\Mdef{\bV}{\mathbb{V}}
\Mdef{\bW}{\mathbb{W}}
\Mdef{\bX}{\mathbb{X}}
\Mdef{\bY}{\mathbb{Y}}
\Mdef{\bZ}{\mathbb{Z}}

\Mdef{\cA}{\mathcal{A}}
\Mdef{\cB}{\mathcal{B}}
\Mdef{\cC}{\mathcal{C}}
\Mdef{\mcD}{\mathcal{D}} % Something funny about \cD.
\Mdef{\cE}{\mathcal{E}}
\Mdef{\cF}{\mathcal{F}}
\Mdef{\cG}{\mathcal{G}}
\Mdef{\mcH}{\mathcal{H}} % There's something funny about \cH: it 
\Mdef{\cI}{\mathcal{I}}
\Mdef{\cJ}{\mathcal{J}}
\Mdef{\cK}{\mathcal{K}}
\Mdef{\mcL}{\mathcal{L}}% There's something funny about \cL: it 
\Mdef{\cM}{\mathcal{M}}
\Mdef{\cN}{\mathcal{N}}
\Mdef{\cO}{\mathcal{O}}
\Mdef{\cP}{\mathcal{P}}
\Mdef{\cQ}{\mathcal{Q}}
\Mdef{\mcR}{\mathcal{R}}% There's something funny about \cR: it 
\Mdef{\cS}{\mathcal{S}}
\Mdef{\cT}{\mathcal{T}}
\Mdef{\cU}{\mathcal{U}}
\Mdef{\cV}{\mathcal{V}}
\Mdef{\cW}{\mathcal{W}}
\Mdef{\cX}{\mathcal{X}}
\Mdef{\cY}{\mathcal{Y}}
\Mdef{\cZ}{\mathcal{Z}}

\Mdef{\tA}{\tilde{A}}
\Mdef{\tB}{\tilde{B}}
\Mdef{\tC}{\tilde{C}}
\Mdef{\tE}{\tilde{E}}
\Mdef{\tH}{\tilde{H}}
\Mdef{\tK}{\tilde{K}}
\Mdef{\tL}{\tilde{L}}
\Mdef{\tM}{\tilde{M}}
\Mdef{\tN}{\tilde{N}}
\Mdef{\tP}{\tilde{P}}
\Mdef{\tR}{\tilde{R}}

\Mdef{\Ab}{\overline{A}}
\Mdef{\Bb}{\overline{B}}
\Mdef{\Cb}{\overline{C}}
\Mdef{\Db}{\overline{D}}
\Mdef{\Hb}{\overline{H}}
\Mdef{\Ib}{\overline{I}}
\Mdef{\Kb}{\overline{K}}
\Mdef{\Lb}{\overline{L}}
\Mdef{\Mb}{\overline{M}}
\Mdef{\Nb}{\overline{N}}
\Mdef{\Qb}{\overline{Q}}
\Mdef{\Tb}{\overline{T}}
\Mdef{\Ub}{\overline{U}}

\Mdef{\hb}{\overline{h}}
\Mdef{\qb}{\overline{q}}
\Mdef{\rb}{\overline{r}}
\Mdef{\tb}{\overline{t}}
\Mdef{\ub}{\overline{u}}
\Mdef{\vb}{\overline{v}}
\Mdef{\hc}{\hat{c}}
\Mdef{\he}{\hat{e}}
\Mdef{\hf}{\hat{f}}
\Mdef{\hA}{\hat{A}}
\Mdef{\hH}{\hat{H}}
\Mdef{\hJ}{\hat{J}}
\Mdef{\hM}{\hat{M}}
\Mdef{\hP}{\hat{P}}
\Mdef{\hQ}{\hat{Q}}
\Mdef{\hR}{\hat{R}}
\Mdef{\hS}{\hat{S}}
\Mdef{\hT}{\hat{T}}
\Mdef{\hU}{\hat{U}}
\Mdef{\hV}{\hat{V}}
\Mdef{\hW}{\hat{W}}
\Mdef{\hX}{\hat{X}}
\Mdef{\hY}{\hat{Y}}
\Mdef{\hZ}{\hat{Z}}

\Mdef{\bolda}{\mathbf{a}}
\Mdef{\boldb}{\mathbf{b}}
\Mdef{\boldk}{\mathbf{k}}
\Mdef{\boldD}{\mathbf{D}}
\Mdef{\boldR}{\mathbf{R}}

\Mdef{\fm}{\frak{m}}

\Mdef{\eps}{\epsilon}

\newcommand{\EAA}{E(A \times A )}
\newcommand{\cOAxA}{\cO_{A \times A} }
\newcommand{\AxA}{A \times A }

\newcommand{\TxT}{\T \times \T }
\newcommand{\DTtwo}{\Delta \T[2]}
\newcommand{\TxTbar}{\overline{\T \times \T} }

\newcommand{\xihat}{\hat{\xi}}
\newcommand{\units}{\mathrm{Units}}

\renewcommand{\div}{\mathrm{div}}

\newcommand{\tp}{\mathrm{tp}}
\newcommand{\zar}{\mathrm{Zar}}
\newcommand{\an}{\mathrm{an}}
\newcommand{\cOan}{\cO^{an}}
\newcommand{\Omegaan}{\Omega_{an}}

\newcommand{\groj}{\mathrm{Groj}}
\newcommand{\grojtilde}{\widetilde{\mathrm{Groj}}}
\newcommand{\adele}{ad\`ele}
\newcommand{\wtET}{\widetilde{E}^{\T}}

\newcommand{\EAmod}{\mbox{$EA$-mod}}

\newcommand{\cMA}{\cM_A}
\newcommand{\SA}{\mathcal{S}_A}
\newcommand{\SAt}{\mathcal{S}_A^t}
\newcommand{\tpsheafA}{\mathrm{sheaves}/A}
\newcommand{\tpOmod}{\mbox{$\cO_A$-mod}}

\newcommand{\proj}{\mathrm{Proj}}
\newcommand{\ord}{\mathrm{ord}}
\newcommand{\kk}{k}

\newcommand{\tETz}{\tilde{E}_{\T}^0}
\newcommand{\tET}{\tilde{E}_{\T}}
\newcommand{\ETz}{E_{\T}^0}

\newcommand{\ETc}{E_{\T}^*}
\newcommand{\EATc}{EA_{\T}^*}

\newcommand{\tEATi}{\widetilde{EA}_{\T}^i}

\newcommand{\tEAThz}{\widetilde{EA}^{\T}_0}
\newcommand{\tEAThone}{\widetilde{EA}^{\T}_1}

\newcommand{\tEATh}{\widetilde{EA}^{\T}}
\newcommand{\tEATc}{\widetilde{EA}_{\T}}

\newcommand{\EAT}{EA_{\T}^*}

\Mdef{\at}{\cA_t}
\Mdef{\At}{\cA_t}
\Mdef{\As}{\cA_s}
\Mdef{\as}{\cA_s}
\Mdef{\dat}{D(\at)}
\Mdef{\das}{D(\as)}
\Mdef{\atA}{\cA A_t}
\Mdef{\asA}{\cA_s}
\Mdef{\cAsfs}{\cA^{sf}_s}
\Mdef{\cAsft}{\cA^{sf}_t}
\Mdef{\datA}{D(\atA)}
\Mdef{\dasA}{D(\asA)}
\newcommand{\dimC}{\mathrm{dim}_{\C}}

\newcommand{\phiA}{\phi A}

\newcommand{\QAlr}[1]{QA\langle #1 \rangle}
\newcommand{\TAlr}[1]{TA\langle #1 \rangle}

\newcommand{\Alr}[1]{A\langle #1 \rangle}
\newcommand{\Etlr}[1]{\tilde{E}\langle #1 \rangle}

\newcommand{\philr}[1]{\phi_{#1}}%{\phi\langle #1 \rangle}

\newcommand{\ShHom}{\mathcal{HOM}}

\newcommand{\cOiGt}{\cO (\infty \G[tors])}
\newcommand{\Gpointy}[1]{\G \langle #1 \rangle}

\newcommand{\EsubH}{E[\subseteq H]_+}
\newcommand{\TA}{TA}

\newcommand{\siftyzn}{S^{\infty z^n}}
\newcommand{\subH}{[\subseteq H]}

\newcommand{\ek}{E\langle K \rangle}
\newcommand{\efp}{E\cF_+}
\newcommand{\cFlrpi}{\cF \langle \pi \rangle}

\newcommand{\Ghat}{\widehat{\bG}}
\newcommand{\Gnd}{\bG_{nd}}
\newcommand{\Gndlr}[1]{\Gnd \langle #1 \rangle}
\newcommand{\EG}{E\G}
\newcommand{\EGT}{E\G_{\T}}
\newcommand{\EGm}{E\Gm}

\newcommand{\EGa}{E\Ga}

\newcommand{\ET}{E\T}
\newcommand{\ETp}{E\T_+}

\newcommand{\HomcO}{\mathrm{Hom}_{\cO}}
\newcommand{\HomOs}{\mathrm{Hom}_{\cO_s}}
\newcommand{\tensorcO}{\otimes_{\cO}}
\newcommand{\HomAs}{\mathrm{Hom}_{\cA_s}}
\newcommand{\HomK}{\mathrm{Hom}_{\cK}}
\newcommand{\ExtAs}{\mathrm{Ext}_{\cA_s}^1}
\newcommand{\HomAt}{\mathrm{Hom}_{\cA_t}}

\newcommand{\HomEA}{\mathrm{Hom}_{EA}}
\newcommand{\ExtEA}{\mathrm{Ext}_{EA}}
\newcommand{\HomEAtp}{\mathrm{Hom}_{EA,\tp}}
\newcommand{\ExtEAtp}{\mathrm{Ext}_{EA,\tp}}

\newcommand{\cOit}{\cO (\infty tors)}
\newcommand{\GcOit}{\Gamma \cO (\infty tors)}
\newcommand{\GcO}{\Gamma \cO}
\newcommand{\inftors}{(\infty tors)}
\newcommand{\cEi}{\cE^{-1}}
\newcommand{\psibar}{\overline{\psi}}
\newcommand{\Wbar}{\overline{W}}

\newcommand{\LZ}{\Lambda Z}
\newcommand{\MLZ}{M\Lambda Z}

\renewcommand{\L}{\mathbb{L}}
\renewcommand{\tf}{t^{\cF}_*}

\newcommand{\Dtp}{D_{tp}}
\newcommand{\cOtpmod}{\mbox{$\cO^{\tp}_A$-mod}}
\newcommand{\cOzarmod}{\mbox{$\cO^{\zar}_A$-mod}}
\newcommand{\shvtp}{\mathrm{Shv}_A^{\tp}}
\newcommand{\shvzar}{\mathrm{Shv}_A^{\zar}}

\newcommand{\cOtp}{\cO^{\tp}_A}
\newcommand{\cOzar}{\cO^{\zar}_A}
\newcommand{\mapho}[1]{[#1]_{Ho(EA)}}
\newcommand{\elr}[1]{E\langle #1 \rangle}
\newcommand{\XT}{X^{\T}}
\newcommand{\XC}{X^C}
\newcommand{\cKF}{\cK_{\cF}}

\begin{document}
\title[Rational $S^1$-equivariant elliptic cohomology.]
{Rational $S^1$-equivariant elliptic cohomology.\\
(Version 5.2)}
\author{J.P.C.Greenlees}
\address{Department of Pure Mathematics, Hicks Building,
Sheffield S3 7RH. UK.}
\email{j.greenlees@sheffield.ac.uk}
%\author{M.J.Hopkins}
%\address{Department of Mathematics, MIT, Cambridge, MA 02139-4307, USA.}
%\email{mjh@math.mit.edu}
%\author{I.Rosu}
%\address{Department of Mathematics, MIT, Cambridge, MA 02139-4307, USA.}
%\email{ioanid@math.mit.edu}
\thanks{
%JPCG and MJH are grateful to the Mathematisches Forschungsinstitut
% Oberwolfach  for the opportunity to talk at the 1998 Homotopietheorie
% meeting, and to JPCG's audience for tolerating the resulting delay.
The author is grateful to M.Ando, M.J.Hopkins, I.Rosu and N.P.Strickland 
for useful conversations,  and to H.R.Miller and various referees for 
stimulating comments on precursors of this paper.}
\date{}
\begin{abstract}
We give a functorial construction of  a rational $S^1$-equivariant 
cohomology theory from  an elliptic curve $A$ 
equipped with suitable coordinate data. The 
elliptic curve may be recovered from the cohomology theory; indeed, 
the  value of the cohomology theory on the compactification of an 
$S^1$-representation is given by the sheaf cohomology of a suitable 
line bundle on $A$. This suggests the construction: by considering
functions on the elliptic curve with specified poles one may 
write down the  representing $S^1$-spectrum
in the author's algebraic model of rational $S^1$-spectra
\cite{s1q}. 

The construction extends to give an equivalence of 
categories between the homotopy category of module $S^1$-spectra
 over the representing spectrum and a derived category of sheaves
of modules over the structure sheaf of $A$.
\end{abstract}
\maketitle

\begin{center}
Email: J.Greenlees@sheffield.ac.uk
\end{center}
\setcounter{tocdepth}{1}
\tableofcontents

\section{Introduction.}
\label{secIntro}

Two of the most important topological cohomology theories 
are associated to one dimensional group schemes: ordinary cohomology 
is associated to the additive group, and $K$-theory is associated to 
the multiplicative group. This connection is most transparent in
the equivariant context, and because the group schemes are one dimensional 
it is enough to consider a one dimensional group of equivariance: 
the circle group $\T$.

Beginning with   ordinary cohomology,  we use  the Borel construction
to define an equivariant theory for $\T$-spaces $X$ 
 by $H_{\T}^*(X)=H^*(E\T \times_{\T}X)$. The coefficient
ring $H_{\T}^*=H^*(B\T)\cong \Z [x]$ inherits a coproduct from the map 
$B\T \times B\T \lra B\T$ classifying tensor product of line bundles, 
and the resulting Hopf algebra represents the additive group. 

This construction works equally well for any complex oriented theory.
For instance if we let $z$ denote the natural representation of the
circle group on the complex numbers, 
$K$-theory of the Borel construction has coefficient
ring $K^0(BT)=\Z [[y]]$, with $y=1-z$, 
and this represents the multiplicative formal group. 
However,  by working with the correct equivariant theory we may obtain 
the uncompleted version. Indeed, the coefficient ring 
$K_{\T}^0=\Z [z ,z^{-1}]$ of Atiyah-Segal 
equivariant $K$-theory acquires a coproduct from the group 
multiplication $\T \times \T \lra \T$, and the resulting Hopf 
algebra represents the multiplicative group. 

Elliptic cohomology  was first defined \cite{Landweber, LRS} 
as a non-equivariant
complex oriented cohomology theory whose associated formal group 
is the completion of an elliptic curve around the identity.
It is therefore natural to hope for an equivariant cohomology theory 
giving the associated  elliptic curve $A$ without completion.
It is the purpose of the  present paper to construct such a theory 
over the rationals and establish its basic properties. 
The most obvious new feature is that $A$ is not affine, and one of our 
main tasks is to  elucidate the connection between  the cohomology theory 
and the elliptic curve

 A programme to extend this work to higher
dimensional abelian varieties and higher dimensional tori is
underway \cite{atnq1,atnq2,qtnq,GEGA}.

In concrete terms, the main purpose of this paper is to construct a rational
$\T$-equivariant cohomology theory $\EAT (\cdot )$ associated to any elliptic
curve $A$ over a $\Q$-algebra. The construction is compatible with 
base change, and the properties of the cohomology theory 
%%FIELD
when we work over a field may be summarized as follows;
we give full details in Section \ref{sec:EA} below.

\begin{thm}
\label{IntroMainThm}
For any elliptic curve $A$ over a field $k$ of characteristic 0, 
there is a 2-periodic, multiplicative, rational $\T$-equivariant cohomology 
theory $\EAT (\cdot )$. 
The value on the one point compactification
$S^W$ of a complex representation $W$ of $\T$ with $W^{\T}=0$
is given as the sheaf cohomology 
of a line bundle $\cO (-D(W))$. To describe this,  
we write $A[n]$ for the divisor of points of order dividing $n$ in $A$. 
If $W=\sum_n a_n z^n$, we consider the divisor
$D(W)=\sum_n a_n A[n]$, and the  associated line bundle
$\cO (-D(W))$ on $A$. The cohomology of $S^W$ is given by 
$$\tEATc^i (S^{W}) \cong H^i(A; \cO (-D(W)) \mbox{ for } i=0,1$$
and the  homology by 
$$\tEATh_{-i} (S^{W}) \cong H^i(A; \cO (D(W)) \mbox{ for } i=0,1. $$
In particular, the coefficient ring is 
$$EA^{\T}_*=k[u,u^{-1}]\tensor \Lambda (\tau)$$
where $u$ (of degree 2) is a generator of $H^0(A;\Omega)$
(i.e., a nowhere zero, regular differential) and $\tau$
(of degree $-1$) is a generator of $H^1(A;\cO)$.
\end{thm}

\begin{remark}
The above properties do not quite determine the cohomology theory.
The cohomology theory depends on one  auxiliary piece of data: a
{\em coordinate} $t_e$ on $A$. This is a
function vanishing to the first order at the identity, whose zeroes
and poles are all at points of finite order. The construction is 
natural for isomorphisms of the data $(A, k, t_e)$.
All three of the inputs, $A$, $k$ and $t_e$ 
can be recovered from the cohomology of suitable spaces.
\end{remark}

\begin{remark}
For general spaces there is a Hasse long exact sequence 
describing how to calculate elliptic cohomology. A precise 
statement is given in \ref{prop:Hasse}, but the idea is that, 
just as the arithmetic Hasse square 
recovers global data from completions at various primes $p$, compatible
in the rationalization, the Hasse sequence recovers elliptic cohomology 
from Borel cohomology of $C$-fixed points with 
coefficients in completions of local rings at points of 
 order $|C|$, compatible with the cohomology of the $\T$-fixed 
points with coefficients in meromorphic functions.
\end{remark}

The first version of $\T$-equivariant elliptic cohomology was
constructed by Grojnowksi in 1994 \cite{grojnowski}. He was
interested in implications for the representation theory of
certain  elliptic algebras: these implications are the
subject of the work of Ginzburg-Kapranov-Vasserot \cite{GKV} and
the context is explained further in \cite{Ginzburg}. For this
purpose it was sufficient to construct a theory on finite
complexes taking  values in analytic sheaves over the elliptic curve. Later
Rosu  \cite{rosu} used this sheaf-valued theory to give a proof of
Witten's rigidity theorem for the equivariant Ochanine genus of a
spin manifold with non-trivial $\T$-action, and Ando-Basterra do the same
for the Witten genus \cite{AndoBasterra}.  Ando \cite{Ando}
has related the sheaf valued theory to the representation theory of loop 
groups.

However, to exploit the theory fully,  it is essential to have a theory defined
on general $\T$-spaces and $\T$-spectra, and to have a conventional 
group-valued theory represented by a $\T$-spectrum $EA$. 
This allows one to use the
full apparatus of equivariant stable homotopy theory. For example,  twisted 
pushforward maps are
immediate consequences of Atiyah duality; in more concrete terms,  it allows
one to calculate the theory on free loop spaces, and to describe algebras of 
operations. It is also likely to be useful in constructing an integral version
of the theory, and we hope it may also prove useful in the continuing search 
for a geometric definition of  elliptic cohomology.
The theory we construct has these desirable properties, whilst retaining
a very close connection with the geometry of the underlying elliptic curve.

Returning to the geometry, a very appealing  feature  is that although our 
theory is group valued, the original curve can still be recovered from the 
cohomology theory. It is also notable that the  earlier sheaf theoretic 
constructions work over larger
rings and certainly require the coefficients to contain roots of unity:
the loss of information can be illustrated by comparing the
rationalized representation ring $R(C_n)=\Q [z]/(z^n-1)$ (with
components corresponding to {\em subgroups} of $C_n$) to the complexified
representation ring, isomorphic to the character ring  $\map (C_n, \C)$
(with components corresponding to the {\em elements} of $C_n$).

Finally,  the ingredients of the model are very natural invariants
of the curve given by sheaves of functions with specified poles at
points of finite order: Definition \ref{defn:EA}  simply writes down the 
representing object in terms of these, and readers already familiar 
with elliptic curves and the model of \cite{s1q} may wish to look at 
this  immediately. In fact the algebraic model of \cite{s1q}
gives a generic de Rham model for all $\T$-equivariant theories,
and the models of elliptic cohomology theories highlight this
geometric structure. These higher de Rham models should allow
applications in the same spirit as those made for de Rham models
of ordinary cohomology and $K$-theory \cite{guillemin}.
%There is a conjectural model for rational $T^n$-spectra for an $n$-dimensional
%torus $T^n$, and a programme to verify that it is indeed a model is underway.
%It is easy to identify the model of a $T^n$-equivariant cohomology theory
%in this algebraic category.

In fact, we are able to go beyond constructing a particular cohomology theory 
$EA_{\T}^*(\cdot )$ and  establish  an equivalence between a derived
category of sheaves over the elliptic curve and 
cohomology theories which are modules over $EA_{\T}^*(\cdot )$.
Because homotopy theory only sees points of finite order, we use the 
torsion point topology on the elliptic curve consisting of complements 
of sets of points of finite order, which  is coarser than the 
Zariski topology, and because a
$\T$-equivariant homotopy equivalence is an equivalence in $H$-equivariant
fixed points for all subgroups $H$, 
the maps inverted in forming $\Dtp (\cOtpmod)$
are those inducing isomorphisms of $H^*(A; \cO_A (D(W)) \tensor_{\cO} (\cdot))$
for all representations $W$ with $W^{\T}=0$.

\begin{thm} 
\label{IntroEquivalence}
The representing object $EA_a$ in the algebraic category may be taken to 
be a commutative ring, and there is an equivalence
$$\Dtp (\cOtpmod) \simeq D_{\T}(\mbox{$EA_a$-mod})$$
between derived categories of sheaves of $\cO_A$-modules on $A$ and 
$EA_a$-modules. These categories both have relative injective dimension 
1, so that maps are calculated by a short exact sequence from Hom 
and Ext groups in an abelian category.

The corresponding $\T$-spectrum $EA$ is a ring up to homotopy, 
and the above equivalence classifies homotopy $EA$-module spectra up 
to equivalence as $\cO_A$-modules up to isomorphism.
Using the result of \cite{qtnq}, that  $EA$ can be realized as a strictly 
commutative ring spectrum, the right
hand side may be replaced the derived category of $EA$-module $\T$-spectra,
and morphisms of module spectra are thereby also classified.
\end{thm}

This is proved in Section \ref{sec:equivalence}.
Our construction directly models the representing ring spectrum $EA$ in 
the author's algebraic model $\cA_s$ of rational $\T$-spectra \cite{s1q}.
We describe the abelian category $\cA_s$ in detail in  
Section \ref{sec:ModelT}, but it can be viewed as a category of sheaves
over the space of closed subgroups of $\T$ \cite{atnq1}. The equivalence
is obtained from functors at the level of abelian categories, and
(Theorem \ref{Grojnowskicons}) 
Grojnowski's sheaf $\groj (X)$ associated to a compact $\T$-manifold
$X$ is obtained by applying the functor to the function spectrum $EA$-module 
$F(X,EA)$, and then changing to the analytic topology. Thus, for a compact
$\T$-manifold $X$,  there is a short exact sequence 
$$0 \lra \Sigma H^1(A; \groj (X))
\lra EA_{\T}^*(X) \lra H^0(A; \groj (X)) \lra 0$$
relating the cohomology of Grojnowksi's sheaf to $EA_{\T}^*(X)$.

By way of motivation,  we will discuss the way that a $\T$-equivariant
 cohomology theory is associated to several other geometric objects.
Perhaps most familiar is the complete case discussed in
Section \ref{secFormal}, where the Borel theory
for a complex oriented cohomology theory is associated to a formal group.
Amongst global groups, the additive and multiplicative ones are the
simplest, and in Appendix \ref{sec:GaGm} we describe how they give rise to
ordinary Borel cohomology and equivariant $K$-theory; the behaviour of the
construction on the non-split torus is also notable.
%This construction is notable in that it gives a  construction of
%equivariant cohomology theories from oriented 1-dimensional group schemes
% which is {\em functorial for isomorphisms}. It is also
%functorial for certain isogenies as explained in \ref{isogeny}.

We have divided the paper into six parts. Part 1 explains how 
equivariant cohomology theories ought to be related to group schemes. 
Part 2 provides prerequisites on rational $\T$-equivariant cohomology 
theories. Part 3 provides prerequisites on elliptic curves. Part 4
is extremely short, and just contains the construction. Part 5 
describes some properties of the theory. Part 6 builds on the 
construction to give an equivalence between a derived category of 
sheaves over $A$ and a derived category of $\T$-spectra.
The appendix re-examines equivariant $K$-theory from the present point
of view. \\[2ex]

\noindent
{\bf Historical note.}
Early versions of this paper were under joint authorship with 
M.J.Hopkins and I.Rosu. This reflected the fundamental influence 
of their ideas, in the expectation that they would continue to 
be part of the project. To the disappointment of all of us, 
circumstances prevented this, and the other authors  withdrew.

Rosu's emphasis on the sheaf associated to a sphere \cite{rosu} was 
significant. When the author first heard it at the 1997 Glasgow workshop 
on elliptic cohomology, he  believed this would necessitate
representing elliptic cohomology by sheaves of spectra. 
However it led  Hopkins towards his vision that a result like
the Theorem \ref{IntroMainThm} proved here should be true. Work on 
the present paper began after a breathless conversation between 
the author and Hopkins in Oberwolfach at the 1998 Homotopietheorie meeting.

The present paper is  Version 5.2 of the preprint.

\part{Equivariant cohomology theories and group schemes.}
In Part 1 we describe how equivariant cohomology theories
and group schemes are related in ideal circumstances. We 
begin with the familiar
example of formal groups and complex oriented theories, and then
explore how this correspondence should be extended. 

\section{Formal groups from complex oriented theories.}
\label{secFormal}

The purpose of this section is to recall that any complex
orientable cohomology theory $E^*(\cdot )$ determines a one
dimensional, commutative formal group $\Ghat$ and to explain how
the cohomology of various spaces can be described in terms of the
geometry of $\Ghat$. This is well known (see especially
\cite{AHS}) but it introduces the
geometric language, and motivates our main construction, which
{\em over the rationals} reverses the process by using 
this geometric data to construct the cohomology theory.
Indeed, we will show that the machinery of \cite{s1q} permits  a
 construction of a 2-periodic rational $\T$-equivariant
cohomology theory $E\G^*_{\T}(\cdot )$ from a  one dimensional
group scheme $\G$ over a $\Q$-algebra, functorial in $\G$ with some
additional data. Furthermore, the construction is
reversible in the sense that $\G$ can be recovered from
$E\G^*_{\T}(\cdot )$. The most interesting case of this is when
$\G$ is an elliptic curve, but the affine case is treated in 
Appendix \ref{sec:GaGm}.

\subsection{Geometry of formal groups.}
Before bringing the cohomology theory into the picture, we
introduce the geometric language. When all schemes are affine,
the geometric language is equivalent to the ring theoretic language, 
and all geometric  statements can be given meaning by translating them 
to algebraic ones. It  is traditional in topology to stick to algebra, 
but to prepare for the case  of an elliptic curve, 
we will use the geometric language.

A one dimensional commutative formal group law over a ring $k$ is
a commutative and associative coproduct on  the complete
topological $k$-algebra $k[[y]]$. Equivalently, it is a complete topological
Hopf $k$-algebra $\cO$ together with an element $y \in \cO$ so that
$\cO =k[[y]]$. A topological Hopf $k$-algebra $\cO$ for which such
a $y$ exists  is the ring of functions on a one dimensional
commutative formal group $\Ghat$. The counit $\cO \lra k$, is viewed as
evaluation of functions at the identity $e \in \Ghat$, and the augmentation
ideal $I$ consists of functions vanishing at $e$. The element $y$ generates
the ideal $I$, and is known as a {\em coordinate} at $e$.

%We consider a group object $\G$ in a category of affine
%schemes over a ring $k$. We may consider the ring of functions
%$\cO$ on $\G$. Since $\G$ is a group, $\cO$ is a Hopf algebra.
%The identity element
%$e \in \G$ defines an ideal $I$ of functions vanishing at $e$.
We also need to discuss locally free sheaves $\cF$ over $\Ghat$, and in the present
affine context these are
specified by the $\cO$-module $M=\Gamma \cF$ of global sections. In particular,
line bundles $L$ over $\Ghat$ correspond to  modules $M$ which are
submodules of the ring of rational functions and free  of rank 1.
Line bundles can also be described in terms of
the zeroes and poles of their generating section: we only need this
in special cases made explicit below.
The generator $f$ of the $\cO$-module $M$ is
a section of $L$, and as such it defines a divisor $D=D_+-D_-$, where
$D_+$ is the subscheme of $\Ghat$ where $f$ vanishes (with multiplicities),
and $D_-$ is the subscheme of $\Ghat$ where $f$ has poles
(with multiplicites). This divisor determines $L$, and  we  write
$L=\cO(-D).$ For example, 
$$M=I=(y) \mbox{ corresponds to }\cO (-(e)), $$
and 
$$M=I^a=(y^a) \mbox{ corresponds to } \cO (-a(e)) .$$
Next we may consider the $[n]$-series map $[n]: \cO \lra \cO$, which
corresponds to the $n$-fold sum map
$n: \Ghat \lra \Ghat$. We write $\Ghat [n]$ for the kernel of $n$,
and its ring of functions is $\cO/([n](y))$. Hence, since
$n^*y=[n](y)$ by definition,
$M=([n](y))$ corresponds to $\cO (-\Ghat [n])$,  and
$M=(\; ([n](y))^a\; )$ corresponds to $\cO (-a\Ghat [n])$.
Finally, if $M$ corresponds to $\cO (-D)$
and   $M'$ corresponds to $\cO (-D')$ then
$M^{\vee}:=\Hom (M, \cO)$ corresponds to $\cO (D)$ and
$M \otimes M'$ corresponds to $\cO (-D-D')$. This gives sense to
enough line bundles for our purposes.

%In this section $\G$ will be a one dimensional
%commutative formal group, so that we may choose  a coordinate $y$ at $e$
%(i.e. there is an element $y\in \cO$ so that  $I=(y)$) and $\cO =k[[y]]$.

\subsection{Complex oriented cohomology theories.}
Now  suppose that $E$ is a 2-periodic ring valued  theory with coefficients
$E^*$ concentrated in even degrees. The collapse of the Atiyah-Hirzebruch
spectral sequence for $B\T $ shows that  $E$ is  complex orientable.
We may  define the \T -equivariant Borel
cohomology  by $\ETc (X)=E^*(\ET \times_{\T} X)$.
We work over the ring  $k=\ETz (\T)=E^0$, and
 view $\ETz =E^0(B\T  )$ as the ring of functions on a formal group
$\Ghat$ over $k$. The tensor product and duality
of line bundles makes $B\T $ into a group object, so
$E^0(B\T )$ is a topological Hopf algebra and $\Ghat$ is a group.
From this point of view, the augmentation ideal
$I=\ker (\ETz \lra E^0) $ consists of functions vanishing at the
 identity $e \in \Ghat$. We may also define
the module of cotangent vectors at the identity by 
$$\omega:=I/I^2=\tE^0(S^2)=E^{-2}=E_2. $$
This allows us to recover the graded cohomology ring from the 
ungraded ring since 
$$E_{\T}^{-2n}(X)=\ETz (X) \tensor \omega^n.$$

Now, if $W$ is a complex representation of the circle group $\T$
with $W^{\T}=0$,
we also let $W$ denote the  associated bundle over $B\T $ and
 the Thom isomorphism shows
% that $\tEz (X^{\eta})$ is a rank 1 free module over $E^0(X)$. In particular
$\tE^0 ((B\T)^W)=\tETz (S^W)$ is a rank 1 free module over $\ETz$,
and hence corresponds
to a line bundle $\L (W)$ over $\Ghat$, whose global sections are naturally
isomorphic to the module
$$\Gamma \L (W)=\tETz (S^W) .$$
From the fact that Thom isomorphisms are transitive we see that
$\L (W \oplus W')=\L (W) \otimes \L (W')$. The values of all these line
bundles can be deduced from those of powers of $z$.

\begin{lemma}
\be
\item $\L (0) =\cO $ is the trivial bundle.
\item $\L (z) =\cO (-(e))$ is the sheaf of functions vanishing at $e$, and
its module of sections $I$ is generated by the coordinate $y$.
\item $\L (z^n)= \cO (- \Ghat [n])$ is the sheaf of functions vanishing
on $\Ghat [n]$, and its module of sections is generated by the multiple
$[n](y)$ of the coordinate $y$.
\item $\L (az^n)= \cO (- a\Ghat [n])$ is the sheaf of functions vanishing
on $\Ghat [n]$ with multiplicity $a$, and its module of sections is generated
$([n](y))^a$.
\ee
\end{lemma}
\ppf\ The first statement is clear since $\tETz (S^0) =\ETz$. For the
second we use the equivalence $(B\T)^z\simeq (B\T)^0 /(\mathrm{pt})^0$.
The third statement follows from the Gysin sequence since
$z^k$ is the pullback of $z$ along the $k$th power map $B\T  \lra B\T$.
The final statement follows from the tensor product property.
\qqed

This gives the fundamental connection between the equivariant cohomology of
a sphere and sections of a line bundle.

\begin{cor} If $E_{\T}^*(\cdot )$ is a complex oriented 2-periodic
cohomology theory with associated formal group $\Ghat$ then 
for any $a \in \Z$, $n \neq 0$ we have
$$\tETz (S^{az^n})= \cO (-a\Ghat [n]).\qqed$$
\end{cor}

\section{What to expect when the group is not affine.}
\label{sec:nonaffine}

This  section  discusses what happens if we replace the formal group
$\Ghat$ (which is affine) in Section \ref{secFormal} by 
a (one dimensional) group $\G$ with higher cohomology. 

\subsection{Odd cohomology.} The main point is that 
we cannot expect a cohomology theory entirely in even degrees. 
Now that the group is not affine, $\cO$ denotes  the structure {\em sheaf}
of $\G$. 
This is reconciled to the above usage since in the affine case,  the structure
sheaf is determined by its ring of global sections.
In the non-affine case,  the cofibre sequence 
$$S^{az} \sm \T_+ \lra  S^{az} \lra S^{(a+1)z} $$
of based $\T$-spaces forces there to be odd cohomology. Indeed, we expect  a 
corresponding short exact sequence of sheaves
$$\cO (-ae)/\cO (-(a+1)e) \lla \cO (-ae) \lla \cO (-(a+1)e). $$
Any satisfactory cohomology theory will be functorial, and
applying $\tETz(\cdot )$ will give sections of the associated
sheaves. However the global sections functor on sheaves is not usually right
exact, and the sequence of sections continues with the sheaf
cohomology groups $H^1(\G ; \cdot )$.
It is natural to hope that the long exact cohomology sequence
 induced by  the sequence of spaces should be
 the long exact cohomology sequence induced by the sequence of
sheaves. This gives  a natural candidate for the odd cohomology:
$$\tET^i(S^{az})=H^i(\G ; \cO (-a(e))) \mbox{ for } i=0,1. $$
This explains why it is possible for complex orientable cohomology
theories to have coefficient rings in even degrees (formal groups are
affine), and how their values on all complex spheres can be the same
(formal groups have a regular coordinate). It also explains
why we cannot expect either property for a theory associated to
an elliptic curve. 

\subsection{The definition of type.} We are now ready to formalize
the relationship between  group schemes and cohomology theories. 

\begin{defn}
\label{typeG}
(i) Given a virtual complex representation  $W$ with $W^{\T}=0$ we define
an associated divisor $D(W)$ as follows. We write 
$W=\sum_na_nz^n$, and then take   $D(W)=\sum_na_n\G [n]$, where
$\G [n]=\ker (n: \G \lra \G)$.\\
(ii) We say that a 2-periodic  $\T$-equivariant cohomology theory
$E^*_{\T}(\cdot )$ is of {\em type $\G$} if, for 
any complex representation $W$, 
$$\widetilde{E}_{\T}^i(S^W)\cong H^i(\G ; \cO (-D(W))$$
and 
$$\widetilde{E}^{\T}_{-i}(S^W)\cong H^i(\G ; \cO (D(W)).$$
for $i=0,1$. 

We also require these  isomorphisms to be natural for inclusions $j:W \lra W'$ 
of representations. To describe this, first note that 
such a map induces a map $S^W \lra S^{W'}$ of based $\T$-spaces 
and hence  maps 
$$j^*: \widetilde{E}_{\T}^i(S^{W'}) \lra \widetilde{E}_{\T}^i(S^W)$$
and 
$$j_*:\widetilde{E}^{\T}_{-i}(S^{W}) \lra \widetilde{E}^{\T}_{-i}(S^{W'}).$$
On the other hand we have inclusion of divisors $D(W) \lra D(W')$, 
inducing maps
$$\cO (-D(W')) \lra \cO (-D(W))$$ 
and
$$\cO (D(W)) \lra \cO (D(W')). $$ 
The induced maps in sheaf cohomology are required to be $j^*$ and $j_*$.
\end{defn}

\begin{remark}
The naturality requirement really allows us to  identify the homology and
cohomology of spheres with spaces of functions or their duals. 
For example, all the sheaves $\cO (-D(V))$ are subsheaves of the constant sheaf 
$$\cK =\{ f \st f \mbox{ is a function on $\G$ with poles only at points of finite order }\}, $$
of meromorphic functions. Thus  the naturality requirement shows we may actually identify 
$\widetilde{E}_{\T}^0(S^{-W})$ with a set of meromorphic functions. 
In the presence of Serre duality (see Section \ref{sec:mult}), the first cohomology groups may 
similarly be identified with duals of spaces of functions.
\end{remark}

\begin{remark}
We also need to discuss the appropriate behaviour for 
for representations $W$ with trivial summands. The convention that 
$$E^{\T}_{i+2}(X)=E^{T}_i(X) \tensor \omega$$
or 
$$E_{\T}^{i-2}(X)=E_{T}^i(X) \tensor \omega $$
leads to an appropriate formula, where $\omega$ is the cotangent space at the identity 
of $\G$. However to obtain a properly natural identification it is better to use
sheaves and the fact that 
$$H^i(\G ; \mcL \tensor_{\cO} \Omega )=H^i(\G ; \mcL) \tensor  \omega , $$
where $\tensor $ denotes tensor product over $k$, and $\Omega$ is the sheaf of K\"ahler 
differentials on $\G$.

This leads to the requirements
$$\widetilde{E}_{\T}^i(S^W)\cong H^i(\G ; \cO (-D(W/W^{\T})) \tensor_{\cO} \Omega^{\dimC (W^{\T})})$$
and 
$$\widetilde{E}^{\T}_{-i}(S^W)\cong H^i(\G ; \cO (D(W/W^{\T}))\tensor_{\cO} \Omega^{-\dimC (W^{\T})})$$
for $i=0,1$ (here and elsewhere $\Omega^n$ denotes the $n$th tensor power
of $\Omega$). The answer for other values of $i$ follows easily.
\end{remark}

\begin{remark}
The use of differentials to give suspensions  means that a cohomology theory of type $\G$
contains data about Thom isomorphisms. For example, if 
$S^{\infty W}:=\colim_{U^{\T}=0}S^U$ then $S^z \lra S^{\infty W}$ induces
$$\begin{array}{ccc}
\wtET_2(S^z )& \lra & \wtET_2(S^{\infty W} )\\
\cong \downarrow && \downarrow \cong\\
H^0(\G ; \cO ((e)) \tensor_{\cO} \Omega ) & \lra & H^0(\G ; \cK \tensor_{\cO}\Omega).
\end{array}$$
This picks out a $k$-subspace of the constant  sheaf
$\cK \tensor_{\cO }\Omega= H^0(\G; \cK \tensor_{\cO }\Omega) $.
When $\G$ is an elliptic curve, this is the one dimensional space of 
invariant differentials.
\end{remark}

\subsection{The affine case revisited.}
\label{subsec:affine}
It is worth pointing out that if $\G$ is affine and has a good coordinate, 
any cohomology theory  of type $\G$ is complex orientable and in even degrees
(we construct a number of such theories in Appendix \ref{sec:GaGm}).
More precisely, we require that $\G$ has a regular coordinate
function  $y$ in the sense that the identity $e \in \G$ is defined
by the vanishing of $y$ and $y$ is a regular element of the ring 
$\cO$ of functions on $\G$. The multiplication by $n$ map is 
also required to be flat for $n \geq 1$.

First, since $\G$ is affine, there is no higher cohomology. Thus, 
the condition that $E^*_{\T}(\cdot)$ is of type $\G$ states
that the cohomology of spheres of complex representations is in 
even degrees, and that if $W^{\T}=0$,   
$$\widetilde{E}_{\T}^{-2n}(S^W)=\cO (-D(W)) \tensor \omega^{n}, $$
where we have identified the sheaf with its space of global
sections. It remains to observe that $\cO (-D(W)) $ is a 
free module on one generator. 
Indeed, $\G [n]$ is defined by the vanishing of $n^*(y)$
the pullback of $y$ along the multiplication by $n$  map of 
$\G$. Since this map is flat, $n^*y$ is a regular element.

Since we have a complex oriented theory we also have
Thom classes and Euler classes, and these  depend on the coordinate, $y$.
For example, the Thom class of $z^n$ is the chosen generator 
of $\cO (-\G[n])$, and the Euler class is its pullback to 
$\cO$, namely 
$$\chi_y (z^n)=[n](y):=n^*(y).$$
Thus we have the idea that the Euler class of $z^n$ is a function whose
vanishing defines $\G [n]$. 

In characteristic 0 it is elementary to go one step further and decompose the 
divisor $\G [n]$:
$$\G [n]=\sum_{s|n}\Gpointy{s}$$
where $\Gpointy{s}$ is the divisor of points of exact order $s$.
In fact, we define a function  $\philr{s} (y)$
vanishing to the first order on $\Gpointy{s}$ recursively by the condition
$$\chi_y (z^n) =\prod_{s|n} \philr{s} (y):$$
the formula for $n=1$ defines $\philr{n}(y)$ directly for $n=1$, 
and for larger values of $n$, $\philr{n}(y)$ is defined by dividing 
$\chi_y (z^n)$ by the previously defined $\philr{s}(y)$. Each $\phi_n(y)$
is regular by the regularity of $y$ and the flatness requirement.

\subsection{Summary.}
We may summarize the correspondence between algebra and topology.

\begin{itemize}
\item The suspension $S^{az^n}\sm \EG$ corresponds to the sheaf $\cO (a\G[n])$
and more generally, suspension by $z^n$ corresponds
to tensoring with $\cO (\G[n])$.
\item The subgroup $\T [n]$ of order $n$ (kernel of $z^n$) corresponds
to the subgroup $\G [n]$ of elements of order dividing $n$ (defined by 
the vanishing of $\chi (z^n)$).
\item The inclusion $S^0 \lra S^{z^n}$ which induces multiplication by the
Euler class (in the presence of a Thom isomorphism)
corresponds to $\cO \lra \cO (\G[n])$.
\item We extend the notation, so that 
$$S^{\infty z^n}:=\colim_aS^{az^n}
\mbox{ corresponds to the sheaf }
\cO (\infty \G[n]):=\colim_a\cO (a\G[n])$$
and
$$\etf :=\colim_{U^{\T}=0} S^U
\mbox{ corresponds to the sheaf }
\cOiGt :=\colim_{a,n}\cO (a\G[n]).$$
\item The family $\cF$ of finite subgroups corresponds to the set $\G [tors]$
of elements of torsion points.
\end{itemize}

\part{Background on rational 
$\protect\T$-equivariant cohomology theories.}

The method of this paper is only practical because there is a complete
algebraic model for rational $\T$-equivariant cohomology theories \cite{s1q}.
In Part 2 we describe this model and explain how to make relevant 
calculations in it.

\section{The model for rational $\protect \T$-spectra.}
\label{sec:ModelT}

For most of the paper we work with the representing objects of
$\T$-equivariant cohomology theories, namely  $\T$-spectra \cite{LMS(M)}.
Thus we prove results about the representing spectra,
and deduce consequences about the cohomology theories.
More precisely,  any suitable
$\T$-equivariant cohomology theory $E_{\T}^*(\cdot)$
is  represented by a $\T$-spectrum $E$ in the sense that for a based
$\T$-space $X$, 
$$\tET^*(X)=[X,E]_{\T}^*.$$
This enables us to define the associated homology theory
$$\widetilde{E}^{\T}_*(X)=[S^0,E \sm X]^{\T}_*$$
in the usual way. We shall make use of the elementary
fact that the Spanier-Whitehead dual of the sphere $S^W$ is $S^{-W}$, as one
sees by embedding $S^W$ as the equator of  $S^{W\oplus 1}$.
Hence, for example
$$\tETz (S^W)=[S^W,E]^{\T}=[S^0,S^{-W}\sm E]^{\T}=
\pi^{\T}_0(S^{-W} \sm E)=\tilde{E}^{\T}_0(S^{-W}). $$

We say that a cohomology theory is {\em rational} if its values
are graded rational vector spaces. A  spectrum is rational if
the cohomology theory it represents is rational. It suffices to
check the values on the homogeneous spaces $\T/H$ for closed subgroups
$H$, since all spaces are built from these up to weak equivalence.

\begin{convention}
Henceforth  all spectra and the values of all cohomology theories
are rationalized whether or not this is indicated in the notation.
\end{convention}

Our results are made possible because there is a complete algebraic
model of the category of {\em rational} $\T$-spectra, and hence of
rational $\T$-equivariant cohomology theories \cite{s1q}. For the
convenience of the reader we spend the rest of this section summarizing
 the relevant results from  \cite{s1q} in a convenient form. 
There are two models for rational \T -spectra, as derived categories of
abelian categories. 
\begin{thm}\cite[5.6.1, 6.5.1]{s1q}
There are equivalences
$$\T\mbox{-Spectra} \simeq \das \simeq \dat. $$
of triangulated categories. 
\end{thm}
The {\em standard} abelian category  $\as$ has injective dimension 1,
and the {\em torsion} abelian category $\at$  is of injective  dimension 2. 
The derived category $\das$ is formed by taking differential graded objects
in $\as$ and inverting homology isomorphisms, and similarly for $\dat$.
It is usually easiest to identify the model for a \T -spectrum in \dat , 
at least providing its model has homology
of injective dimension 1. This is then  transported to the
standard category, where calculations are sometimes easier.
We describe what we need about the categories in the following
subsections. 

\subsection{Rings of functions.}
To describe the categories, we need some ingredients. The information 
is organized by isotropy group, and we let  
$\cF$ denote the discrete set of finite subgroups of $\T$. 
On this we consider the constant
sheaf $\mcR$ of rings with stalks $\Q [c]$ where $c$ has degree $-2$. We need
to consider the ring 
$$R=\map (\cF , \Q [c])\cong \prod_{H \in \cF}\Q [c]$$ 
of global sections, where maps and product are  graded. 
For each subgroup $H$, we let $e_H \in R$ 
denote the idempotent with support $H$.

To avoid confusion about  grading we introduce the requisite suspensions.
In topology we may suspend by   complex representations $W$; these
enter the theory through  the dimension function $w(H):=\dimC (W^H)$. 
Notice that $w$ takes only finitely many values, and is equal to 
$w(\T)$ for almost all finite subgroups $H$.

\begin{defn}
\label{defn:suspRmod}
Suppose $w: \cF \lra \Z$ is an almost constant  
function. We  divide the set of finite subgroups into sets 
$$\cF_{w,i}=\{ H \st w(H)=i\}$$
on which $w$ is constant; only finitely many of these are non-empty, 
and  all but one are finite. We write $w(\T )$ for
the value of $w$ on the infinite set.

Let $e_{w,i}\in R$ be the idempotent supported on $\cF_{w,i}$, and
introduce the suspension functor on $R$-modules by 
$$\Sigma^wN=\bigoplus_i \Sigma^{2i} e_{w,i}N.$$
\end{defn}

Now if  $w: \cF \lra \Z_{\geq 0}$ is zero almost everywhere, we write 
$c^w $ for the {\em universal Thom class} of $w$, defined by  
$c^w(H)=c^{w(H)}$. Since it is not homogeneous, $c^w$ is not an  
element of $R$, but nonetheless it is natural to consider the  $R$-module
$$c^w R :=\Sigma^{-w}R=\prod_H c^{w(H)} \Q [c],  $$
viewed as an  $R$-submodule of $\prod_H \Q [c,c^{-1}]$; since $c^w$ is
a generator in some sense,  we call $c^w$ a Thom class (further explanation 
is given at the end of the section).  Classical 
Thom classes give rise to Euler classes by restriction to the coefficient ring.
We now create a ring in which  the Euler classes
corresponding to the Thom classes $c^w$ belong. First, let
$$\cE =\{c^w \st w:\cF \lra \Z_{\geq 0} \mbox{ of finite support}\}; $$
thinking of this as if it generates a multiplicatively closed subset, 
we make an adelic construction by  forming the $R$-submodule
$$\tf =\cEi R := \colim_w \Sigma^w R = \bigcup_w c^{-w}R $$
of $\prod_H \Q [c,c^{-1}]$.   
Observe that $\tf$ is a graded $R$-algebra. As a graded vector space $\tf$ is
$\bigoplus_H \Q$ in positive degrees and $\prod_H \Q$ in degrees
zero and below.

\begin{remark}
(i) Note that if $w(\T)=0$, there is a natural degree 0 {\em isomorphism}
$$c^w:\tf \stackrel{\cong} \lra \Sigma^w \tf , $$
which in the $s$th component is 
$$c^{w(s)}: \Q [c,c^{-1}] \lra \Sigma^{2w(s)}\Q [c,c^{-1}]. $$
It is natural to see this as multiplication by the Euler class.\\ 
(ii) Given a complex representation $W$ of $\T$ with $W^{\T}=0$
we may define an associated function $w: \cF \lra \Z_{\geq 0}$ zero
almost everywhere by $w(H)=\dimC (W^H)$. We sometimes write
$c^W$ for this element of $R$, and we note that $\cE$ is generated
as a multiplicative subset by elements of this form.\\ 
(iii) Viewing $R$ as the ring of functions on the discrete space $\cF$, 
 the universal Euler classes can be used to define finite subsets. 
Indeed, we may view $c^w$ as a non-homogeneous section of the structure
sheaf, or as a homogeneous section of a line bundle 
$\mcR (-w)$ with global sections $\Sigma^{-w}R$. 
Now, for any finite subset $\mcH \subseteq \cF$ we may consider its 
characteristic function $\chi (\mcH)$.  The
associated universal Euler class $c^{\chi (\mcH)}$ is the function 
vanishing to the first order on $\mcH$.
\end{remark}

%$e(V)$ for the representations $V$ of $\T$ with $V^{\T}=0$.
%These are defined by $e(V)=c^v$, where $v(H)=\dimC(V^H)$.
%In particular for  $V=z^n$
%we have $e(z^n)=c^{\mathrm{sub} (n)}$ where $\mathrm{sub} (n)(H)=1$
%if $H \subseteq \T [n]$ and $0$ otherwise.
%Equivalently,

\subsection{Description of the abelian categories.}
The objects of the standard model
\as\ are triples $(N,\beta , V)$ where $N$ is an $R$-module (called the
{\em nub}), $V$ is a graded rational vector space (called the {\em vertex})
and  $\beta: N \lra \tf \otimes V$ is a morphism
of $R$-modules (called the {\em basing map})
 which becomes an isomorphism when $\cE$ is inverted.
When no confusion is likely, we simply say that $N \lra \tf \otimes V$
is an object of the standard abelian category.
An object of $\cA_s$ should be viewed as the module $N$ with the  additional
structure of a trivialization of $\cEi N$.
A morphism $(N, \beta , V) \lra (N', \beta', V')$ of objects is given
by an $R$-map $\theta : N \lra N'$ and a $\Q$-map $\phi : V \lra V'$
compatible under the basing maps.

 Since the standard abelian category  has injective dimension 1,
homotopy types of objects of the derived category $D(\cA_s)$
are classified by their homology in $\cA_s$, so that
homotopy types correspond to isomorphism classes of  objects
of the abelian category $\cA_s$. In the sheaf theoretic approach
\cite{atnq1}, $N$ is the space of global sections of a sheaf on the space of
closed subgroups $\T$, the vertex $V$ is the value of the sheaf at
the subgroup $\T$ and the fact that the basing map
$\beta : N \lra \tf \otimes V$ is an isomorphism away from
$\cE$ is the manifestation of the patching condition for sheaves.

The objects of the torsion abelian category
\at\ are triples $(V,q,T)$ where $V$ is a graded rational vector space,
$T$ is an $\cE$-torsion $R$-module  and
$q: \tf \otimes V \lra T$ is a morphism of $R$-modules.
The condition on $T$ is equivalent to requiring (i) that
$T$ is the sum of its idempotent
factors $T(H)=e_HT$ in the sense that $T=\bigoplus_H T(H)$
and (ii) that  each $T(H)$ is a torsion $\Q [c]$-module.
When no confusion is likely, we simply say that $\tf \otimes V \lra T$
is an object of the torsion abelian category.
In the sheaf theoretic approach,  the module $T(H)$ is the
cohomology of the structure sheaf with support at $H$.
By contrast with the standard abelian category, the torsion abelian category
has injective dimension 2. Thus not every object $X$ of the
derived category $D(\cA_t)$ is determined up to equivalence
by its homology $H_*(X)$ in the abelian category $\cA_t$.
We say that $X$ is {\em formal} if it is equivalent to its homology
(considered as a differential graded object with zero differential), and
that it is {\em intrinsically formal} if it is equivalent
to any object with the same homology. Evidently, an intrinsically formal object
is formal. The Adams spectral sequence
shows immediately  that  $X$ is intrinsically 
formal if its homology has injective dimension 0 or
1 in $\cA_t$. In general, if  $H_*(X)=(\tf \otimes V \lra T)$,
the object $X$ is equivalent to  the fibre of a
 map $(\tf \otimes V \lra 0) \lra (\tf \otimes 0 \lra \Sigma T)$
(in the derived category) between objects in $\cA_t$
 of injective dimension 1.  This map is classified by an element of
$\Ext_R (\tf \otimes V ,\Sigma T)$, so that $X$ is formal
if the Ext group is zero in even degrees. Thus $X$ is 
intrinsically formal
if both $V$ and $T$ are in even degrees or (\cite[5.3.1]{s1q})
if  $T$ is injective
in the sense that each $T(H)$ is an injective $\Q [c]$-module.

\begin{lemma}
\label{mapsq}
An  $R$-map $q: \tf \otimes V \lra \bigoplus_s T_s$ 
is determined by its idempotent pieces 
$q_s : \Q [c,c^{-1} ] \otimes V \lra T_s$. 

Conversely, any sequence  of $\Q [c]$-maps $q_s$ 
so that, for each $f \in V$, only finitely 
many of the values $q_s(c^0 \otimes f)$ are  non-zero, 
 determines an $R$-map $q$. 
\end{lemma}

\begin{proof}
To see that the idempotent pieces determine $q$, note that if all idempotent 
pieces are zero we may argue that $q=0$: if $q(1 \tensor v) \neq 0$
some idempotent piece would be non-zero, hence $q$ vanishes on $R \tensor V$, 
and hence induces a map 
$$\overline{q} : (\tf /R) \tensor V =
\bigoplus_s (\Q[c,c^{-1}]/\Q [c] )\tensor V \lra \bigoplus_s T_s, $$
which is the direct sum of its idempotent pieces. 

The converse statement is easily checked.
\end{proof}

\subsection{Spheres, suspensions and Euler classes.}
Spheres are important because they are invertible objects, and
 therefore play a role corresponding to that of 
line bundles in categories of sheaves. 
We introduce the appropriate apparatus to discuss them.

We described the suspension $\Sigma^w$ on the category of 
$R$-modules in \ref{defn:suspRmod}. 

\begin{defn}
\label{defn:suspAs}
The suspension functor on objects of the standard abelian 
category $\cA_s$ is defined by 
$$\Sigma^w(N \lra \tf \tensor V)=
(\Sigma^wN \lra \Sigma^w \tf \tensor V \stackrel{c^w}\lra
\tf \tensor \Sigma^{2w(\T )}V).$$
Thus, the  basing map for the suspension is obtained by multiplying 
the original one by the appropriate Euler class, which is 
$c^{w(i)-w(\T)}$ on $e_{w,i}N$. 
\end{defn}

\begin{defn}
\label{spheres} \cite[5.8.2]{s1q}
The algebraic 0-sphere is the object
$$S^0=(R \lra \tf \tensor \Q)$$
where $R$ is the submodule of $\tf \tensor \Q$ generated
by $1 \tensor 1$.

Given an almost constant function $w: \cF \lra \Z$ the algebraic 
$w$-sphere is the object of $\cA_s$ defined by
$$S^w=\Sigma^wS^0 =(R (w) \lra \tf \tensor \Sigma^{2w(\T)}\Q) $$
where 
$$R(w)=\Sigma^w R =c^{-w}R \subseteq 
\Sigma^w \tf \stackrel{\cong} \lra \Sigma^{2w(\T)}\tf$$
as above. Note that different parts of this diagram have been 
shifted by different amounts, so that both the grading and the
structure maps are different for different spheres. 
\end{defn}

If $X$ is a $\T$-spectrum we write $M_s(X)$ and $M_t(X)$ for the models
of $X$ in $\As$ and $\At$ respectively. In fact,  if
$\Phi^{\T}X$ denotes the geometric fixed point spectrum of $X$, and 
 $E\cF$ denotes  the universal almost free $\T$-space, we have 
$$H_*(M_t(X))=(\tf \tensor V \lra T), $$
where
$$V=\pi_*(\Phi^{\T}X)), $$
and 
$$T=\pi^{\T}_*(\Sigma E\cF_+ \sm X).$$
 Since $\At$ is 
of injective dimension 2, this does not always determine $M_t(X)$.
On the other hand, 
since $\As$ is of injective dimension
$1$, we may take $M_s(X)$ to be an object of the underlying abelian 
category $\As$ (i.e., to have zero differential). In fact,  
$$M_s(X) \simeq H_*(M_s(X))=(N \lra \tf \tensor V)$$
where $V$ is as above and $N$ lives in a long exact sequence
$$\cdots \lra N \lra \tf \tensor V \lra T \lra \cdots $$
This at least makes clear that $V$ is to do with $\T$-fixed points of $X$, 
$T$ is to do with the almost free part of $X$ and $N$ is an appropriate 
mixture. It also suggests the relationship between $\As$ and $\At$.
This amount of detail is more than we need for the present
paper. Finally,  we need to record that spheres  in the
algebraic and topological contexts correspond.

\begin{lemma}
\label{suspensions}
 \cite[5.8.3]{s1q}
Suppose $W$ is a virtual complex representation, and let 
$w=\mathrm{dim}_{\C}(W)$.\\
(i) The object modelling the sphere $S^W$ in $\cA_s$ is
the algebraic sphere $S^w$:
$$M_s(S^W)=S^w=(R(w) \lra \tf).$$
(ii) Algebraic and topological suspensions coincide 
in the sense that 
$$M_s(\Sigma^W X)=\Sigma^w M_s(X).$$
\end{lemma}

\ppf\ Part (ii) follows from Part (i) since the algebraic suspension
is tensor product with $S^w$ and $S^w$ is flat.
\qqed

\begin{warning} We are modelling 
{\em complex} representations $W$. Thus if $\eps$ is the trivial 
representation of $\T$ on $\C$, we have $S^{\eps}=S^2$. We thus need
to be careful when discussing a single suspension 
(smash product with the circle). 
We use the same method to resolve this conflict in algebra as
in topology: an integer has its usual meaning, whereas the function 
$\cF \lra \Z$ with constant value 1 will be denoted $\eps$. 
\end{warning}

We are now in a position to justify calling the function 
$c^w$ a universal Euler class when $w(\T)=0$. 
In the topological context,  the Euler class of a complex representation 
$W$ with $W^{\T}=0$ in a complex oriented cohomology theory is
 defined by pulling back a Thom class along $S^0 \lra S^W$; 
equivalently, in the associated homology theory we take the
image of the Thom class under the map $S^{-W} \lra S^0$.
In the algebraic context we do precisely the same. The Thom
class of $S^{-W}$ is the `generator' of $R(-w)$, namely the 
`element' $c^w$, which is the image of $1 \in \tf$ under
the isomophism $c^w :\tf \lra \tf$. The two obstructions to a 
universal Thom isomorphism are the two linked facts that $c^w$ is 
not homogeneous and that the putative isomorphism is not compatible
with basing maps.

Consider the subgroup $\T [n]$  of order $n$, and the representation 
$z^n$. If we take the $K$-theory Euler class we have
$$e(z^n)=1-z^n =\prod_{s|n}\phi_s, $$
where $\phi_s$ is the $s$th cyclotomic function, independent of $n$.
Similarly, the dimension function corresponding to $z^n$, is 
the characteristic function $\mathrm{sub}(n)$ for the subgroups of
$\T [n]$. Hence the universal Euler class defining the subgroups of $\T [n]$ is
$$c^{z^n}=c^{\mathrm{sub}(n)}=\prod_{s|n}c_s , $$
where $c_s$ is the universal Euler class for the characteristic
function of the singleton $\{ \T [s] \}$. It is therefore natural  
to view $c_s$ as a universal cyclotomic function.

\section{Cohomology of spheres.}
The main point of contact between topology and geometry is
through the cohomology of spheres and line bundles. We therefore
describe how this works in the standard model for $\T$-spectra. 
We shall only need to discuss $\T$-spectra with  particularly nice 
algebraic models, so we begin by describing them.

\subsection{Rigidity.}

Given a $\T$-spectrum $E$ with torsion model 
$M_t(E)$ with homology $H_*(M_t(E))=(\tf \tensor V \lra T)$
in the abelian category $\cA_t$, it is not hard to calculate
$V= E^{\T}_*(\etf )$ or $T=E_T^*(\Sigma^{-1}\efp)$. However, 
if this is to determine $E$ we must show in addition
that $M_t(E)$ is formal.

\begin{defn}
We say that a $\T$-equivariant cohomology theory $E$ is 
{\em rigid} if the following two equivalent conditions
hold 
\be
\item $H_*M_t(E)=(\tf \otimes V \stackrel{q}\lra T)$
has  surjective structure map $q$.
\item $H_*M_s(E)=(N \stackrel{\beta}\lra \tf \tensor V)$ has injective
structure map $\beta$.
\ee
We say that a rigid spectrum  $E$ is {\em even} if $V$, $T$ 
and $N$ are concentrated in even degrees. 
\end{defn}

\begin{lemma}
\label{surjMtgivesMs}
If $E$ is rigid then $M_t(E)$ is intrinsically formal, and
if $H_*M_t(E)=(\tf \otimes V \stackrel{q}\lra T)$
then 
$$M_t(E)\simeq (\tf \otimes V \stackrel{q}\lra T)$$
and 
$$M_s(E)\simeq (N \stackrel{\beta}\lra \tf \otimes V)$$
where
$$N=\ker (\tf \otimes V \lra T),  $$
and the basing map $\beta$ is the inclusion.
Furthermore we have the explicit injective resolution
$$0\lra M_s(E )
\simeq \left(
\begin{array}{c}
N\\
\downarrow\\
\tf \otimes V
\end{array}
\right)
\lra
\left(
\begin{array}{c}
\tf \otimes V\\
\downarrow\\
\tf \otimes V
\end{array}
\right)
\lra
\left(
\begin{array}{c}
T\\
\downarrow\\
0
\end{array}
\right)
\lra 0 $$
in $\cA_s$.
\end{lemma}

\ppf\ To see that $M_t(E)$ is formal, it is only necessary to remark
that $T$ is the quotient of
an $\cE$-divisible group and therefore injective \cite[5.3.1]{s1q}.
\qqed

\begin{lemma}
\label{rigidisflat}
If $E$ is rigid, the  corresponding object $M_s(E)=(N \lra \tf \tensor V)$ in 
$\cA_s$ is flat. 
\end{lemma}

\begin{proof} Tensor product on $\As$ is defined termwise.
First, note that  $\tf \otimes V$ is exact for tensor product
with objects $P$ with $\cEi P \cong \tf \otimes W$ for some $W$, 
so the tensor product is exact on the vertex part. 

For the nub, we use the fact that  
the category $\As$ is of flat dimension 1 by \cite[23.3.5]{s1q}, 
together with the fact that  $N$ is a submodule of $\tf \otimes V$.
\end{proof}

\subsection{Homomorphisms out of $S^0$.}
For an object $X$ of $\cA_s$ there is an exact sequence
$$0 \lra \Ext_{\cA_s} (S^{1+w} , M) \lra [S^w,M]
 \lra \Hom_{\cA_s} (S^{w} , M) \lra 0, $$
so we shall need to calculate these  Hom and Ext groups.
For the present we restrict ourselves to the Hom groups.
We avoid confusion about grading by restricting
to the case $w=0$ using $[S^w,M]=[S^0 , \Sigma^{-w}M]$.

\begin{lemma}
\label{htpygroups}
For an object $M=( N \stackrel{\beta}\lra \tf \otimes V)$
of the abelian category $\cA_s$
$$\Hom_{\cA_s}(S^0,(N \lra \tf \otimes V)) = N(c^{0}):=
\{ n \in N \st \beta (n) \in c^{0} \otimes V \}.$$
\end{lemma}

\begin{proof}
A homomorphism $f: S^0 \lra M$ of degree 0 is given by a square
$$\begin{array}{ccc}
R & \stackrel{\theta} \lra  & N\\
\downarrow && \downarrow\\
\tf \tensor \Q & \stackrel{1 \tensor \phi } \lra & \tf \tensor V.
\end{array}$$
Thus $f$ is determined by the $R$-map $\theta$, and 
  $\Hom_R (R,N)=N$. On the other hand, 
the image of $1 \in R$ under the basing map is $1 \tensor 1$, 
which imposes the stated condition, since $\phi (1) \in 
V_{0}$.
\end{proof}

% \begin{lemma}
% \label{htpygroups}
% For an object $M=( N \stackrel{\beta}\lra \tf \otimes V)$
% of the abelian category $\cA_s$ and a function 
% $w: \cF \lra \Z$ taking finitely many values and 
% the value $w(\T)$ almost everywhere, we have
% $$\Hom_{\cA_s}(S^w,(N \lra \tf \otimes V))_0 = \Sigma^wN(c^{-w}):=
% \{ n \in \Sigma^w N \st \beta (n) \in c^{-w} \otimes V_{2w(\T)} \}.$$
% \end{lemma}

% \begin{proof}
% A homomorphism $f: S^w \lra M$ is given by a square
% $$\begin{array}{ccc}
% R(w) & \stackrel{\theta} \lra  & N\\
% \downarrow && \downarrow\\
% \tf \tensor \Sigma^{w(\T)}\Q & \stackrel{1 \tensor \phi } \lra & \tf \tensor V.
% \end{array}$$
% Thus $f$ is determined by the $R$-map $\theta$, and 
%   $\Hom_R (R(w),N)=\Sigma^{-w}N$. On the other hand, 
% the image of $c^{-w} \in R(w)$ under the basing map is $c^{-w} \tensor 
% \Sigma^{2w(\T)}1$, 
% which imposes the stated condition, since $\phi (\Sigma^{2w(\T)}1) \in 
% V_{2w(\T)}$.
% \end{proof}

\subsection{Cohomology of spheres.}

The aim of the present section is to make explicit
the calculation of $E^{\T}_*(S^W)$ in terms of $H_*(M_t(E))
=(q: \tf \tensor V \lra T) $ assuming that $E$ is rigid and
even.

\begin{lemma}
\label{suspofrigid}
Suppose $w: \cF \lra \Z$ is zero almost everywhere.
If $E$ is rigid and even then the $w$th suspension $\Sigma^wE$ is 
rigid and even.

If $M_t(E)=(\tf \otimes V \stackrel{q}\lra T)$
then 
$$M_t(\Sigma^w E)\simeq (\tf \otimes V \stackrel{q^w}\lra \Sigma^wT), $$
where the structure map is given by 
$$q^w(c_s^{i(s)} \tensor \alpha  )=q(c_s^{i(s)+w(s) }\tensor \alpha)
\in e_s(\Sigma^wT)_{2n-2i(s)}$$
 for $\alpha \in VA_{2n}$. Thus 
 
$$ M_s(\Sigma^w E) =(\Sigma^w N \lra \tf \tensor V)$$
where 
$$\Sigma^w N = \ker( \tf \tensor V \stackrel{q^w} \lra \Sigma^wT). \qqed$$
\end{lemma}

\begin{remark}
\label{cwconv}
A natural mnemonic is to write 
$$q(xc^w \tensor \alpha ) = q^w(x \tensor \alpha), $$
despite the fact that  $xc^w$ is not an element of $\tf$.
\end{remark}

We may now assemble the information to calculate the homology 
of spheres. 
\begin{cor}
\label{rigidcohom}
Suppose that $E$ is rigid and even, so that 
$H_*(M_t(E))=(q:\tf \tensor V \lra T)$ is surjective
and $V$ and $T$ are in even degrees. 
For any function $w: \cF \lra \Z$ zero almost everywhere
$$\widetilde{E}_0^{\T} (S^w)=
   \ker (q : c^w \otimes V_0 \lra (\Sigma^w T)_0)$$
and 
$$\widetilde{E}_{-1}^{\T} (S^w)=
    \cok (q : c^w \otimes V_0 \lra (\Sigma^w T)_0)$$
\end{cor}

\begin{proof}
To calculate the homology we use the short exact sequence
$$0 \lra \Ext_{\cA_s}(S^1, M_s(\Sigma^w E))
\lra E^{\T}_0(S^w) \lra \Hom_{\cA_s}(S^0, M_s(\Sigma^w E)) \lra 0. $$
We may calculate the Hom and Ext groups by applying 
$\Hom_{\cA_s}(S^0, \cdot )$ to the injective resolution 
of $\Sigma^wE$ given in \ref{surjMtgivesMs}.
\end{proof}

% \subsection{Grading again.}
% In the present paper we are interested in cohomology
% theories with a periodicity element $u$ of degree 2. 
% We may therefore shift even degree elements  into degree zero.
% For example $uc$ is the degree 0 counterpart of $c$.
% %For the rest of the paper we use $c$ to denote the degree 0 version.

\part{Background on elliptic curves.}
In Part 3 we summarize relevant facts about elliptic curves, and
make some easy deductions that we will need for the construction
of rational $\T$-equivariant elliptic cohomology. 

%%FIELDS
\section{Elliptic curves.}
\label{secElliptic}

In this section we record the well known facts about elliptic
curves that will play a part in our construction. We use
\cite{Silverman} as a  basic reference
for facts about elliptic curves, and \cite{Hartshorne} as background
from algebraic geometry.

Let $A$ be an elliptic curve
(i.e.,  a smooth projective curve of genus 1 with a specified point $e$)
over a field $\kk$ of characteristic 0
and let $\cO=\cO_A$ be its  sheaf of regular functions. Note that
$\Gamma \cO =\kk$, so the sheaf contains a great deal more information
than its ring of global sections. A divisor on $A$ is a finite $\Z$-linear
combination of points defined over the algebraic closure $\overline{\kk}$
of $\kk$, and associated to any rational function $f$ on $A$
we have the divisor $\div (f)=\Sigma_P \ord_P(f)(P)$, where
$\ord_P(f)\in \Z$ is the order of vanishing of $f$ at $P$. If a divisor
is fixed by $\mathrm{Gal}(\overline{\kk}/\kk)$ it is said to be defined over
$\kk$, and all the divisors we consider will be of this sort. 
In the usual way, if $D$ is a divisor on $A$,
we write $\cO (D)$ for the associated invertible sheaf.
Its global sections are given by
$$\Gamma \cO (D)=\{ f \; |\; \div (f)  \geq -D \} \cup \{ 0\} , $$
so that for a point $P$,  the global sections of $\cO (-P)$ are
the functions vanishing at $P$.

We also have $\cO (D_1) \otimes_{\cO} \cO (D_2)=\cO (D_1 +D_2).$

Since the global sections functor is not right exact, we are
led to consider cohomology, but since $A$ is one-dimensional
this only involves $H^0(A; \cdot )=\Gamma (\cdot )$ and
$H^1(A; \cdot )$, which are related by Serre duality. This
takes a particularly simple form since the canonical
divisor is zero on an elliptic curve:
$$H^0(A; \cO (D) )=H^1(A; \cO (-D) )^{\vee}, $$
where $(\cdot )^{\vee}=\Hom_k(\cdot ,k)$ denotes  vector space duality.

From the Riemann-Roch theorem we deduce that the canonical divisor is 0
and the cohomology of each line bundle:
$$ \dim (H^0(A; \cO (D))=\left\{
\begin{array}{ll}
\deg{D} & \mbox{ if } \deg (D) \geq 1\\
0       & \mbox{ if } \deg (D) \leq -1
\end{array}
\right.
$$
and
$$ \dim (H^1(A; \cO (D))=\left\{
\begin{array}{ll}
|\deg{D} |& \mbox{ if } \deg (D) \leq - 1\\
0       & \mbox{ if } \deg (D) \geq 1 .
\end{array}
\right.
$$
For the trivial divisor one has
$$\dim (H^0(A;\cO ))=\dim (H^1(A;\cO ))=1.$$
Now if  $D=\Sigma_P n_P(P)$ is a divisor of degree 0, we
may form the sum $S(D)=\Sigma_P n_P P$ in $A$, and
$D$ is linearly equivalent to $(S(D))-(e)$. If $S(D)=e$ 
then the sheaf $\cO (D)$ has the same cohomology as $\cO$.
Otherwise,  since no function vanishes to order
exactly 1 at $P$, we find
$$H^0(A;\cO (D))=H^1(A;\cO (D))=0. $$

We may recover $A$ from the graded ring
$\Gamma (\cO (*e))=\{ \Gamma \cO (ne)\}_{n \geq 0}$.
Indeed, this is the basis of the proof in \cite[III.3.1]{Silverman} that
any elliptic curve is a subvariety of $\bP^2$ defined  by a Weierstrass
equation. We choose a basis $\{ 1,x\}$ of $\Gamma \cO (2e)$ and a
extend it to a basis $\{ 1,x,y\}$ of $\Gamma \cO (3e)$. Now observe that
since $\Gamma \cO (6e)$ is 6-dimensional, there is a relation between
the seven 
elements $1,x,x^2,x^3,y,xy$ and $ y^2$: this is the Weierstrass equation,
and it may be verified that $A$ is the closure in $\bP^2$ of the plane
curve it defines.   The graded ring $\Gamma (\cO (*e))$
has generator $Z$ of degree 1
corresponding to the constant function $1$ in $\Gamma \cO (e)$,
$X$ of degree
2 corresponding to $x$,  and $Y$ of degree 3 corresponding to $y$.
These three variables satisfy the homogeneous form of the Weierstrass
equation. The statement that $A$ is the projective closure of
the plane curve defined by the Weierstrass equation may be restated
in terms of $\proj$:
\label{Aisproj}
$$A=\proj (\Gamma (\cO (*e))) .$$

%%FIELDS
\section{Torsion points and topology.}
\label{sec:torstop}

On the one hand, equivariant topology only gives counterparts 
to torsion points, but on the other it gives them greater importance.
This gives two  significant variations of the standard theory: we need
to use a different topology and we need to invert different sets of 
morphisms in forming the derived category.

\subsection{The torsion point topology.}
Because the topological model only gives counterparts of torsion points, 
we restrict sheaves to open sets which are complements of sets of points of 
finite order. This means that for us meromorphic functions
are  only allowed poles at points of finite order, and this entails
a number of other small effects that need attention.

The divisor $\Alr{n}$ of points of exact order $n$ will 
play a central role. Note  that
$$ A[n]=\sum_{s|n} \Alr{s}. $$

\begin{defn}
(i) Any divisor of the form $\sum_s a_s \Alr{s}$ (with $a_s \in \Z$)
is called  a {\em torsion point divisor}. \\
(ii) The {\em torsion point topology} on $A$ is the topology 
whose proper closed sets are specified
by a finite set $F$ of positive integers
$$V_F=\bigcup_{s \in F}\Alr{s}.$$
The non-empty open sets are thus $U_F := A \setminus V_F$.
\end{defn}

Since the sets $V_F$ are closed in the Zariski topology, we have 
a change of topology map $i: A_{\zar} \lra A_{\tp}$, and the
usual adjoint pair of functors
$$\adjunction{i^{-1}}{\shvtp}{\shvzar}{i_*} $$
between categories of sheaves. The restriction 
of topology functor $i_*$ is defined on Zariski presheaves $\cF$ by 
$i_*(\cF)(V):=\cF (V)$, which evidently takes sheaves to sheaves and
is exact. The extension of topology functor 
$i^{-1}$ is defined on torsion point presheaves $\cG$ by 
$(i^{-1}\cG)(U)=\cG (\hU)$, where $\hU =U_F$ where 
$F:=\{ n \st \Alr{n} \cap U=\emptyset\}$; this functor does not
preserve sheaves, so to obtain the sheaf level functor we pass
to associated sheaves. 

\begin{lemma}
The unit of the adjunction gives an isomorphism $i_*i^{-1}\cG \cong \cG$.
\qqed
\end{lemma}

To describe stalks it is convenient to use the notation 
$$\cF_{\infty}:= \colim_{n}\cF (A \setminus \bigcup_{p \leq n}\Alr{p}),$$
for sections with poles at any points of finite order, and 
$$\cF_s:= \colim_{n}\cF (A \setminus \bigcup_{p \leq n, \; p \neq s}\Alr{p})$$
for sections regular on points of exact order $s$ but with 
poles at any points of any other finite order. Note that these are
not Zariski stalks, but if we use the corresponding notation for
a torsion point sheaf $\cG$ we find $\cG_P =\cG_s$, where $s$ is the 
order of $P$. A short calculation then 
gives
$$(i_*\cF)_P =
\dichotomy{\cF_{\infty}& \mbox{ if $P$ is of infinite order}}
{\cF_{s}& \mbox{ if $P$ is of order $s$}}$$
and
$$(i^{-1}\cG)_P =
\dichotomy{\cG_{\infty}& \mbox{ if $P$ is of infinite order}}
{\cG_{s}& \mbox{ if $P$ is of order $s$, }}$$
so that  $i^{-1}$ preserves stalks.

Note that this means Zariski sheaves of the form $i^{-1}\cG$ are 
very rare, since the stalks at points of the same order are 
identical. In particular,  all stalks at points of infinite order are 
the same, suggesting there are no continuous families of sheaves of 
this sort. 

\begin{example}
We may restrict the Zariski structure sheaf $\cOzar$
to the torsion point topology, and we take $\cOtp:=i_*\cOzar$. 

Similarly,  our ring of meromorphic functions is
$$\cK =
\{ f \st f \mbox{ has poles only at points of finite order } \}, $$
with associated  constant sheaf $\cOit$.
Note that functions vanishing at points of infinite order are not 
invertible in $\cK$.

The local rings of the structure sheaf are thus
$$(\cOtp)_P =
\dichotomy
{\cK & \mbox{ if $P$ is of infinite order } }
{ \{ f \in \cK \st f \mbox{ is regular at points of exact order $s$ }\} 
& \mbox{ if $P$ is of finite order $s$} }$$
\end{example}

\begin{lemma}
The functors $i_*$ and $i^{-1}$ are both exact. 
\end{lemma}

\begin{proof}
The exactness of $i^{-1}$ follows since it preserves stalks.
For $i_*$, note that taking sections over $A \setminus F$
is exact  for any non-empty set $F$  of torsion points since it is 
affine; the stalks of $i_*$ are calculated as direct limits of
such functors.
\end{proof}

\begin{cor}
For Zariski sheaves $\cF$, the cohomology in the Zariski and 
torsion point topologies agree:
$$H^*_{\zar}(A;\cF) =H^*_{\tp}(A;i_*\cF). $$
\end{cor}

\begin{proof}
Since $\cF$ and $i_*\cF$ have the same global sections, and $i_*$ 
is exact, it suffices to note that if $\cI$ is flabby then 
$i_*\cI$ is a fortiori flabby too. 
\end{proof}

In future we will simply write $H^*(A;\cF)$ for the common value
of cohomology. Note that  this applies to the sheaves $\cOzar (D(V))$ of most
concern to us, and we will usually omit notation for the topology, 
writing simply $\cO (D(V))$.

% \begin{lemma}
% For Zariski sheaves $\cF$, the cohomology in the Zariski and 
% torsion point topologies agrees:
% $$H^*_{\zar}(A;\cF) =H^*_{\tp}(A;i_*\cF). $$
% For torsion point divisors $D$, 
% the cohomology in the torsion point topology agrees with the usual 
% Zariski cohomology.
% \end{lemma}

% \begin{proof}
% If  $D$ is a divisor of the specified form, the sheaf
%  $\cO (D)$ embeds in the constant sheaf $\cK$. 
% It therefore suffices to note that the constant sheaf $\cK$ is flabby
% for the torsion point topology.
% \end{proof}

\subsection{Torsion point equivalences.}
\label{subsec:tpequiv}
The previous subsection dealt with the change of topology, 
but there is the second issue of what set of morphisms are 
inverted to form the derived category. In equivariant topology
one does not usually invert all equivariant maps which are 
non-equivariant weak equivalences (since this gives only the homotopy 
theory of free actions). Instead, we invert only those 
equivariant maps which are equivalences in all fixed points.

We may transpose these considerations to sheaves of modules.
More precisely, $\cOzar$ is a sheaf of rings in the Zariski
topology and $\cOtp$ is a sheaf of rings in the torsion point
topology, and we may consider their respective categories  of 
modules, $\cOzarmod$ and $\cOtpmod$. These are both abelian categories, 
and related by the adjoint pair
$$\adjunction{i^*}{\cOtpmod}{\cOzarmod}{i_*}, $$
where 
$$i^*N:=i^{-1}N \tensor_{i^{-1}(\cOtp)}\cOzar.$$

\begin{lemma}
The unit of the adjunction gives an isomorphism 
$i_*i^*N \cong N$, so $\cOtpmod$ may be viewed as a subcategory 
of the category $\cOzarmod$.\qqed
\end{lemma}

\begin{lemma}
The functor $i^*$ is exact.
\end{lemma}

\begin{proof}
It suffices to prove that $\cOzar$ is flat over $i^{-1}\cOtp$, which 
we may verify at the level of stalks. This is straightforward
since $\cOzar (U)$ is flat over $i^{-1}\cOtp (U) =\cOzar (\hU)$ 
for any open set $U$.
\end{proof}

Derived categories are formed from abelian categories by taking 
a category of differential graded objects and inverting a suitable 
collection of morphisms. If all homology isomorphisms
are  inverted we obtain $D(\cOzarmod)$ and $D(\cOtpmod)$, 
but we wish to invert  fewer morphisms. The torsion point 
homology isomorphisms are those 
which induce isomorphisms of $H^*(A; \cdot \tensor \cO (D))$ 
for all torsion point divisors $D$, and we denote the derived 
categories obtained by inverting these $\Dtp (\cOzarmod)$ 
and $\Dtp (\cOtpmod)$.

To actually construct the derived categories we use cellular approximation.
This  is determined by specifying a set of spheres
$(\sigma_{\alpha})_{\alpha \in A}$ which must be small objects. An object 
is {\em cellular} if it is built from the spheres $\sigma_{\alpha}$ using 
arbitrary coproducts and triangles. A map $X \lra Y$ is a 
{\em weak equivalence} if it induces
an isomorphism of $[\sigma_{\alpha}, \cdot ]_*$ for all $\alpha$.
A {\em cellular approximation} of an object $X$ is then a weak equivalence
 $\Gamma X \lra X$ where $\Gamma X$ is cellular. We then work with the 
actual homotopy category of cellular objects.
For us the underlying category is the category of 
differential graded sheaves of $\cO$-modules in the appropriate topology
and the cells are the sheaves $\cO (D)$ where $D$ runs through torsion point 
divisors. 

 For clarity we display the relationship with the conventional derived
category of sheaves on $A$.

\begin{prop}
The derived categories are related by functors in the commutative
diagram 
$$\begin{array}{ccc}
\Dtp (\cOtpmod) &\lra &D(\cOtpmod)\\
i_*\uparrow \downarrow i^*&&i_*\uparrow \downarrow i^*\\
\Dtp (\cOzarmod) &\lra &D(\cOzarmod),
\end{array}$$
where the verticals are adjoint pairs with counits giving 
equivalences $i_*i^*N \simeq N$.
\end{prop}

\begin{proof}
The horizontals are elementary, since any torsion point homology isomorphism
is a homology isomorphism.

Since $i_*$ and $i^*$ are  exact, they preserves homology isomorphisms, and 
therefore induces  maps of derived categories. For torsion point 
homology isomorphisms we make additional arguments.
Indeed, 
$i_*\ShHom (M,N)=\ShHom (i_*M, i_*N)$ so that, taking $M=\cO (-D)$
we see that $i_*(N(D))=(i_*N)(D)$ and so $i_*$ preserves 
torsion point homology isomorphisms. Finally, 
$i^*(i_*M \tensor_{\cOtp} N) \cong M \tensor_{\cOzar}i^*N$, so taking
$M=\cO (D)$
 we see that $i^*$ preserves torsion point homology isomorphisms
as required.
\end{proof}

As remarked before, there is a far greater change in character in 
the vertical maps changing the topology than in the horizontal maps
changing the inverted morphisms. Even in $D(\cOzarmod)$
there are continuous families $\cO (P)$ of distinct objects.

%%FIELDS
\section{Coordinate data}
\label{sec:coordinatedata}

Our main theorem constructs a cohomology theory of type $A$ for
an elliptic curve $A$. The construction depends on a choice of function vanishing
at the identity, and the purpose of this section is to make clear the 
exact extent of this dependence. 

\subsection{The coordinate.}
Because the local ring $\cO_e$ in the torsion point topology 
is not quite the usual Zariski local ring, we make explicit
the  the properties we need. 

\begin{lemma}
The ideal 
$$I_e=\{ f \in \cO_e \st f(e)=0\}$$ 
of functions vanishing at the identity in $\cO_e$ is principal. 
The  generators of $I_e$ are exactly the  functions $t_e$
vanishing to first order at $e$ whose zeroes and poles are all at
points of finite order. 

If $t_e$ is a generator of $I_e$ then for any non-zero $f \in \cK$ there
is an integer $n$ such that $ft_e^n \in \cO_e$ and $ft_e^n$
does not vanish at $e$.
\end{lemma}

\begin{proof}
Suppose that $t_e$ is a function whose zeroes and poles are
at points of finite order with $t_e(e)=0$. Certainly $t_e \in I_e$;
on the other hand,  if $f \in I_e$,  
then  $f(e)=0$ so that  $f/t_e$ is still regular
at $e$, and only has poles at points of finite order. Hence
$f=t_e\cdot f/t_e \in (t_e)$ and $I_e =(t_e)$ as required. 
To see that this exhausts the set of generators, we note that 
a function $s \in I_e$ with a zero at a point $P$ of infinite
order is not a generator. Indeed, $I_e$ contains functions $f$ 
which do not vanish at $P$, and whenever $f=s g$, the function 
$g$ has a pole at $P$.

The final statement is clear since $t_e$ is a uniformizing element
in the Zariski local ring. 
\end{proof}

\begin{defn}
(i) A {\em coordinate} on $A$ (at the identity) is a generator
$t_e$ of the ideal $I_e$ in $\cO_e$ of functions vanishing at $e$. 
%This  also specifies a  uniformizer $t_P$ at $P$ by translating $t_e$.

(ii) A {\em coordinate divisor} is a divisor $Z_e$ of the form $\div (t_e)$ 
for some coordinate $t_e$.  By Abel's theorem, a torsion point divisor 
$Z_e=\Sigma_P n_P(P)$ with $n_e=1$ is a coordinate divisor
if and only if $\Sigma_P n_P=0$ and $\Sigma_P n_PP=0$.
\end{defn}

\begin{remark}
The ring $\cO_e$ is not a local ring in the sense of commutative
algebra: although $I_e$ is maximal, not all functions outside $I_e$
are invertible. However, the following lemma will provide the good behaviour
we need.
\end{remark}

\begin{lemma}
For any $s \geq 0$ the quotient $I_e^s/I_e^{s+1}$ is one dimensional 
over $k$, generated by the image of $t_e^s$. Hence
 $\cO_e/I_e^s$ is $s$-dimensional, generated by the images of
$1,t_e, \ldots , t_e^{s-1}$.\qqed 
\end{lemma}

We briefly discuss a special way of choosing coordinates.

\begin{defn}
A  {\em Weierstrass parametrization} of an elliptic curve is a choice of two functions
$x_e$ with a pole of order 2 at the identity and nowhere else, and
$y_e$ with a pole of order 3 at the identity and nowhere else. 
Because we work with the torsion point topology, we also require 
that  $x_e$ and $y_e$ only vanish at torsion points.
This Weierstrass parametrization  determines a coordinate $t_e=x_e/y_e$ 
of $\cO_e$.
\end{defn}

\begin{remark}
\label{choiceofcoord}
(i) The function $x_e$ is specified up to scalar multiplication by a
pair of non-identity points $A, B$ of finite order  with $A+B=e$ by 
the condition 
$\div (x_e)=-2(e)+(A)+(B)$. 
The function $y_e$ is specified up to scalar multiplication by 
three  non-identity points $C,D,E$ of finite order with $C+D+E=e$ by
$\div (y_e)=-3(e)+(C)+(D)+(E)$. This gives the coordinate divisor
$$\div (t_e)=(e)+(A)+(B)-(C)-(D)-(E).$$ 

(iii) One popular choice of Weierstrass parametrization involves choosing a point 
$P$ of order 2. This determines a  choice of $x_e$ and $y_e$   up to a constant multiple
by the conditions
$$\div (x_e)=-2(e)+2(P) \mbox{ and } \div (y_e)=-3(e)+(P)+(P')+(P'')$$
where  $A[2]=\{ e, P, P', P''\}$. Thus we obtain the coordinate divisor
$$\div (t_e)=(e)+(P)-(P')-(P'').\qqed$$
\end{remark}

\subsection{The cyclotomic functions.}
Once we have chosen a coordinate, this determines the choice of
 a function defining the points of exact order $s$.

\begin{lemma} 
\label{functionsts}
Given a choice $t_e$ of coordinate on the elliptic curve $A$, 
for each $s \geq 2$, there is a unique function $t_s$
with the properties
\be
\item $t_s$ vanishes exactly to the first order on $\Alr{s}$, 
\item $t_s$ is regular except at the identity $e \in A$
where it has a pole of order $|\Alr{s}|$, 
\item $t_e^{|\Alr{s}|}t_s$ takes the value 1 at $e$
\ee
Furthermore, the function $t_s$ only depends on the image
of $t_e$ in $\omega:=I_e/I_e^2$, and multiplying $t_e$ by 
a scalar $\lambda$ multiplies $t_s$ by $\lambda^{|\Alr{s}|}$.
\end{lemma}
\begin{proof} Consider the divisor $\Alr{s}-|\Alr{s}|(e)$.
Note that the sum of the points of $\Alr{s}$ in $A$ is 
the identity: if $s\neq 2$  this is because points occur
in inverse pairs, and if $s=2$ it is because the $A[2]$
is isomorphic to $C_2 \times C_2$. It thus follows from  
the Riemann-Roch theorem that 
there is a function $f$ with $\Alr{s} -|\Alr{s}|(e)$ as its
divisor. This function (which satisfies the first two properties
in the statement) is unique up to multiplication by 
a non-zero scalar. The third condition fixes the scalar, and replacing
the coordinate $t_e $ by $t_e+ft_e^2$ has no effect since
$t_st_e^{2|\Alr{s}|}$ vanishes at $e$.
\end{proof}

\begin{remark}
If we choose any finite collection $\pi =\{s_1, \ldots, s_s\}$ of
orders $\geq 2$, there is again a unique function 
$t_{\pi}$ with analogous properties.
Indeed, the good multiplicative property of the normalization 
means we may take 
$$t_{\pi} =\prod_{i} t_{s_i}.$$
This applies in particular to the set $A[n] \setminus \{e\}$. \qqed
\end{remark}

For some purposes, 
it is convenient to have a basis for functions with specified
poles. We already have the basis  $1,x,y,x^2, xy, \ldots $ if
all the poles are at the identity. Multiplication by a function 
$f$ induces an  isomorphism
$$ f\cdot  : \Gamma \cO (D) \stackrel{\cong} \lra \Gamma \cO (D-(f))  $$
so we can translate the basis we have.

\begin{lemma}
\label{shifttoe}
For the divisor $D=\Sigma_{s \geq 1} n(s) \Alr{s}$ 
let $t^*(D):= \prod_{b \geq 2} t_b^{n(b)}$. 
Multiplication by $t^*(D)$ gives an isomorphism
$$t^*(D) \cdot :  H^0(A;\cO (D))
\stackrel{\cong} \lra
H^0(A; \cO (\deg (D) \cdot (e))).$$

A basis of $H^0(A; \cO (D))$ is given by 
 $1/t^*(D)$ if $\deg (D)=0$, and by the first 
$\deg (D)$ terms in the sequence
$$ 1/ t^*(D), x/t^*(D), y/t^*(D), x^2/t^*(D), xy/t^*(D), \ldots $$
otherwise. \qqed
\end{lemma}

%\begin{remark}
%It is essential to be aware of the exceptional nature
%of the degree zero case.
%\end{remark}

\subsection{Differentials.}

On any elliptic curve we may choose an invariant differential, also 
characterized by the fact that it has no poles or zeroes. 
This is well defined up to scalar multiplication, and we would like to make a 
canonical choice.  Since $t_e$ vanishes to the first order at $e$, 
its differential is regular and non-vanishing at $e$, so we may 
take $Dt$ to be the invariant differential agreeing with $dt_e$ at $e$.

We shall be 
considering the space $\Omega \tensor_{\cO}\cK$ of meromorphic differentials: 
those which can be written in the form $fDt$ for a meromorphic function $f$.

\begin{warning} The differentials $dt_s$ are not generally meromorphic. 
To give an explicit example, suppose $A$ is defined by $y^2=x^3+ax +b$.
In this case,  the invariant differential is a scalar
multiple of  $dx/y$, and we may take $t_2$ to be a scalar multiple of $y$, 
so that  the zeroes of $dt_2$ are those of $dy=(3x^2+a)dx/y$. The four points 
at which $3x^2+a$ vanishes will not generally be torsion points.  
\end{warning}

It would be nice to make a construction which depends only on the coordinate
divisor and not the coordinate itself, but we only know how to do this for
a generic curve. We shall see that for such a construction,
it suffices to construct for each $s$ 
a meromorphic differential with poles to the first order on each point of
order $s$ which does not change if $t_e$ is multiplied by a scalar. 

 For $s=1$ the expression   $Dt/t_e$ gives a suitable  meromorphic 
differential. For $s \geq 2$, the situation is less straightforward.
To start with, by the last clause of \ref{functionsts}, the differential
$Dt/t_s$ does change if $t_e$ is multiplied by a scalar. 
Our next attempt is to note that  the differential $dt_s$ is again regular 
and non-vanishing
at each point $P$ of exact order $s$, and its value at $P$ is thus a
nonzero multiple $\lambda_P$ of that of $Dt$, but in general $\lambda_P$ does
depend on $P$. The differential $\lambda_PDt/t_s$ is suitable, but it 
involves making a choice of a particular point $P$ of order $s$. 
The alternative is to consider the average value 
$$\lambda_s =\frac{1}{|A\lr{s}|}\sum_{P \in \Alr{s}}\lambda_P$$
of the scalars and use the  differential $\lambda_s Dt/t_s.$
Provided $\lambda_s$ is non-zero, this gives a suitable differential 
depending only on the coordinate divisor.
However, for each $s$ there is a finite number of curves with $\lambda_s=0$, 
so it is only for a generic curve that this is legitimate. To avoid this
restriction we prefer to make a choice of coordinate rather than 
coordinate divisor.

\section{Principal parts of functions on elliptic curves.}
\label{secQ}

The point of this section is to analyze the sheaf 
$\cOit / \cO$ of principal parts of functions with poles
at torsion points. We repeat that we are working with sheaves
in the torsion point topology, so that $\cOit$ is the constant
sheaf corresponding to the ring $\cK$ of functions with arbitrary 
poles at points of finite order. 

For any effective torsion point divisor $D$ we may use the short exact sequence
$$0 \lra \cO \lra \cO (aD) \lra Q(aD) \lra 0$$
of sheaves to define the quotient sheaf $Q(aD)$ for $0 \leq a \leq \infty$.
%Furthermore the inclusions
%$\cO (aD) \lra \cO ((a+1)D)$ induce comparisons, and
%$$Q(\infty D )=\colim_a Q(aD) .$$
The cohomology of $Q(\infty D)$ is the cohomology of $A$ with
support  defined by  $D$.

In fact we may reduce constructions to the case when the divisor
$D=\Alr{s}$ for some $s$. Evidently, $Q(\infty \Alr{s})$ is a skyscraper
sheaf concentrated on $\Alr{s}$, so we may localize at $\Alr{s}$ to obtain
$$0 \lra \cO_{s} 
\lra \cO (\infty \Alr{s})_{s} \lra Q(\infty \Alr{s}) \lra 0.$$
Because we use  the torsion point topology, 
$$\cO (\infty \Alr{s})_s=\cO (\infty \Alr{s})_{\Alr{s}}=\cO (\infty tors)
=\cK .$$
Since $t_s$ is an invertible meromorphic function vanishing
 to the first order on $\Alr{s}$,  the sequence may be written
$$0 \lra \cO_s \lra \cO_s[1/t_s] \lra \cO_s/t_s^{\infty} \lra 0.$$
This gives the basis of a Thom isomorphism for the homology of almost 
free spectra. 

%Since $A$ is a smooth curve, the local ring $\cO_P$ is a discrete
%valuation ring, and if we choose a local uniformizer $t_P$ 
%any element of $\Gamma (Q(\infty P))$
%may be represented by an element of the form
%$$a_{-1}t_P^{-1}+a_{-2}t_P^{-2}+a_{-3}t_P^{-3}+ \cdots + a_{-n}t_P^{-n}$$
%for suitable scalars $a_{-i}$.  

\begin{lemma}
A choice of coordinate  gives isomorphisms
$$\cO ((a+r)\Alr{s})/\cO (r\Alr{s})=Q((a+r)\Alr{s})/Q(r\Alr{s}) \cong 
Q(a\Alr{s}), $$
induced by multiplication by $t_s^r$and hence
$$Q(\infty \Alr{s}) \otimes \cO (r\Alr{s}) \cong Q(\infty \Alr{s}).$$
If $s \geq 2$ the dependence is only through the image of $t_e$ in 
$\omega =I_e/I_e^2$.
\end{lemma}

\begin{proof}
Since the sheaves are all skyscraper sheaves over $\Alr{s}$, 
it suffices to observe that for any $a$, multiplication by $t_s$ 
induces an isomorphism
$$t_s: \cO ((a+1)\Alr{s})_{s} \stackrel{\cong}\lra 
\cO (a\Alr{s})_{s}.$$
To see this,  view the rings as subrings of the ring $\cK$ of meromorphic
functions. Since $t_s$ vanishes on $\Alr{s}$ and its  poles
are at points of finite order other than $s$, the image lies in the
stated subring. Multiplication by any non-zero function is injective, 
and to see the map is surjective, we observe that if $f \in \cK$ has
 no pole of order more than $a$ on $\Alr{s}$ then $f/t_s$ is a meromorphic 
function  no pole of order more than $a+1$ on $\Alr{s}$.
\end{proof}

Note that it is immediate from the Riemann-Roch formula that
for $0 \leq a \leq \infty$ the cohomology group
$H^0(A;Q(a\Alr{s}))$ is $a|\Alr{s}|$ dimensional, and  $H^1(A;Q(aP))=0$.

Now we may assemble these sheaves. Indeed,
we have a diagram
$$\begin{array}{ccccc}
\cO        &\lra & \cO (\infty D )      &\lra& Q(\infty D)\\
\downarrow &     & \downarrow           &    &\\
\cO        &\lra & \cO (\infty (D+D'))  &\lra& Q(\infty (D +D'))\\
\end{array}$$
of sheaves,  and hence a  map $Q(\infty D) \lra Q(\infty (D+D'))$.

\begin{prop}
If $s, s' \geq 1$ are distinct, then  the natural map
$$Q(\infty \Alr{s} ) \oplus Q(\infty \Alr{s'}) \stackrel{\cong}\lra 
Q(\infty (\Alr{s}+\Alr{s'}))$$
 is an isomorphism.
\end{prop}

\begin{proof}  
We apply the Snake Lemma to  the diagram
$$\begin{array}{ccccc}
\cO \oplus \cO & \lra & \cO (\infty \Alr{s} ) \oplus \cO (\infty \Alr{s'})& \lra &
Q(\infty \Alr{s}) \oplus Q(\infty \Alr{s'})\\
\downarrow && \downarrow &&\downarrow\\
\cO & \lra & \cO (\infty(\Alr{s}+\Alr{s'} )) & \lra &
Q(\infty (\Alr{s}+\Alr{s'}))
\end{array}$$
in the abelian category of sheaves on $A$.
The first vertical is obviously surjective with kernel $\cO$.
The kernel of the second vertical is also $\cO$, since if $f$ and $f'$ are
local sections of $\cO (\infty \Alr{s})$ and $\cO (\infty \Alr{s'})$
(i.e., $f$ only has poles on $\Alr{s}$ and $f'$ only on $\Alr{s'}$) then   
$f+f'=0$
implies that $f$ and $f'$ are regular. Finally  we must show that
$\cO (\infty (\Alr{s}+\Alr{s'}))$ is the sheaf quotient of $\cO \lra
\cO (\infty \Alr{s})\oplus \cO (\infty \Alr{s'})$. However, this may be verified
stalkwise, where it is clear.
\end{proof}

\begin{cor}
\label{splittingThom}
(i) The natural map gives an isomorphism
$$\bigoplus_s Q(\infty \Alr{s}) \stackrel{\cong}\lra Q (\infty tors). $$
(ii)  A choice of coordinate  $t_e$  gives an isomorphism
$$T_s: Q(\infty \Alr{s}) \otimes \cO (\Alr{s}) \stackrel{\cong} \lra
Q(\infty \Alr{s}).$$
(iii) The sheaf $Q(\infty \Alr{s})$ has no higher cohomology and its
global sections are 
$$\Gamma Q(\infty \Alr{s}) =\cK /\{ f \st f \mbox{ is regular on } \Alr{s} \}.
\qqed$$
\end{cor}

% \begin{remark}
% This corresponds to the fact that there is a rational splitting
% $$\Sigma \efp \simeq \bigvee_{H} \Sigma \eh$$
% where $\eh=cofibre(E[\subset H]_+ \lra E\subH_+)$
% \cite[2.2.3]{s1q}.
% \end{remark}

\part{The construction.}
In Part 4 we show that the structure of the algebraic model 
for rational $\T$-equivariant cohomology theories matches the
structure of sheaves of functions on an elliptic curve so 
neatly that the construction of a cohomology theory is effortless.
Short as it is, this is the core of the paper. 

%%FIELDS
\section{A cohomology theory associated to an elliptic curve.}
\label{sec:EA}

We are now ready to state and prove the main theorem.

\begin{thm}
\label{mainthm}
Given an  elliptic curve $A$ over a field $k$ of characteristic 0, and a 
coordinate $t_e$, there is an associated 2-periodic rational $\T$-equivariant 
cohomology theory
$EA_{\T}^*(\cdot)=E(A,t_e)_{\T}^*(\cdot)$
of type $A$, so that for any representation $W$ with $W^{\T}=0$ we have
$$\tEATi(S^W)=H^i(A;\cO (-D(W)))$$
and
$$\widetilde{EA}^{\T}_{-i}(S^W)=H^i(A;\cO (D(W)))$$
for $i=0,1$, where the divisor $D(W)$ is defined by taking
$$D(W)=\sum_{n}a_nA[n] \mbox{ when } W=\sum_na_nz^n \mbox{ with } a_0=0.$$
This association is invariant under base extension and 
 functorial for isomorphisms of the pair $(A, t_e)$.

The construction  is also natural for quotient maps $p: A \lra A/A[n]$ in the 
sense that  if the multiplicity of $p(t_e)$ at $e$ is 1 (for example if 
$\div(t_e)$ contains no points of order dividing $n$), there is a map 
$p^*: \infl_{\T /\T [n]}^{\T} E(A/A[n], p(t_e)) \lra  E(A,t_e)$ 
of $\T$-spectra, where $E(A/A[n])$ is viewed as a $\T /\T [n]$-spectrum 
and inflated to a $\T$-spectrum.
\end{thm}

\begin{remark}
(i) The elliptic curve can be recovered from the cohomology theory.
Indeed,  we may form the graded ring
$$\tEAThz(S^{*z}):=\{ \tEAThz(S^{az})\}_{a\geq 0}$$
from the products $S^{az} \sm S^{bz} \lra S^{(a+b)z}$, and
the  elliptic curve can be recovered from the cohomology theory via
$$A=\proj ( \tEAThz(S^{*z})), $$
as commented in Section \ref{secElliptic}. Furthermore, this reconstruction 
is functorial in that any multiplicative
natural transformation of cohomology theories will induce a map of 
elliptic curves. \\
(ii) In fact the coordinate can also be recovered from the 
cohomology theory, by evaluating the theory on suitable spaces
(see Proposition \ref{coorddiv} below).\\
(iii) A Weierstrass parametrization   of $A$ can be specified
by elements of homology:
$$x_e \in \widetilde{EA}_0^{\T}(S^{2z}) \mbox{ and } 
y_e \in \widetilde{EA}_0^{\T}(S^{3z}).$$
\end{remark}

\begin{remark}
\label{rationality}
(i) We have not required that $k$ is an  algebraically closed field.
 To see the advantage of this,  note that even for the multiplicative group,
the individual points of order $n$ are only defined over $k$ if $k$ contains
appropriate roots of unity. However $\G_m[n]$ (defined by $1-z^n$)
 and hence also $\G_m \langle n \rangle$ (defined by the cyclotomic polynomial
$\phi_n(z)$) are defined over $\Q$. Hence
 equivariant $K$-theory itself is defined over $\Q$. \\
%For an elliptic
%curve $A$ we require that there is a basis for $\Gamma \cO (a \Glr{d}))$
%consisting of functions defined over $k$. 
%%FIELDS
%%FIELDS
%%FIELDS
(ii) It is useful to generalize the construction to allow $k$ to 
be an arbitrary $\Q$-algebra so as to include various universal cases. 
There is no obstacle to making the construction in this generality, 
 provided functions $t_e$ and $t_s$ can be specified, but the analysis of
the resulting cohomology theory is  more problematic. Since
the entire construction is invariant under base change (provided we use 
corresponding coordinate functions), the case of a field already gives
significant information. The present methods are intrinsically restricted to 
$\Q$-algebras.
\end{remark}

\begin{remark}
\label{Galois}
One use for the naturality is that any automorphism of the elliptic curve
preserving the coordinate $t_e$ induces an automorphism of the
cohomology theory. For example if $t_e$ is defined using a point $P$ of order
2 as in Example \ref{choiceofcoord} (iii), any rigid Galois automorphism 
fixing $P$ gives an automorphism of the theory.
\end{remark}

\begin{proof} 
The basic ingredients of the torsion model of a the cohomology
theory associated to an elliptic curve $A$  are  analogous to the affine case
described in Appendix \ref{sec:GaGm}. We will write down a rigid,  even object
$$M_t(EA)=(\tf \tensor VA \stackrel{q}\lra TA)$$
of the torsion category $\cA_t$ (i.e., the structure map $q$ is
 surjective and $VA$ and $TA$ are in even degrees).
By \ref{surjMtgivesMs} this is intrinsically formal and therefore determines
$$M_s(EA )=(NA \lra \tf \tensor VA)$$
with $NA=\ker (q)$, and the representing spectrum $EA$.

We divide the proof into three parts: (1)
construction of $VA$ and $TA$, (2) construction of the map $q$
and (3) verification that the cohomology of spheres is correct.\\[1ex]

\noindent
{\em (1) The vertex and nub:}
Exactly as in the affine case, the degree 0 part of the  vertex
$$VA_0=\GcOit = \cK$$
consists of rational functions whose poles are all at torsion
points, however the torsion module is not simply the quotient of
this by regular functions, but rather
$$TA_0= \Gamma (\cOit/\cO ) =\Gamma (Q(\infty tors)).$$
Now we  use the splitting 
$$Q(\infty tors) \cong \bigoplus_s Q(\infty \Alr{s}) $$
of \ref{splittingThom} to separate points of different orders. 
This gives
$$TA_0=\Gamma Q(\infty tors) \cong \bigoplus_s \Gamma Q(\infty \Alr{s}) $$
where
$$\Gamma Q(\infty \Alr{s})= \cK / \{ f \st f \mbox{ is regular on } \Alr{s} \}
=\cK /\cO_s.$$

Both $VA$ and $TA$ are zero in odd degrees, and in
 other even degrees we take
$$VA_{2n}=\cK \tensor_{\cO} \Omega^n \cong VA_0 \tensor \omega^n \mbox{ and } 
TA_{2n}=\Gamma (\cK /\cO \tensor_{\cO} \Omega^n) \cong TA_0 \tensor \omega^n , $$
where $\Omega$ is the sheaf of K\"ahler differentials and $\omega$ 
is the cotangent space at the identity, and where  exponents refer to 
tensor powers (rather than exterior powers). 
We may now  describe the  $R$-module structure on $TA$. 
The direct sum splitting
$$TA=\bigoplus_s \Gamma Q(\infty \Alr{s})$$ 
corresponds  to the splitting 
$$R \cong \prod_s\Q [c],  $$
and 
$$e_s\TA=\Gamma Q(\infty \Alr{s})$$ 
is a $\Q [c]$-module where $c$ acts as multiplication 
by $t_s/Dt$, where $t_s$ defines $\Alr{s}$ as described in \ref{functionsts}. 
For $s=1$ this structure does not change if $t_e$ is multiplied by a 
non-zero scalar, so depends only on the coordinate divisor $Z_e$; for
$s \geq 2$ this depends only on the image of $t_2$ in $\omega =I_e/I_e^2$.
Since the order of any 
pole is finite, $e_sTA$ is a torsion $\Q [c]$-module. \\[2ex]
%Notice
%that the definition of the Thom isomorphism is arranged so that
%the composite
%$$\phi_{\Alr{s}} : Q(\infty \Alr{s}) =Q(\infty \Alr{s}) \otimes \cO
%\lra Q(\infty \Alr{s}) \otimes \cO (\Alr{s}) \cong Q(\infty \Alr{s})$$
%is multiplication by $t_s$.

\begin{remark}
The torsion module $TA$ may be described without using the coordinate
data. Indeed, we may define $TA'$ by giving its idempotent pieces
$$e_s(TA')_{2n}=\cK /\{ f \in \cK \st \ord_s(f) \geq n \}, $$
 and define the $\Q [c]$-action to be projection. A $\Q [c]$-isomorphism
 $TA' \cong TA$ is given by the coordinate:
$$(\frac{Dt}{t_s})^n: e_s(TA')_{2n} \lra e_s(TA)_{2n}.$$
We have used $TA$ rather than $TA'$ because the coordinate data does need
to be used somewhere, whilst differentials are used in a more uniform way 
in $TA$.
\end{remark}

\noindent
{\em (2) The structure map $q$:}
By \ref{mapsq} a map $q$ is determined by its  idempotent summands, 
which can be easily written down.\\[1ex]

\begin{defn}
\label{defn:q}
\label{defn:EA}
We define
$$q:\tf  \otimes VA \lra TA =\bigoplus_se_sTA$$
by specifying its $s$th component
$$q(c_s^{w(s)} \otimes \alpha )= \overline{
(\frac{t_s}{Dt})^{w(s)} \alpha };$$
up to normalization, this picks out the part of $\alpha$ with poles
of order $> w(s)$ on points of order $s$. 
\end{defn}

\begin{remark}
 Any  $\alpha \in V_{2n}$ may be written 
$$\alpha = f \cdot (Dt)^{\tensor n}$$
for some meromorphic function $f \in \cK$. The formula then becomes
$$q(c_s^{w(s)} \otimes f \cdot (Dt)^{\tensor n} )_s= \overline{
{t_s}^{w(s)} f  }\cdot (Dt)^{\tensor (n-w(s))}.$$
\end{remark}

\begin{lemma}
The definition does determine an $R$-map $q : \tf \otimes VA \lra TA$.
\end{lemma}

\begin{proof} 
Since any function is regular at all but finitely many points, the map $q$ 
maps into the sum. 
Thus $q(c^u \tensor \alpha)$ is well defined, and we need to check 
that taken together they specify an $R$-map. For this, we apply \ref{mapsq}.
Taking $V=VA$ and $T=TA$ we note 
that \ref{defn:q} does determine maps $q_s$, and that they 
 satisfy the condition. It follows that there is an $R$-map 
$q$ with these  idempotent pieces.
\end{proof}

\noindent
{\em (3) Cohomology:}
Now we can check that the resulting homology and cohomology of spheres agrees 
with the cohomology of the corresponding divisors on the elliptic curve.
Because the use of differentials is uniform, it is enough to prove the 
result for representations $W$ with $W^{\T}=0$.

Since we have decided to use the isomorphism
 $[S^{-w} , M]=[S^0, \Sigma^{w}M]$, we need to identify the 
suspension of the representing object $EA$. Applying 
\ref{suspofrigid} in this case we obtain the following.

\begin{lemma}
Suppose $w: \cF \lra \Z$ is zero almost everywhere.
 The $w$th suspension of $EA$ is given by 
$$ \Sigma^w EA =(\Sigma^w NA \lra \tf \tensor VA)$$
where 
$$\Sigma^w NA = \ker( \tf \tensor VA \stackrel{q^w} \lra \Sigma^wTA)$$
and for $\alpha \in VA_{2n}=\cK \tensor \omega^n$ 
$$q^w(c_s^{i(s)} \tensor \alpha  )=
\overline{ \alpha (\frac{t_s}{Dt})^{w(s)+i(s)}} \in 
(\cK /\cO_s) \tensor \omega^{n-w(s)-i(s)}=e_s(\Sigma^wTA)_{2n-2i(s)}.$$

We also use the mnemonic
$$q(xc^w \tensor \alpha ) = q^w(x \tensor \alpha), $$
despite the fact that  $xc^w$ is not an element of $\tf$.\qqed
\end{lemma}

Consider the complex representation $W$ with $W^{\T}=0$
and the corresponding function  $w(H)=\dimC (W^H)$. By 
\ref{rigidcohom} the homology is given by 
$$\tEATh_0 (S^W)=\ker (q : c^w \otimes VA_0 \lra (\Sigma^w TA)_0)$$
and 
$$\tEATh_{-1} (S^W)=\cok (q : c^w \otimes VA_0 \lra (\Sigma^w TA)_0)$$
and similarly with $W$ replaced by $-W$. 
Since the kernel and cokernel are  vector spaces over $k$, it is no loss
of generality to extend scalars to assume it is algebraically closed. 
This is convenient because it is simpler to treat separate points
of order $n$ one at a time. 

The following two lemmas complete the proof.
\end{proof}

\begin{lemma}
\label{Hzero}
If $W$ is a representation with $W^{\T}=0$ then 
$$\tEAThz (S^W)=H^0(A; \cO (D(W))), $$
and if $W \neq 0$, 
$$\tEAThz (S^{-W})=0.$$
\end{lemma}

\begin{proof} By definition 
$$q (c^w \otimes f)_s=\overline{(\frac{t_s}{Dt})^{w(s)}f }.$$
First note that $Dt$ is regular and non-vanishing on $\Alr{s}$, 
so the differential
can be ignored for the purpose of calculating the kernel.
Since the function $t_s$ vanishes to exactly the first order
on  $\Alr{s}$,  the condition that $f$ lies in the kernel is that 
$\ord_P(f)\geq -w(s)$ for each point $P$ of exact order $s$.
Since $D(W)=\Sigma_P w(s_P) (P)$ we have
$$\ker (q : c^w \otimes VA_0 \lra (\Sigma^wTA)_0)=
\{ f \in VA \st \div (f) + D(W) \geq 0\} $$
as required.

Replacing $W$ by $-W$, the second statement is immediate.
\end{proof}

The calculation of the odd cohomology is less elementary.

\begin{prop}
\label{firsthomology}
If $W$ is a representation with $W^{\T}=0$ then 
$$\tEATh_{-1}(S^{-W})=H^1(A; \cO (-D(W))), $$
and if $W \neq 0$, 
$$\tEATh_{0}(S^W)=0. $$
\end{prop}

\begin{proof}
We have to calculate $\cok (q : c^{-w} \otimes VA_0 \lra (\Sigma^w TA)_0)$.
First we give the concrete description of  $H^1(A; \cO (-D(W)))$ 
using ad\`eles from  \cite[Proposition II.3]{Serre}.

The  exact sequence of sheaves
$$0 \lra \cO (-D(W)) \lra \cK \lra Q(-D(W)) \lra 0.$$
induces a cohomology exact sequence ending
$$\cK \stackrel{\phi}\lra H^0(A;Q(-D(W))) \lra H^1(A;\cO (-D(W)))\lra 0. $$
However the definition of $Q(-D(W))$ shows that it is a skyscraper 
sheaf concentrated on the support of $D(W)$. 
Its space of sections is $\bA / \bA (-D(W))$, where 
$$\bA  =\{ (x_s)_s \st x_s \in \cK, \mbox{ and almost all } x_s \in k\}$$
is the space of ad\`eles  and 
$$\bA (-D(W))=\{ (x_s) \in W \st \ord_P(x_s) + \ord_s(-D(W)) \geq 0\}.$$
Thus $\cok (\phi )=\bA /(\bA (-D(W))+\cK)$. 

To complete the proof we  construct an isomorphism $m$ so that the 
left hand square in the diagram
$$\begin{array}{ccccccc}
\cK & \stackrel{\phi}\lra & H^0(A;Q(-D(W))) & \lra & H^1(A; \cO (-D(W)))
& \lra & 0\\
=\downarrow && m \downarrow \cong && \downarrow &&\\
c^{-w}\tensor (VA)_0 & \stackrel{q}\lra & (\Sigma^{-w}TA)_0 & \lra & 
\tEATh_{-1}(S^{-W}) & \lra & 0
\end{array}$$
commutes; the result follows from the 5-lemma. 
Both the domain and codomain of $m$ split into pieces corresponding to 
the divisors $\Alr{s}$. If $a_s=\dimC (W^{\T [s]})$, and we define
$m$ by taking  the $s$th term
$$m_s: \bA /\bA (-a_s \Alr{s})=\cK /\cO (-a_s \Alr{s}) \lra 
\cK /\cO_s \tensor \omega^{a_s}=e_s(\Sigma^{-w}TA)_0$$
to be 
$$m_s(\overline{f})=f \cdot (\frac{Dt}{t_s})^{a_s}.$$
Indeed, the definition is forced by the requirement that the 
square commute, but since the vanishing of $t_s$ defines $\Alr{s}$, 
$m_s$ is an isomorphism.
\end{proof}

\begin{remark}
It is possible to give a more explicit proof of \ref{firsthomology}
as follows. First, one checks  any element $(g_1,g_2, \ldots ) \in 
\bigoplus_s e_sTA$ is congruent (modulo the image of $q^w$)
 to one with $g_2=g_3=\cdots =0$.
Now, using \ref{shifttoe}, identify a subspace of the correct codimension 
in the image. Using
divisors one sees  the cokernel must be at least this big. Finally, 
the cokernel is naturally dual to $H^0(A;\cO (D(W))$, 
and hence naturally isomorphic to $H^1(A;\cO (-D(W)))$ by Serre 
duality. 
\end{remark}

\part{Properties of $\protect\T$-equivariant elliptic cohomology.}
Now that we have defined the cohomology theory $EA_{\T}^*(\cdot )$ 
associated to an elliptic curve $A$, we discuss some of its properties, 
including multiplicativity and a structure reflecting the addition on 
$A$.

\section{Homotopical multiplicative properties.}
\label{sec:mult}

For the rest of this section we identify $EA$ with 
the corresponding object in $\cA_s$, so that $EA=(NA \lra \tf \otimes VA)$,  
and  there is a short exact sequence
$$0 \lra NA \stackrel{\beta}\lra \tf \otimes VA \stackrel{q}\lra TA \lra 0.$$

\subsection{The ring structure on $EA$.}
Note that $VA=\bigoplus_n \cK \tensor_{\cO} \Omega^n$ has a commutative and 
associative product, which therefore induces such a product
on $\tf \tensor VA$.
%(so we may refer to it as an {\em algebra} of sections). 

\begin{thm}
\label{EAisring}
The product of functions and differential forms induces a
commutative and associative product on on the algebraic model for $EA$, 
so $EA$ is a commutative ring spectrum up to homotopy. Using results of
\cite{qtnq} we may choose $EA$ to be a strictly commutative ring $\T$-spectrum.
\end{thm}

\begin{proof}
First, note that by \ref{rigidisflat}  $EA$ is flat, 
so that  tensor product with $EA$ models the smash product.
It therefore suffices to show that the product on 
$\tf \tensor VA$ restricts to a product on $NA$.

Suppose $a, b \in \tf \tensor VA$; we must show that if
$q(a)=0$ and $q(b)=0$ then $q(ab)=0$. It suffices to concentrate
on the component mapping into $e_sTA$ for each $s$. The key to this is 
that for fixed $s$ we may give $VA$ the structure of 
a $\Q [c]$-module by letting $c$ act as $t_s/Dt$.  With this
definition, $c$ acts invertibly, so that we have a ring homomorphism
$$i_s: \Q [c,c^{-1}] \lra VA.$$
Now $q_s$ factors as
$$\Q [c,c^{-1}] \tensor VA \stackrel{i_s\tensor 1}\lra 
VA \tensor VA \lra VA \lra e_sTA.  $$
The fact that $q_s(ab)=0 $ if $q_s(a)=0$ and $q_s(b)=0$ 
now follows since the product of two functions regular
at a point is also regular there. 
 \end{proof}

\subsection{Duality.}
Now that we have a product structure we can tie up topological
and geometric duality in a satisfactory way.

\begin{lemma}
Spanier-Whitehead duality for spheres corresponds to 
Serre duality in the sense that the Serre duality 
pairing
$$\begin{array}{ccc}
H^1(A;\cO (-D(W))) \otimes H^0(A;\cO (D(W))) &\lra &H^1(A;\cO)\\
\| &&\|\\{}
[S^0,S^{-W}\sm \Sigma EA]^{\T} \otimes [S^0,S^{W}\sm  EA]^{\T} &&
[S^0, \Sigma EA]^{\T}
\end{array}$$
is induced by the algebraically obvious  Spanier-Whitehead 
pairing 
$$S^{-W} \sm EA \sm S^W \sm EA \simeq S^{-W} \sm S^W \sm EA \sm EA
\lra S^0 \sm EA \sm EA \lra  EA.$$
\end{lemma}

\begin{proof}
Both maps can be taken to be induced by multiplication of functions and 
a residue map (see \cite[Chapter II]{Serre}).
\end{proof}

\section{Reflecting the group structure of the elliptic curve.}
\label{sec:addition}

The group multiplication on an affine algebraic group $\G$ gives its ring
of functions $\cO$ a diagonal, and thus $\cO$ becomes a Hopf algebra.
When we say that $K$-theory corresponds to the multiplicative 
group  $\Gm$ we mean that not only is $K_{\T}^0=\Z [z,z^{-1}]$ the representing
ring for $\Gm$ but also that the diagonal also has a  topological source. 
Indeed, the 
multiplication map  $\mu : \T \times \T \lra \T$ induces a map
$$K^0_{\T}\stackrel{\mu^*}\lra K^0_{\T \times \T }= K^0_{\T}\tensor K^0_{\T}, $$
which turns out to be the coproduct on the ring of functions on $\Gm$.
%
%The counterpart in topology is that a non-equivariant
%complex oriented theory $E^*(\cdot )$ the map $B\T \times B\T 
%\lra B\T$ gives the  ring $E^*(B\T )$ a coproduct, and 
%makes it into a a complete topological Hopf algebra. It is in 
%this sense that the complex oriented theory gives rise to a formal 
%group.
The corresponding situation for  formal groups and complex oriented theories
is even more familiar. 

When we work with an elliptic 
curve, we again expect the group structure on $A$ to give 
additional structure on spaces of functions. However the 
structure is not just a coproduct, and we extract the relevant
information from Mumford's work \cite{Mumford}. Indeed, choosing a line
bundle $L$ to control the behaviour of functions, 
the multiplication $\mu : A \times A \lra A$
would give a map $\mu^*: H^0(A;L) \lra H^0(A \times A ; \mu^*(L))$, 
but since  $\mu^*(L)$ does not decompose as a tensor product, this is not 
very helpful.  Instead Mumford considers the map 
$$\xi : A \times A \lra A \times A$$
given by $\xi (x,y)=(x+y,x-y)$. It then turns out that if
we let $M=p_1^*L \tensor p_2^*L$, by the see-saw principle and the theorem 
of the square that  $\xi^*M \cong M^2$ (see \cite[p. 320]{Mumford}).
Using the K\"unneth isomorphism, we obtain a map 
$$\phi_L: H^0(A; L )\tensor H^0(A; L ) =H^0(A\times A ; M)
\stackrel{\xi^*} \lra H^0(A\times A ; M^2) =
H^0(A; L^2 )\tensor H^0(A; L^2 ) . $$
Applying this when $L=\cO (D(W))$ for  a representation $W$ 
with $W^{\T}=0$ we see that this is a map 
$$\phi_W : \tEATh_0(S^W)\tensor \tEATh_0(S^W) \lra
\tEATh_0(S^{2W})\tensor \tEATh_0(S^{2W}). $$
By choosing $W$ sufficiently large we can evidently find 
 $\xi^*(f_1,f_2)$ for an arbitrary meromorphic functions $f_1,f_2$, 
and since $\xi^*(f_1,f_2)(x,y)=(f_1(x+y),f_2(x-y))$, we recover 
$f_1(x+y)$ by suitable restriction.

We now describe how $\phi_W$ should be  realised at the level  of spectra.
The realization involves using $\T \times \T$-equivariant spectra, so 
proofs lie outside the scope of the present paper. However the picture
is sufficiently compelling to merit a brief account.

Suppose there exists a $\TxT$-equivariant cohomology theory
of type $\AxA$. Constructing such a theory is significantly 
easier than constructing a $\TxT$-equivariant theory for an arbitrary
abelian surface. To the representation $w^i\tensor z^j$ of $\TxT$ we associate
the divisor 
$$D(w^i\tensor z^j)=\ker (\AxA \stackrel{(i,j)} \lra \AxA), $$
and extend this to arbitrary representations so that  
$$D(V \oplus W )=D(V) + D(W). $$
The 2-periodic theory $\EAA^{\TxT}_*(\cdot )$ should then come with 
 a spectral sequence
$$H^*(\AxA ; \cOAxA (D(W))\Rightarrow \widetilde{\EAA}^{\TxT}_*(S^W).$$
  Since some line bundles have cohomology in degree 2, this
does not determine $\widetilde{\EAA}^{\TxT}_*(S^W)$ in general. 
However when $\cOAxA (D(W))$ has no cohomology in dimension 2
we find
$$\widetilde{\EAA}^{\TxT}_0(S^W)=H^0(\AxA ; \cOAxA (D(W))). $$
Next, the map $\xi : \AxA \lra \AxA$ is an isogeny with kernel
$$\Delta A[2]=\{ (a,a) \st a+a=e\}. $$
We also consider  the corresponding group homomorphism
$$\xihat: \TxT \lra \TxT, $$
defined by $\xihat (w,z)=(wz,w/z)$, which is surjective with kernel
$$\Delta \T [2]=\{ (z,z) \st z^2=1 \}. $$
To minimize confusion,  we identify the second $\TxT$ with 
$\TxTbar =(\TxT) /\DTtwo$.
The map    $\xi$ should correspond to a map 
$$\xi^*_i: \infl_{\TxTbar }^{\TxT} \EAA \lra \EAA$$
($i$ for inflation) of $\TxT$-spectra or,  adjointly, to a map 
$$\xi_f^*:  \EAA \lra \EAA^{\DTtwo} $$
($f$ for fixed point) of $\TxTbar$-spectra. 

\begin{lemma}
For any representation $\Wbar$ of $\TxTbar$, the map $\xi_f$ induces
$$\xi_f^*: \EAA^{\TxTbar}_0(S^{\Wbar}) \lra
\EAA^{\TxT}_0(S^{\Wbar}) .$$
\end{lemma}

\begin{proof} The map  $\xi_f^*$ induces
$$\xi_f^*: [S^0,S^{\Wbar} \sm \EAA ]^{\TxTbar}_0\lra 
[S^0,S^{\Wbar} \sm \EAA^{\DTtwo} ]^{\TxTbar}_0, $$
so it suffices to identify the domain and codomain. 
By definition $\widetilde{\EAA}^{\TxTbar}_0(S^{\Wbar}) 
=[S^0,S^{\Wbar} \sm \EAA ]^{\TxTbar}_0$
so we turn to the codomain and calculate
$$\begin{array}{rcl}
[S^0,S^{\Wbar} \sm \EAA^{\DTtwo} ]^{\TxTbar}_0
&=&[S^{-\Wbar},\EAA^{\DTtwo} ]^{\TxTbar}_0\\
&=&[S^{-\Wbar},\EAA ]^{\TxT}_0\\
&=&[S^0 , S^{\Wbar}\sm \EAA ]^{\TxT}_0\\
&=&\widetilde{\EAA}^{\TxT}_0(S^{\Wbar}).
\end{array}$$
\end{proof}

To model $M=p_1^*L \tensor p_2^*L$ with $L=\cO (D(W))$ 
 for a representation $W$ of  $\T$ we  
take $\Wbar = (W \tensor 1)\oplus (1\tensor W)$. 
Direct sum of representations
corresponds to tensor product of line bundles and to sums of divisors, so if 
$$W \mbox{ corresponds to the line bundle } L  \mbox{ and the divisor }
 D(W), $$
then 
$$\Wbar \mbox{ corresponds to the line bundle } p_1^* L \tensor p_2^*L
  \mbox{ and the divisor }  [D(W) \times A] + [A \times D(W)]. $$
Viewed as  a representation of $\TxT$ by pullback along $\xihat$ we find
$$\xihat^*(\Wbar)=\xihat_1^*W \oplus \xihat_2^*W.$$
In particular if $W=z^n$ we find
$$\xihat^*(\Wbar)=(w^n \tensor z^n)  \oplus (w^n \tensor z^{-n}).$$
Finally, we need to observe that for any $n$, the bundles 
associated to 
$$(w^n \tensor z^n)  \oplus (w^n \tensor z^{-n})
\mbox{ and } (w^{2n} \tensor 1)\oplus (1\tensor z^{2n})$$
are isomorphic: this is precisely the same argument 
as showed $\xi^*M\cong M^2$ above. With $L=\cO (D(W))$, 
we thus expect a commutative diagram
$$\begin{array}{ccccccc}
H^0(A; L)^{\tensor 2} &=&
H^0(A \times A; p_1^*L \tensor p_2^*L)& \stackrel{\xi^*}\lra & 
H^0(A \times A; p_1^*L^2 \tensor p_2^*L^2)&=&
H^0(A; L^2)^{\tensor 2}\\
&&\downarrow &&\downarrow &&\\
&&\widetilde{\EAA}^{\TxTbar}_0(S^{\Wbar}) &\stackrel{\xi_f^*}\lra&
\widetilde{\EAA}^{\TxT}_0(S^{\Wbar}) .&&
\end{array}$$

\section{The completion theorem.}

By formal completion around the identity, we may associate a formal 
group $\hA$ to an  elliptic curve $A$. In favourable circumstances
there is a (non-equivariant) $2$-periodic complex oriented 
cohomology theory $E\hA^*(\cdot )$ associated to $\hA$, and a 
Borel theory 
$$E\hA_{\T}^*(X):=E\hA^*(E\T \times_{\T} X). $$

The purpose of this section is to make explicit the relationship 
between the equivariant theory $EA_{\T}^*(X)$ associated to the 
elliptic curve $A$ and the Borel theory associated to the
formal group $\hat{A}$.

\begin{prop}
\label{completionthm}
The cohomology of $\ET$ is concentrated in even degrees, and 
in degree 0 it is the  completion of $\cO_e$ at the ideal $I_e$ 
of functions vanishing at $e$:
$$EA_{\T}^0(E\T )= \ilim_k \cO_e /I_e^k.$$
%  \Hom_{\cO}(\Sigma^{-1}Q(\infty (e)), \cO)=
%\Hom_{\Q [t_e]}^*(\Q [t_e]^{\vee}, \cK /\cO_e )
%=\ilim_k ( \cO ((k+1)(e))/\cO_e, t). $$
\end{prop}

\begin{proof}
Indeed, we may make the calculation
$$\tEATc^*(\ETp )=[\ETp , EA]_{\T}^*=[\ETp , EA \sm \ETp ]_{\T}^*
=\Hom_{\Q [c]}^*(\Q [c]^{\vee}, e_1 TA ).$$
Now, shifting into degree 0 we replace the action by $c$ with 
the action by $t_e$ and find this is 
$$\Hom_{\Q [t_e]}^*(\Q [t_e]^{\vee}, \cK /\cO_e )
=\ilim_k (\ann (\cK /\cO_e , t_e^{k}), t_e)
=\ilim_k ( \cO_e (k(e))/\cO_e, t_e). $$
Now multiplication by powers of $t_e$ gives an isomormphism between
the inverse system $(\cO_e (k(e))/\cO_e,t_e)$ and the inverse 
system $(\cO_e/I_e^k, projection)$.
\end{proof}

 Since the formal group law on $\hat{A}$ comes from 
$f(a+b)=F(f(a),f(b))$ when $f$ is a coordinate function, 
the formal group law for $E\hat{A}$ 
can be inferred from the map $\xi^*$ for 
$EA$ described in Section \ref{sec:addition}.

% Rather surprisingly, we do not 
% use the model $E\T = S(\infty z)$ together with the based cofibre
% sequence
% $$S(\infty z )_+ \lra S^0 \lra S^{\infty z}. $$
% The point is that we have good control over the cohomology of 
% spheres.

% \begin{lemma}
% $$\tEATc^i(S^{\infty z})=
% \dichotomy
% {\Gamma \cO (\infty (e))^{\vee} & \mbox{ if } i=0}
% {0 & \mbox{ if } i=1 }$$
% \end{lemma}

% \begin{proof}
% Since $S^{\infty z}=\colim_k S^{kz}$ we have a Milnor exact sequence
% for the cohomology. However, $\tEATc^{0}(S^{kz})=0$ for $k \geq 1$, 
% so we need only consider the cohomology in odd degrees. 
% Here we have
% $$\tEATc^{1}(S^{kz})=H^1(A;\cO (-k(e))) =H^0(A;\cO (k(e)))^{\vee} $$
% The map $S^{kz} \lra S^{(k+1)z}$ induces the dual of the inclusion
% $H^0(A;\cO (k(e))) \lra H^0(A;\cO ((k+1)(e)))$. It is therefore
% surjective, so that  there are no $lim^1$ terms. 
% \end{proof}

% This is easily enough to show that $EA_{\T}^*(ET)$ is in even 
% degrees, but only gives the desired cohomology up to extension. 

There is another less natural approach involving comparison with 
the Borel theory of the periodic theory represented by 
$$HP =\bigvee_{n \in \Z}\Sigma^{2n} H. $$
This has  coefficients
$$HP_{\T}^0=\Q [[y]] .$$

\begin{lemma}
\label{EAHP}
There is an equivalence $EA \sm E\T_+ \simeq HP \sm E\T_+$, 
and therefore
$$EA_{\T}^*(X \times E\T )\cong HP^*(E\T \times_{\T}X) , $$
so that in the notation above,  $E\hat{A}\simeq HP$.
\end{lemma}

\begin{remark}
The additional information in $E\hat{A}$ is in the comparison with 
$EA$, and hence in the relationship between the formal group law 
and the addition on $A$.
\end{remark}

\begin{proof}
First, to see the equivalence we need only show the two theories 
give  homology of $E\T$ isomorphic as $\Q [c]$-modules \cite[4.4.1]{s1q}. 
Since, both theories are 2-periodic and 
$EA^{\T}_*(E\T)$ and $HP_*(B\T)$  are divisible, 
it suffices to observe that the two theories have 
isomorphic non-equivariant coefficients.

Now for a based space $Y$, 
$$F( E\T_+ \sm Y,EA) \simeq F( E\T_+ \sm Y,E\T_+ \sm EA)
\simeq F( E\T_+ \sm Y,E\T_+ \sm HP)
\simeq F( E\T_+ \sm Y, HP). $$
\end{proof}

\section{The homology and cohomology of universal spaces.}

From the point of view of equivariant topology, the completion 
theorem of the previous section is just one example of a family 
of calculations. For other universal spaces we obtain analogous results
by the same proof. For simplicity we restrict the statement to the
value on a point.

Suppose then that  $\pi$ is a finite set of positive
integers and let $\cF (\pi)$ denote the family of subgroups with orders
dividing elements of $\pi$ and $A[\pi]$ denote the set of points with 
orders dividing elements of $\pi$.

\begin{thm} \label{Complct}
(i) {\em (Completion theorem.)}
The cohomology of $E\cF (\pi)$ is in even degrees and 
$$EA_{\T}^0(E\cF (\pi))=H^0(A;\cO_{A[\pi]}^{\wedge})$$
where $A[\pi]$ is the set of points with orders dividing
elements of $\pi$. Since $\cO_{A[\pi]}^{\wedge}$ is a skyscraper
sheaf, this is just the sum of the completed local rings at 
the points of $A[\pi]$.

(ii) {\em (Local cohomology theorem.)}
The homology of $E \cF (\pi)$ is in odd degrees and 
$$EA^{\T}_1(E\cF (\pi))=H^1_{A[\pi]}(\cO), $$
where the cohomology on the right is $A[\pi]$-local cohomology.
\end{thm}

\begin{proof}
The proof of Part (i) follows that of \ref{completionthm}.

For Part (ii) we may use the model
$$S(\infty V(\pi))=E\cF (\pi) \mbox{ where }
V(\pi)=\bigoplus_{n|\pi} z^n. $$
The cofibre sequence
$$S(\infty V(\pi))_+ \lra S^0 \lra S^{\infty V(\pi)}$$
and the fact that the Euler class of $z^n$ defines $A[n]$ 
give the result. 
\end{proof}

The calculation of the cohomology of $E\cF (n)=E(\T/\T [n])$ corresponds to the
fact that one may obtain a $\T [n]$-equivariant  formal group law
in the sense of \cite{Afgl} by formal completion of the curve $A$ along
$A[n]$, as described in \cite{LAm}.

\section{The Hasse square.}

We want to combine the localization and completion theorems to give
a method of calculation of elliptic cohomology in terms of Borel theories
combined using the geometry of the curve. 

The localization theorem is elementary. 
\begin{lemma} {\em (Localization theorem)}
\label{Localizationthm}
For any $\T$-space $X$ we have 
$$EA_*^{\T}(X \sm \etf)=H_*(X^{\T} ; \Omega_A^* \tensor_{\cO} \cK),$$
where the grading on the right is that for homology with graded
coefficients (i.e.,  total degree). 
A similar result holds in cohomology for finite complexes $X$.
\end{lemma}

\begin{proof}
Since $\colim_V \cO (D(V))=\cK$, and $\etf = \colim_{V^{\T}=0}S^V$ we have
$$EA_{2d}^{\T}(\etf)=(\Omega_A^1)^{\otimes d}\otimes_{\cO} \cK.$$
\end{proof}

We want to apply the completion theorem for the family of all finite
subgroups. To do this  for arbitrary complexes it is convenient 
to introduce the notation 
$$H_{\T}^*(X^C;I):=\Hom_{H^*(B\T_+)}(H_*^{\T}(X^C);I)$$
for any $H^*(B\T_+)$-module $I$, where the grading is that of 
homomorphisms of $H^*(B\T_+)$-modules. If $I$ is injective,  
this is a cohomology
theory in $X$, and if $H_*^{\T}(X^C)=H_*(X^C)\otimes H_*(B\T_+)$ then 
$H^*_{\T}(X^C; I)=H^*(X^C; \Hom_{H^*(B\T_+)}(H_*(B\T_+), I))$.

\begin{lemma}
For any $\T$-space $X$
$$EA_{\T}^*(X\sm \efp)=\prod_CH_{\T}^*(X^C; T_CA\tensor \omega_A^*).$$ 
If $H_*^{\T}(X^C)=H_*(X^C)\otimes H_*(B\T_+)$ for all $C$ then 
$$EA_{\T}^*(X\sm \efp)=
\prod_CH^*(X^C; \cO_C^{\wedge}\tensor \omega_A^*), $$
where $\cO_C^{\wedge}$ is the ring obtained as the formal completion of 
$\cO$ at $\Alr{s}$ if $C$ is of order $s$. 

\end{lemma}

\begin{proof}
The first statement amounts to the fact that 
$EA \sm \Sigma \efp$ is injective, with coefficients $TA \tensor \omega_A^*$. 
Now
we use the fact that there is a rational splitting 
$\efp \simeq \bigvee_C \elr{C}$ corresponding
to $TA\simeq \bigoplus_C T_CA$, and that $[X , \elr{C} \sm Y]^{\T}=
[X^C , \elr{C} \sm Y]^{\T}$. Passing to the summand corresponding to $C$, the 
$H^*(B\T_+)$-module  structure on rings of functions is through $t_s/Dt$.
The second statement follows since the short exact sequence
$$0 \lra \cK_C \lra \cK \lra T_CA \lra 0$$
gives an isomorphism
$$\Hom_{H^*(B\T_+)}(H_*(B\T_+), T_CA\tensor \omega_A^*)=
\Ext_{H^*(B\T_+)}(H_*(B\T_+), 
\cK_C\otimes \omega_A^*)=\cO_C^{\wedge}\tensor \omega_A^*.$$
\end{proof}

We express the homotopy level Hasse square via the associated
Mayer-Vietoris long exact sequence.

\begin{prop} (Hasse square)
\label{prop:Hasse}
For any $\T$-space $X$ there is a long exact sequence
\begin{multline*}
\cdots \lra EA^n_{\T}(X) \lra
H^n(\XT ; \cK \tensor_{\cO} \Omega_A^*) \times 
\prod_CH_{\T}^n(\XC; T_CA\tensor \omega^*_A) 
\lra H^n(\XT ; \cKF \tensor \omega_A^*)\\
\lra EA^{n+1}_{\T}(X) \lra \cdots , 
\end{multline*}
natural in $X$, where $\cKF = \prod_C \cO_C^{\wedge} \tensor \cK$.
If $H_*^{\T}(X^C)=H_*(X^C)\otimes H_*(B\T_+)$  then 
$$H_{\T}^n(X^C; T_CA\tensor \omega^*_A) \cong  
H^n(X^C; \cO_C^{\wedge}\tensor \omega^*_A).$$
\end{prop}

\begin{remark}
 Since $X$ is a space, two of the maps in the above long exact sequence
give a diagram of rings
$$\begin{array}{ccc}
EA^*_{\T}(X) &\lra&H^*(\XT ; \cK \tensor_{\cO} \Omega_A^*) \\
\downarrow&&\downarrow\\
\prod_CH_{\T}^*(\XC; T_CA\tensor \omega^*_A) &\lra& 
 H^*(\XT ; \cKF \tensor \omega_A^*).
\end{array}$$
When the connecting homomorphism in the long exact sequence is zero, this
is a pullback diagram of rings. For example, this applies  if both  
$H^*(X^{\T})$ and $H_{\T}^*(X^C)$ are in even degrees for all $C$.
\end{remark}

\begin{proof}
Any $\T$-spectrum $E$ occurs in the Tate homotopy pullback square
$$\begin{array}{ccc}
E &\lra & E \sm \etf \\
\downarrow &&\downarrow\\
F(\efp, E)&\lra & F(\efp, E ) \sm \etf
\end{array}$$
where $\cF$ is the family of proper subgroups, and applying $F(X, \cdot )$ we obtain the homotopy 
pullback square
$$\begin{array}{ccc}
F(X,E) &\lra & F(X,E \sm \etf) \\
\downarrow &&\downarrow\\
F(X \sm \efp, E)&\lra & F(X,F(\efp, E ) \sm \etf).
\end{array}$$
Note that $[X,Y \sm \etf]^{\T}_*=[X^{\T}, \PT Y]_*$, so that both the right hand
terms can be expressed in terms of the geometric fixed points of $X$.
Now take $E=EA$ and apply the localization theorem 
\ref{Localizationthm} to see 
 that $\piT_* (F(X,EA \sm \etf))=H^*(X^{\T}; \cK \otimes \omega_A^*)$
and the completion theorem \ref{Complct}(i)  to see that 
$$\piT_*(F(X \sm \efp , EA))=EA^*(X \sm \efp)=\prod_CH_{\T}^*(X^C; \cO_C^{\wedge}).$$
\end{proof}

\section{Recovering the coordinate.}

By showing that  the coordinate  used in Section 
\ref{sec:EA} can be recovered from the cohomology theory we show
that it is necessary to make such a choice. 

To give a full algebraic model of Type $A$ theories in the sense of 
\ref{typeG} we would need 
to show that  if $E_{\T}^*(\cdot)$ is a cohomology 
theory of Type $A$ then there is a unique coordinate so that 
$E_{\T}^*(\cdot )=E(A,t_e)_{\T}^*(\cdot )$. However it certainly requires
certain additional structure on the cohomology theory to do this.
First, we need to asssume that the theory is multiplicative
(this will mean it is specified by a collection of differentials
$\omega_s$ vanishing to first order at points of exact order $s$). 
However to relate the points of different orders we need to take into
account the group structure on $A$ and its reflection in cohomology.
We restrict ourselves to showing the required
uniqueness for theories constructed by the procedure of 
Section \ref{sec:EA}.

\begin{prop}
\label{coorddiv}
If $EA$ is constructed as in Section \ref{sec:EA}, the coordinate
 $t_e$ may be recovered from the cohomology theory. 
\end{prop}

\begin{proof} First we will recover the coordinate {\em divisor}, 
by concentrating
on point with trivial isotropy, and then return to find a suitable coordinate 
with this divisor by considering isotropy of order 2 and 3.

We evaluate the cohomology on suitable objects
$B=    (M \lra \tf \tensor U)$ of $\cA_s$ (depending on a number $n$ and
a representation $W$). 
These are certain  wide spheres  in  the sense of \cite[23.3]{s1q}, but
we give a self-contained description here.

Since our concern is mainly with what happens at
the identity, we separate the behaviours at and away from 
$e$ using idempotents. Indeed, we adopt the convention that 
$M'=e_1M$, $M''=(1-e_1)M$ and so forth. 
Away from the identity we take $B$ to be an 
ordinary  wedge of spheres
$$B''=(S^{0}\vee \Sigma^2 S^{-W})''$$
with $W^{\T}=0$. By choosing suitable 
representations $W$ this allows us to permit poles
away from the identity, for which we write
%$$\cK (W_0)''=\{ f \in \cK \st \ord_s(f) \geq - \dimC (W_0^{C_s})
%\mbox{ for }  s \geq 2\}$$
%and
$$\Omega (W)'':=\{ \alpha \in \Omega \tensor_{\cO}\cK \st 
\ord_s(\alpha ) \geq - \dimC (W^{C_s})
\mbox{ for }  s \geq 2\} . $$
The interesting part of $B$ is what happens at the identity
$$B'=(M' \lra \Q [c,c^{-1}] \tensor U). $$
First we take $U=\Q \oplus \Sigma^2 \Q$ with basis $b_0, b_2$
in degrees $0$ and $2$ (as forced by $B''$). 
Now take $M'$ to be the $\Q [c]$-submodule
of $\Q [c,c^{-1}]\tensor U$ generated by $a_0=1 \tensor b_0$
and $a_{2n+2}=c^{-(n+1)} \tensor b_0 + c^{-n} \tensor b_2$.

\begin{lemma}
The cohomology of the object  $B$ defined
 above (depending on $W$ and $n$) is given by 
$$\tEATc^0(B)=\{ (\lambda ,\alpha ) 
\in k \times \Omega (W)'' \st 
\ord_e( \lambda \frac{Dt}{t_e} + \alpha )\geq n\}. $$
\end{lemma}

\begin{remark}
Since $B$  has geometric $\T$-fixed points $S^0\vee S^2$, 
the identification with a subset of $\cK \times (\Omega \tensor_{\cO}\cK)$
is intrinsic.
\end{remark}

\begin{proof}
Consider a map $B \lra EA$ given by the diagram
$$\begin{array}{ccc}
M & \stackrel{\theta} \lra & NA\\
\downarrow && \downarrow\\
\tf \tensor U &\stackrel{1 \tensor \phi}\lra & \tf \tensor VA
\end{array}$$
Since $NA \subseteq \tf \tensor VA$, the map
is determined by the $R$-map $\theta : M \lra NA$.
Since $M \subseteq \tf \tensor (\Q \oplus \Sigma^2\Q)$ the map 
$\theta$  is determined by $f=\phi (b_0) \in VA_0=\cK$ and
 $\alpha =\phi (b_2) \in VA_2=\Omega\tensor_{\cO}\cK$. However
in order for $(f , \alpha)$ to determine such a map we need
to know the generators of $M$ map into $NA =\ker (q)$.

Exactly as in \ref{Hzero}, the condition away from $e$ is that 
$f$ is regular away from  $e$ and  $\alpha  \in \Omega (W)'' $. The condition 
at $e$ imposes the two conditions that 
 $\theta (c^{0} \tensor f) \in NA'_{0}$ and that 
$\theta (c^{-(n+1)} \tensor f + c^{-n} \tensor \alpha) \in NA'_{2n+2}$.
The first of these shows $f=\lambda$ is constant, and the second gives
the stated condition on $\alpha$.
\end{proof}

Now fix $\lambda =1$ (say), and consider the set 
$$\Lambda_{n,W}(t_e)=
\{ \alpha  \in \Omega (W)'' \st 
\ord_e(  \frac{Dt}{t_e} + \alpha )\geq n\}. $$
Finally,  suppose $t$ and $\tb$ are two choices of coordinate with
 $\Lambda_{n,W}(t)=\Lambda_{n,W}(\tb)$, then provided the two 
sets are non-empty (as we may assume by choice of $W$), we deduce 
$$\ord_e(\frac{Dt}{t} -
\frac{D\tb}{\tb})\geq n. $$
Expressing $t$ and $\tb$ in terms of a fixed coordinate $t_0$
we have $t=ut_0$ and $\tb =\ub t_0$ 
the condition is equivalent to requiring that 
$\frac{u(e)}{u} -\frac{\ub (e)}{\ub} $ vanishes to order $n$. Now, since
this is true for all $n$ and $u$ and $\ub$ are both non-zero at $e$, it
follows that $u /\ub$ is the scalar $u(e)/\ub (e)$. This shows that 
$EA$ determines the coordinate divisor.

Now choose a coordinate $t_0$ with the appropriate divisor, and consider
which multiple $t=\mu t_0$ gives the correct cohomology theory. For this 
we use a similar argument to the above with $n=0$, once with the idempotent
$e_1$ replaced by $e_2$,  and once with $e_1$ replaced with $e_3$.
Using $e_2$, we may pick out $\alpha$ satisfying the condition 
$$\ord_2(\frac{Dt}{t_2}+\alpha) \geq 0, $$
(this determines $\mu^3$). 
Using $e_3$,  we may pick out $\alpha$ satisfying the condition 
$$\ord_3(\frac{Dt}{t_3}+\alpha) \geq 0$$
(this determines $\mu^8$). These two together give $\mu$ as required. 
\end{proof}

\part{Categories of modules.}

We are working towards a comparison between the derived category 
$\Dtp (\cOtpmod)$ of $\tp$-sheaves of $\cO_A$-modules and a category 
of $EA$-modules. Before this can be useful, we need to describe methods
of calculation, and settle a number of technical difficulties. 

\section{Algebraic categories of modules.}
To start with,  we  work entirely in the algebraic category $\As$ with 
the strictly commutative ring $EA$ in $\As$.

\subsection{Modules over $EA$.}
 We may consider the  category $\EAmod$ of left modules over over the
algebraic model of $EA$. In fact a left $EA$-module 
$M=(P \lra \tf \otimes W)$  is given by a map $EA \otimes M \lra M$, 
or more explicitly, a diagram
$$\diagram
NA \otimes_R P \rto \dto &P\ddto\\
(\tf \otimes VA ) \otimes_R (\tf \otimes W ) \dto^=&\\
 \tf \otimes (VA \otimes W) \rto &\tf \otimes W
\enddiagram$$
From examples we see that we do not wish to require the structure
map $P \lra \tf \otimes W$  to be monomorphic, so we view
the $NA$-module $P$ as the basic object. Compatibility with 
the $VA$-module structure on $W$ imposes a further condition. 

Thus an $EA$-module is given by a suitably restricted
$NA$-module $P$. We view $P$ as a module of sections over 
the algebra $NA$ of regular sections.

It is worth making this more explicit for special types of 
object. If $M=e(W)$, then the module structure is simply 
the structure of a $\cK$-module on $W$.

If $M$ is torsion so that $M=f(T)$
 then $P=\bigoplus_sT_s$
where each module $T_sP$ is a module over 
$NA_s$, which is spanned by elements $c_s^i \tensor f$ with 
$t_s^if$ regular on $\Alr{s}$. Furthermore, the action of 
$NA_s$ factors through 
$\cO_s=\{ f \st f \mbox{ is regular on } \Alr{s}\}$.

\subsection{Homological algebra of the category of modules.}
\label{subsec:DTEAmod}
The purpose of this section is to describe the derived category 
$D_{\T}(\EAmod)$ of the algebraic category of modules, where the 
subscript $\T$ refers to the fact that only the counterparts of 
equivariant equivalences are inverted. We classify 
its objects up to isomorphism and give a means of calculating
maps. Since the $\tp$-derived category is formed by inverting 
maps which are homology 
isomorphisms for all twists, the maps are calculated in terms
of the corresponding relative Ext groups, which we now describe.

With sheaves it is convenient to work with flabby objects
rather than injective objects because we invert cohomology 
isomorphisms (i.e., isomorphisms of the derived functors
of global sections, or $\Ext_{\cO}^*(\cO, \cdot)$). 
There are enough flabby objects for homological dimension to be
visible at the level of abelian categories. We will work with a
corresponding class of $EA$-modules.

First we introduce the relevant test objects, namely the spheres and
torsion modules
$$\cT :=\{ EA \sm S^V \st V \mbox{ a complex representation}  \}
\cup \{ M \st \PT M=0\}.$$
The {\em $\tp$-flabby} objects are then given by 
$$\cI_F:= \{ I \st \ExtEA^s(T,I)=0 \mbox{ for all } T \in \cT, s \geq 1\}.$$
We next form an injective class by a process of saturation; the 
{\em $\tp$-monomorphisms} are 
$$\cM := M(\cI_F):=
 \{ f:X \lra Y \st f^*: \HomEA (Y,I) \lra \HomEA (X,I) 
         \mbox{ is epi for all } I \in \cI_F\}, $$
and the {\em $\tp$-injectives} by 
$$\cI := I(\cM):=\{ I \st f^*: \HomEA (Y,I) \lra \HomEA (X,I) 
\mbox{ is epi for all } f \in \cM \}.$$

First we need some examples of $\tp$-flabby objects.

\begin{lemma}
\label{injEA}
If $W$ is any
$\cK$-module, then  $e(W)$ is a $\tp$-flabby $EA$-module.

If $I=\bigoplus_sI_s$ with $I_s$ a divisible $\cO_s$-module, then 
$f(I)$ is a $\tp$-flabby $EA$-module.
\end{lemma}

\begin{proof}%Grisedale lemma
First note that modules of the form $e(W)$ admit 
injective resolutions of the same form, and similarly for those
of form $f(I)$. This means we can settle the question by considering
just Hom.

Next, we note that the  case $U=0$ of the condition holds
(i.e.,  $\ExtEA^s(EA,N)=0$ for $s>0$ for $N$ of the
specified forms). Indeed, 
$$\HomEA (EA,N)=\HomAs (S^0, N), $$
so it suffices for $N$ to be injective in $\As$, which is certainly 
the case for both $N=e(W)$ and $N=f(I)$ with $I_s$ being $c$-divisible.

For the modules $e(W)$  we use the adjunction 
$$\HomEA (M,e(W))=\HomK (V,W)$$
where $V$ is the vertex of $M$. The result when  $\PT M \simeq 0$ is clear
since it has zero vertex. The vertex of $S^U \sm EA$ 
is independent of $U$, the result follows from the case $U=0$.

For the modules $f(I)$,  we use the adjunction 
$$\HomEA (M,f(I))=\prod_s \HomOs (M_s,I_s)$$
where $M_s$ is the $s$th idempotent summand of the nub of $M$.
The result is clear since $I_s$ is injective by hypothesis.
%The result when  $\PT M \simeq 0$ follows since $I_s$ is injective.
%Since the idempotent summands of the nubs of  $S^U \sm EA$ 
%are suspensions of those for $U=0$, the result follows from 
%the case $U=0$.
\end{proof}

\begin{lemma}
The objects $\cI$ and the morphisms $\cM$ form an injective class and
a monomophic class.
\end{lemma}

\begin{proof}
By definition $\cI =I(\cM)$, and by saturation $\cM = M(\cI)$. 
It remains to show that for any $EA$-module $N$ there is a map
$f:N \lra F$ in $\cM$ with $F \in \cI$.

For an arbitrary $EA$-module $N=(L \lra \tf \tensor V)$ we have a map 
$N \lra \cEi N = e(V)$. The kernel $K$ is of the form $f(T)$
for a torsion module $N$, and we may embed this in a divisible
module $I$, giving a short exact sequence
$$0 \lra N \stackrel{i}\lra e(V) \oplus f(I) \lra f(J) \lra 0, $$
where $f(J)$ is divisible and hence also $\tp$-flabby. 

The fact that the map $i$ is $\tp$-monomorphic follows since
$f(J)$ is a test object. 
\end{proof}

% Indeed, if $\cI' \subseteq \cI$ is any set of $\tp$-injectives, 
% we may form the map 
% $$N \lra \prod_{I \in \cI'}\prod_{f:N \lra I} I$$
% with $f$th component being $f: N \lra I$. Since $\cI$ is closed under
% products, the codomain lies in $\cI$. The map is tautologically 
% surjective in $\HomEA (\cdot , I)$ for any $I \in \cI'$. It remains
% to observe that for any fixed $N$ we can choose a set $\cI'$
% so that any map $f: N \lra J$ with $J \in \cI$ factors through 
% an object of $\cI'$, and hence that the map is surjective in 
% $\HomEA (\cdot , J)$ and therefore in  $\cM$.
%[[FINISH THIS]]

This means that we can do relative homological algebra, and form
$\ExtEAtp^*(M,N)$. Better still, 
the proof supplied $\tp$-injective resolutions of length 1.

\begin{cor}
The $\tp$-injective dimension of any $EA$-module is $\leq 1$, so 
that $\ExtEAtp^s(M,N)=0$ for $s \geq 2$. Furthermore
$\HomEAtp (M,N)=\HomEA (M,N)$.\qqed
\end{cor}

This makes the category very accessible to calculation.

\begin{thm}
(i) All objects of $D_{\T}(\EAmod)$ are formal, in that $M \simeq H_*(M)$.
Thus homotopy types in $D_{\T}(\EAmod)$ correspond to isomorphism classes
of $EA$-modules. 

(ii) For $EA$-modules $M$ and $N$ there is a short exact sequence
$$0 \lra \ExtEAtp^1 (\Sigma H_*(M), H_*(N))
\lra [M,N]_{EA} \lra  \HomEA (H_*(M), H_*(N))\lra 0.$$
\end{thm}

The method of proof is standard, and slightly simplified 
by the fact that any $EA$-module can be considered as 
an object of $D_{\T}(\EAmod)$ by using the zero differential.

 We consider the map 
$$\nu: [M,N] \lra \HomEA (H_*(M), H_*(N))$$
given by taking homology. We will show that it is an isomorphism
for good $\tp$-flabby modules $N$. We have seen that any $EA$-module 
may be embedded in  a good $\tp$-flabby module with $\tp$-flabby quotient.  
Now, for an arbitrary differential graded $EA$-module 
$N$ we choose a $\tp$-resolution 
$$0 \lra H_*(N) \lra I_0 \lra I_1 \lra 0$$ 
of its homology, where $I_0$ and $I_1$ are both good $\tp$-flabby modules.
Now let $N \lra I_0$ be the  map corresponding to the first map 
in the resolution and note that the mapping cone has homology $I_1$.
Up to isomorphism we therefore have a cofibre sequence
$$ N \lra I_0 \lra I_1, $$
and applying $[M,\cdot ]_{EA}$ we obtain Part (ii) of the theorem.
Part (i) now follows, since if $H_*(M) \cong H_*(M')$ we may lift
this isomorphism to a map $M \lra M'$, which, being a homology 
isomorphism, is an equivalence. In particular $\tp$-flabby objects are
classified by their homology, so it was reasonable to call the 
cofibre $I_1$. It remains to prove that our good $\tp$-flabby modules
have the right properties. 

\begin{lemma}
If $N$ is one of  the modules $e(W)$ and $f(T)$ in \ref{injEA}, 
the map $\nu$ is an isomorphism.
\end{lemma}

\begin{proof}
By definition the functor $\HomEA ( \cdot , H_*(N))$ is exact
when $N$ is $\tp$-flabby, so we have a natural transformation of cohomology
theories and it suffices to check it is an isomorphism on a collection 
of $EA$-modules which generate all modules using direct sums and
cofibre sequences. By Adams's projective resolution argument, it
suffices to use the objects 
$EA \sm S^V$, since they are small and detect weak equivalences.
The objects $EA\sm S^V$ are extended by construction, and we
have
$$\piA_*(EA \sm S^V) \cong \piA_*(EA) \tensor \piA_*(S^V), $$
and hence a commutative diagram
$$\diagram
[EA \sm S^V, N]_{EA} \rto \dto_{\cong} & 
\HomEA (\piA_*(EA\sm S^V), \piA_*(N))\dto^{\cong}\\
[S^V, N]\rto  & \HomAs (\piA_*(S^V), \piA_*(N)).
\enddiagram$$
The result follows from the fact that the objects are injective
in $\As$ together with the corresponding statements there
 \cite[5.6.7,5.6.8]{s1q}.
\end{proof}

\section{Homotopy modules.}
\label{sec:homod}

The equivalence of \cite{s1q} is only defined at the homotopy 
level and the equivalence of \cite{shipley} is not known to 
be monoidal at the model category level. The results of \cite{qtnq}
do show  that we may choose $EA$ to be a strictly commutative ring 
spectrum, and hence there is a model category of $EA$-module $\T$-spectra, 
but since this is not yet published, it seems worth including a brief account
of what can be said about modules up to homotopy:
this section will discuss how good 
a model of $D_{\T}(\EAmod)$ can be obtained by working with 
rings and modules up to homotopy.

Modules up to homotopy  have notoriously bad 
formal behaviour, but  the low homological dimension of the algebraic 
categories means we can nonetheless obtain some useful information.
The idea is to use the category of homotopy modules and homotopy 
module maps as a  model for the homotopy category of modules.
To see the effectiveness of this, we continue to work in the 
algebraic category.

At the level of objects, the model is good. 

\begin{lemma}
Every homotopy $EA$-module is represented by a strict $EA$-module.
Two homotopy $EA$-modules are equivalent if and only if their
strict representives are equivalent. 
\end{lemma}

\begin{proof}
Any $EA$-module $M$ is obviously a homotopy module. Since the orginal 
module may be recovered via the action
of $EA=\piA_*(EA)$ on $\piA_*(M)$, the forgetful map is injective
on objects. Furthermore, every object of $\As$ is formal, so 
there is an equivalence $M \simeq \piA_*M$, and the action passes
to $\piA_*(M)$. Thus any homotopy module $M$ is equivalent to 
the strict module $\piA_*(M)$, and  the forgetful map 
is surjective. 
\end{proof}

Given two homotopy modules $M,N$, we define the group of 
homotopy module maps by 
$$\mapho{M,N}:= \{ f \in [M,N] \st f \mbox{ is a module map up to homotopy }
\}.$$
The main point to make is that this is a subset of the maps 
ignoring $EA$-module structure, so that it is unlikely to model 
phenomena of positive filtration. As usual the cofibre of a 
homotopy module map has no canonical structure as a homotopy module.
Taking homotopy module maps need not be exact, even if applied to a
cofibre sequence of strict modules. 

The best we can do is to attempt to detect homotopy module maps.
Given homotopy modules, we choose strict modules $M, N$ representing them, 
and to simplify the notation, we assume they have zero differential. 
We then have a diagram
$$\diagram
0 \rto & \ExtEAtp^1 (\Sigma M, N) \ddto \rto &
 [M,N]_{EA} \dto \rto &  \HomEA (M, N) \ddto \rto &0\\
       &                                           &
\hspace*{6ex}\mapho{M,N}\dto       &                                     &\\
0 \rto & \ExtAs (\Sigma M, N) \rto &
 [M,N] \rto &  \HomAs (M, N) \rto &0.
\enddiagram$$
Since every module map in homology is represented by a strict module map 
$M\lra N$, it is represented by a homotopy module map. Subtracting
this, the remaining issue is how to decide when a map inducing zero 
in homology is a homotopy module map. Certainly it suffices for it
to be in the image of $\ExtEAtp^1 (\Sigma M,N) \lra \ExtAs (\Sigma M,N)$.
When no other elements of $\ExtAs (\Sigma M,N)$ represent
homotopy module maps, the forgetful map  
$$[M,N]_{EA}\lra \mapho{M,N}$$
is surjective, but even then its kernel is 
$$\ker \left[\; \ExtEAtp^1 (\Sigma M,N) \lra \ExtAs (\Sigma M,N) \; \right],$$ 
which may be non-trivial.

\part{An equivalence between derived categories of sheaves and spectra.}
Having shown that the structure sheaf $\cO_A$ of the elliptic curve
$A$ gives rise to a commutative ring $EA$ in $\As$, we show in this part 
that this extends to an equivalence between their derived categories of 
modules. 

The discussion of modules up to homotopy in Section \ref{sec:homod} shows 
how much of the resulting information can be transported to the category 
of spectra without using further technology. However, the results of 
\cite{qtnq} show that the strictly commutative ring in $\As$ gives
a {\em strictly} commutative ring $\T$-spectrum, and using this additional 
technology, the present account applies without change to categories of 
equivariant $EA$-module {\em spectra}.

%\section{Sheaves of $\protect \cO$-modules and $\protect EA$-module 
%$\protect \T$-spectra.}

\section{Sheaves from spectra.}
\label{sec:shfromsp}

We describe a natural construction of a sheaf over $A$ from a $\T$-spectrum. 
In Section \ref{sec:Grojnowski} we  show how it
is related to Grojnowski's construction \cite{grojnowski}. 

\subsection{Sheaves associated to $R$-modules.}
\label{sec:sheavesfromRmod}
An object of $\As$ is a {\em based} $R$-module in a suitable sense, 
but it will clarify the later construction to begin with a construction 
on arbitrary $R$-modules $N$. 

Note first that we have defined suspension functors 
$\Sigma^wN$ for any almost constant function $w: \cF \lra \Z$, 
and if $w(s)\leq w'(s)$ for all $s$ there is  a map 
$\Sigma^wN \lra \Sigma^{w'}N$ which is multiplication by 
$c^{w'(s)-w(s)}$ on the $s$th idempotent summand. 

Recall that,  for any finite set $\pi$ of positive integers, 
$V_{\pi}$ is the set of points of $A$ whose orders are in $\pi$, and
$U_\pi=A \setminus V_{\pi}$.

\begin{defn} Suppose  $N$ is an $R$-module and let
$$\cEi_{\pi}N:=\colim_{w(\pi)=0} \Sigma^wN.$$

 Now define  a presheaf $\tN$ of $R$-modules on $A$ by taking
$$\tN (U_{\pi}):=\cEi_{\pi}N$$
\end{defn}

\begin{lemma}
\label{tNsheaf}
The presheaf $\tN$ is a sheaf.
\end{lemma}

\begin{proof}
First note that, since $\cEi_{\pi}N=\cEi_{\pi}R \tensor_R N$,  we have 
$\tN = \tR \tensor_R N$. 

Now since any cover has a finite subcover, it suffices to check 
the sheaf condition on $U_{\pi \cap \pi'}=U_{\pi} \cup U_{\pi'}$. 
Since $\cEi_{\pi}R$ is flat for any $\pi$, it suffices to deal 
with the special case $N=R$, where  we have an exact sequence
$$0 \lra \cEi_{\pi \cap \pi'}R
\lra \cEi_{\pi}R\oplus \cEi_{\pi'}R
\lra \cEi_{\pi \cup \pi'}R .$$
\end{proof}

\subsection{Construction of the sheaf.}
\label{sec:sheavesfromTspec}
We begin in earnest by defining a  functor 
$$\cMA: \As \lra \tpsheafA$$
at the level of abelian categories. We will show that it restricts
to a functor 
$$\cMA: \EAmod \lra \tpOmod . $$

When $\pi$ is the set of divisors of $n$ we think of $V_{\pi}$ 
as defined by the Euler class of $z^n$. This motivates some 
corresponding definitions in equivariant topology. For each 
subgroup $H$ we need the space 
$E\langle H \rangle =\cofibre (E[\subset H]_+ \lra E[\subseteq H]_+)$, 
whose distinguishing feature is that its $K$-fixed points are 
contractible unless $K=H$, and $S^0$ if $K=H$. We consider 
the set $\cFlrpi$ of subgroups of $\T$ with order in $\pi$ and 
then form the space
$$E\langle \pi \rangle :=e_{\cFlrpi}E\cF_+=
\bigvee_{H \in \cF (\pi)} E\langle H \rangle , $$
where $e_{\cFlrpi} \in \map (\cF ,\Q)$ is the idempotent
with support $\cFlrpi$. We may then form the 
 space $\Etlr{\pi}$ using the cofibre sequence
$$E\langle \pi \rangle \lra S^0 \lra \Etlr{\pi}. $$
The space $E\langle \pi \rangle$ is modelled in $\As$ by 
$$T(\pi)=(\bigoplus_{H \in \cFlrpi} \I (H) \lra 0)$$
and the space $\Etlr{\pi}$ by 
$$L(\pi)=(R(\infty \pi) \lra \tf )$$
where $R(\infty \pi) \subseteq \tf$ consists of elements with poles
only  on $\cFlrpi$.

Next, we associate a sheaf $\cMA (X)$ in the torsion-point topology 
with  an object $X=(P \lra \tf \otimes W)$ of $\As$. First recall the notation 
$$P(c^0)=\{ p \in P \st \beta (p) \in c^0 \otimes W\}=
\Hom_{\As} (S^0, X).$$
Continuing the analogy with sections, we write
$$P(\infty \pi)=P \otimes_R R(\infty \pi) $$ 
so that 
$$X \otimes L(\pi) = (P(\infty \pi) \lra \tf \otimes W). $$

\begin{defn}
\label{defn:MA}
For any object $X=(P \lra \tf \otimes W)$ of $\As$  
the presheaf $\cMA (X)$ is defined by 
$$\cMA (X)(U_{\pi})= \HomAs (S^0, X \tensor L(\pi))= P(\infty \pi) (c^0) . $$
The restriction associated to $U_{\pi'} \subseteq U_{\pi}$ is induced by 
the map $L(\pi) \lra L(\pi')$ which is the identity on the vertex. 
\end{defn}

\begin{lemma}
The presheaf $\cMA (X)$ is in fact a sheaf.  
\end{lemma}

\begin{proof}
It suffices to consider the cover of $U_{\pi \cap \pi'}$ by 
$U_{\pi}$ and $U_{\pi'}$, and we need to 
show there is an exact sequence 
$$0 \lra \cMA (X) (U_{\pi \cap \pi'}) \lra 
\cMA (X) (U_{\pi}) \oplus \cMA (X) (U_{\pi'}) 
\lra \cMA (X) (U_{\pi \cup \pi'}) .$$
This is obtained  from the short exact sequence 
$$0 \lra L(\pi \cap \pi') \lra L(\pi) \oplus L(\pi') \lra L(\pi \cup \pi') \lra 0. $$
Indeed, since $L(\pi \cup \pi')$ is flat, we obtain the desired exact 
sequence by applying the functor \\
$\Hom_{\As}(S^0, X \tensor \cdot )$.
%$$0 \lra \Q \lra \Q \oplus \Q \lra \Q \lra 0$$
%on vertices and  the exact sequence 
%$$0 \lra R(\infty (\pi \cap \pi')) \lra R(\infty \pi) \oplus R(\infty \pi') 
%\lra R(\infty (\pi \cup \pi')) \lra 0 $$
%on nubs. 
\end{proof}

\begin{lemma}
\label{MApreservesvertex}
If $X$ has vertex $V$ then 
$$\cMA (X) (\infty tors)=V.$$
\end{lemma}

\begin{proof}
Since 
$\cMA (X)(U_{\pi})=\HomAs (S^0, X \tensor L(\pi ))$ and $S^0$ is 
small,  we find
$$\cMA (X)(\infty tors)
 = \colim_{U_{\pi}}\HomAs (S^0, X \tensor L(\pi ))
 = \HomAs (S^0, \colim_{U_{\pi}}X \tensor L(\pi )) =V. $$
\end{proof}

For torsion free spectra we can also identify stalks.

\begin{lemma}
\label{stalkMAX}
If $X$ is torsion free then the stalk of $\cMA (X)$ at a point of 
order $s$ is given by 
$$\cMA (X)_s =\ker (c^0 \tensor V \lra T \lra e_sT).  $$

\end{lemma}

\begin{remark}
It is natural to refer to $\cMA (X)_s$ as the space of  
 $X$-meromorphic functions regular at $s$.
\end{remark}

\begin{proof}
To calculate the stalk we take a direct limit over $U_{\pi}$
containing points of order $s$, which are $U_{\pi}$ with $s \not \in \pi$.
For a torsion free $X$ 
$$\cMA (X)(U_{\pi}) =\ker (c^0 \tensor V \lra T \lra 
\bigoplus_{r \not \in \pi} e_rT). $$
\end{proof}

Since direct sums commute with tensor products and $S^0$ is
small, we deduce a useful formal property.

\begin{lemma}
\label{cMApreservessums}
The functor $\cMA$ preserves arbitrary direct sums. \qqed
\end{lemma}

\subsection{The sheaf associated to an $EA$-module.}

We show that applying the functor to an $EA$-module
gives a sheaf of $\cO_A$-modules.

\begin{lemma}
(i)  The functor $\cMA$ takes $EA$ to the structure sheaf:
$$\cMA (EA)=\cO_A.  $$

(ii)  The functor $\cMA$ takes $EA$-modules to $\cO_A$-modules, and
therefore induces a functor
$$\cMA : \EAmod  \lra \tpOmod .$$
\end{lemma}

\begin{proof} Part (i)  is clear from our construction of elliptic 
cohomology. 

For Part (ii), we need to show that there are structure maps
$\cO (U_{\pi}) \tensor \cMA (X)(U_{\pi}) \lra \cMA (X)(U_{\pi})$, or in 
other words, 
$$NA(\infty \pi)(c^0) \tensor P (\infty \pi)(c^0) \lra
P (\infty \pi)(c^0) .$$
However we need only note that, $L(\pi)$ (like $S^0$) is idempotent 
in the sense that $L (\pi ) \tensor L(\pi )=L(\pi)$ so that the required map 
is the composite 
\begin{multline*}
\Hom_{\As} (S^0, EA \tensor L(\pi)) \tensor 
\Hom_{\As} (S^0, X \tensor L(\pi)) \stackrel{\tensor} \lra 
\Hom_{\As} (S^0 \tensor S^0, EA \tensor L(\pi) \tensor 
X \tensor L(\pi)) \\
=\Hom_{\As} ( S^0, EA \tensor X \tensor L(\pi)) 
\lra \Hom_{\As} ( S^0,  X \tensor L(\pi)) .
\end{multline*}
Compatibility with restriction is clear since the restriction 
associated to $U_{\pi'} \subseteq U_{\pi}$ is induced by a map 
$L(\pi) \lra L(\pi')$.
\end{proof}

One more special value plays an important role. 

\begin{lemma}
\label{MAsphereline}
The $EA$-module   $S^W \sm EA$ is taken to the corresponding line bundle
$$\cMA (S^W \sm EA)=\cO_A (D(W)). $$
\end{lemma}

\begin{proof} 
This follows directly from the parallel between topological suspension
\ref{defn:suspAs} and algebraic twisting by line bundles. 
\end{proof}

\section{$\protect\mathbb{T}$-spectra from $\protect \mathcal{O}_A$-modules.}
\label{sec:TspecfromOmod}

In this section we adapt the construction of $EA$ given 
in Section \ref{sec:EA} to associate an object of  $\cA_s$
to an $\cO_A$-module, and hence provide a functor
$$\SA : \tpOmod \lra \EAmod .$$
In the    construction of $EA$ we made fundamental use of 
the fact that the sheaf $\cO_A$ is torsion free in the sense that
$\cO (D)$ is a submodule of $\cO (\infty tors)=\cK$ for all 
torsion point divisors $D$. As a result, 
the nub is  a submodule of $\tf \tensor VA$, where
$VA_0$ consists of the space $\cK$ of  meromorphic functions. For an 
$\cO$-module $\cY$, it often happens for a non-zero sheaf $\cY$, 
 that the sheaf $\cY (\infty tors)$ of meromorphic sections is zero, 
so that the earlier construction would give zero. 
The construction we give here does specialize to construct $EA$, 
but also deals with torsion sheaves.

\subsection{The construction.}

In topology, the object of $\cA_s$ associated to a $\T$-spectrum $X$
is obtained from  the map 
$$X \sm D\efp \lra X \sm D\efp \sm \etf $$
by taking equivariant homotopy groups. The key facts are 
\begin{itemize}
\item $X \sm D\efp \sm \etf \simeq \Phi^{\T} X \sm D\efp \sm \etf$ and
\item the cofibre of the map is the $\T$-free object 
$X \sm \Sigma \efp$
\item there is a  cofibre sequence
$$D\efp \lra D\efp \sm \etf \lra \Sigma \efp.$$
\end{itemize}

We make the analogous construction on sheaves, 
by starting with an analogue of the above cofibre sequence.
Indeed, we consider $\tf \tensor VA$ as the constant sheaf
of $R$-modules, and $Q(\infty tors)$ as the sum of skyscraper sheaves
for the modules $e_sTA$. We then define a sheaf $\mcD$ by 
the short exact sequence of sheaves
$$\mcD \lra \tf \tensor VA \stackrel{q}\lra Q(\infty tors)\tensorcO \Omega^* $$
of $R$-modules. Surjectivity of $q$ follows from the corresponding
fact for $R$-modules. We also note that neither $\tf \tensor VA$
nor $Q(\infty tors)$ have higher sheaf cohomology. 
Thus $\mcD$ encapsulates all the cohomology of spheres. 

\begin{remark}
Unlike the topological case, it appears that $\mcD$ is not the
dual of anything. In particular
$$\Sigma \HomcO (Q(\infty tors ), \cO)\simeq
\HomcO (Q(\infty tors ), Q (\infty tors))$$
is a proper  completion of $\mcD$. The 0th idempotent piece of its space of 
sections is of uncountable dimension, so it is different from $\mcD$.
\end{remark}

The next step in the construction is to tensor the basic short
exact sequence with the $\cO$-module $\cY$ to form
$$\mcD \tensorcO \cY \lra \tf \tensor VA\tensorcO \cY \stackrel{q}\lra 
Q(\infty tors) \tensorcO \cY .$$

To understand the central term we note that 
$VA_0=\cK =\cO (\infty tors)$. 

\begin{lemma}
For any $\cO$-module $\cY$ the sheaf $\cY (\infty tors)=\cY \tensor_{\cO }
\cO (\infty tors)$ is constant, and its cohomology is entirely
in  degree zero. \qqed
\end{lemma}

Similarly, the essential thing about  the last term is that its
cohomology is $\cE$-torsion.

\begin{lemma}
\label{torsionsheaf}
The $R$-module 
$$H^i(\cY \tensor_{\cO} Q(\infty tors) \tensor_{\cO}\Omega^*)$$
is $\cE$-torsion for $i=0$ or 1. 
\end{lemma}

\begin{proof}
Consider the decomposition $Q(\infty tors)=\bigoplus_s Q(\infty \Alr{s})$:
the $s$th term is a direct limit of terms
$Q(k\Alr{s})$ whose cohomology is annihilated by inverting $t_s$.
\end{proof}

\begin{cor}
\label{Eulerinvsheaf}
The map $\mcD \lra \tf \tensor VA$
induces an isomorphism
$$\cEi H^i(\mcD\tensor_{\cO} \cY  )\cong
H^i( \mcD \tensor_{\cO} \cO (\infty tors ) \tensorcO \cY )
=
\dichotomy{\tf \tensor \cY (\infty tors ) \tensorcO \Omega^*&\mbox{ for } i=0}
{0& \mbox{ for } i=1. }$$
\qqed
\end{cor}

\begin{defn}
\label{defn:SA}
We now define the functor
$$\SA : \tpOmod \lra \As $$
at the level of abelian categories. The object 
$$\SA (\cY)=(N\cY \lra \tf \tensor V\cY )$$ 
of $\cA_s$ associated to  a sheaf $\cY$ in degree 0 is 
$$H^*(\cY \tensor_{\cO} \mcD )
\lra H^*(\cY \tensor_{\cO} \mcD  \tensor_{\cO} 
\cO (\infty tors)).$$
To be explicit the nub is
$$N\cY_{ev}=H^0(\cY \tensor_{\cO} \mcD )$$
in even degrees, and 
$$N\cY_{od}=\Sigma^{-1}H^1(\cY \tensor_{\cO} \mcD )$$
in odd degrees. The vertex is entirely in even degrees and
$$V\cY_0=\cY (\infty tors). $$
\end{defn}

\begin{remark}
The fact that this is indeed an object of $\As$ follows from 
\ref{Eulerinvsheaf}.
Furthermore,  for any $\cO$-module $\cY$, the vertex $V\cY$ is 
entirely in even degrees. The odd degree part of the nub
is entirely $\cE$-torsion
\end{remark}

Since  tensor product is compatible with  passage to stalks, 
we may describe the divisible torsion part. 

\begin{cor}
\label{stalks}
The sheaf $\cY \tensor_{\cO} Q(\infty tors)$ is a sum of skyscraper
sheaves. Indeed, the stalk at a point of order $s$ 
is 
$$\cY \tensor_{\cO} Q(\infty tors)_s=\cY_s \tensor_{\cO_s}e_sTA_0. \qqed$$
\end{cor}

Since direct sums commute with tensor products and $\cO$ is
small, we deduce a useful formal property.

\begin{lemma}
\label{SApreservessums}
The functor $\SA$ preserves arbitrary direct sums. \qqed
\end{lemma}

\subsection{Module structure}

The formal nature of the construction gives a module structure
rather simply. 

\begin{lemma}
(i) The functor $\SA$ takes $\cO_A$ to the structure ring spectrum
$$\SA (\cO_A )=EA.$$

(ii) The functor $\SA$ takes $\cO_A$-modules to $EA$-modules, and therefore
induces a functor 
$$\SA : \tpOmod \lra \EAmod .$$
\end{lemma}

\begin{proof} 
(i) It is built into  the definition that, $\SA (\cO_A )=EA$, and
we proved in \ref{EAisring} that $EA$ is a ring. 

(ii) The sheaf level construction preserves tensor products, and 
there is a map 
$$H^i(\cY )\tensor H^i(\cZ) \lra H^i(\cY \tensorcO \cZ ) .$$
\end{proof}

One more special value will be important.

\begin{lemma}
\label{SAlinesphere}
The functor $\SA$ takes the basic line bundles to spheres:
$$\SA (\cO_A (D(V)))=S^V \sm EA.$$
\end{lemma}

\begin{proof}
The correspondence between line bundles and suspensions has
been built into the framework \ref{suspofrigid}. Thus, if we
take $\cY =\cO_A(D(V))$ we note first that 
$\cY (\infty tors)=\cK$ and  $\cY \tensorcO Q(\infty tors)=\cK /\cO (D)$. 
By construction $\cMA (\SA^{ev}\cY)(U_{\pi})$ 
is the space $\cO (D)(U_{\pi})$ of functions regular away from $\pi$. 
\end{proof}

\subsection{The functor $\protect \SA$ on torsion free sheaves.}
Whenever $\cY$ is torsion free
in the sense that it is a subsheaf of the constant 
sheaf $\cY (\infty tors)$ then the spectrum  $\SA (\cY)$ 
can be constructed exactly as we originally constructed $EA$.

\begin{defn}
If $\cY$ is an $\cO$-module we define an object 
$$\SAt (\cY )=(\tf \otimes V\cY  \lra T\cY) $$
of $\cA_t$. Here 
$$V\cY_0 =\cY (\infty tors )= \colim_{\pi} \cY (U_{\pi}), $$
and 
$$T\cY_0 =H^0(A; Q\cY), $$
 where $Q\cY$ is defined by the exact sequence
$$0 \lra \cY \lra \cY (\infty tors) \lra  Q\cY \lra 0. $$
These are made periodic with differentials as usual:
$$V\cY =V\cY_0 \tensor \omega^* \mbox{ and } T\cY =T\cY_0 \tensor \omega^*.$$
Now the structure map is defined exactly as before, 
using the differentials $Dt/t_s$. 
\end{defn}

\begin{lemma}
\label{SAontorsionfree}
If $\cY$ is torsion free then 
$$\SA (\cY) \cong (N\cY \lra \tf \tensor V\cY)$$
where $N\cY=\ker (q: \tf \tensor V\cY \lra T\cY)$.
\end{lemma}

\begin{proof}
This is immediate from the defining triangle
$$\cY \tensorcO \mcD
\lra
\cY \tensorcO \mcD \tensorcO \cO (\infty tors)
\stackrel{q}\lra 
\cY \tensorcO Q(\infty tors)\tensorcO \Omega^*.  $$
Note first that  $\cY $ is flat, being a submodule of the flat
module $\cY (\infty tors)$,  so that this is 
 a short exact sequence of sheaves. It therefore 
induces a long exact sequence in cohomology. 
Since the cohomology of $\cY (\infty tors) $
is in even degrees, it therefore suffices to show that 
the map $q$ induces a surjective map in cohomology.
\end{proof}

\section{Equivalence of $EA$-modules and $\mathcal{O}$-modules.}
\label{sec:equivalence}

We now have functors relating the algebraic model of spectra and sheaves 
over an elliptic curve. In this section we show that these can be 
combined to  give an equivalence between suitable derived categories. 

\subsection{The derived categories.}
We recall the constructions in parallel. In both cases we form the 
derived categories by a process of cellular approximation as in 
Subsection \ref{subsec:tpequiv}.

%A cellular approximation  is determined by specifying a set of spheres
%$(\sigma_{\alpha})_{\alpha \in A}$ which must be small objects. An object 
%is {\em cellular} 
%if it is built from the spheres $\sigma_{\alpha}$ using arbitrary coproducts 
%and triangles. A map $X \lra Y$ is a {\em weak equivalence} if it induces
%an isomorphism of $[\sigma_{\alpha}, \cdot ]_*$ for all $\alpha$.
%A {\em cellular approximation} of an object $X$ is then a weak equivalence
% $\Gamma X \lra X$ where $\Gamma X$ is cellular.

In the topological case, the category $D_{\T}(\EAmod )$ from 
Subsection \ref{subsec:DTEAmod}  is formed from the 
category of differential graded $EA$-modules. It is natural to  use the 
cells $EA \sm \T /H_+$
where $H$ runs through the set of closed subgroups of $\T$. However
the cofibre sequence
$$\T /\T [n]_+ \lra S^0 \lra S^{z^n}$$
shows that it is equivalent to use the cells $EA \sm S^V$ as $V$ runs
through complex representations. 
With either of these collections of cells, a map $X \lra Y$ of $EA$-modules
is a weak equivalence if and only if it induces an isomorphism 
$\pi_*^H(\cdot )$ for all closed subgroups $H$, which is the usual 
notion of  an equivariant weak equivalence of $\T$-spectra
(and equivalent to being a homology isomorphism in $\As$).

In the algebraic case, the category $\Dtp (\cOtpmod)$ from 
Subsection \ref{subsec:tpequiv} is formed from the category of 
differential graded sheaves of $\cOtp$-modules.
Motivated by the topological case, we use the cells $\cO (D(V))$
for representations $V$. It is equivalent to use the line bundles 
$\cO (D)$ where $D$ runs through torsion point divisors as was done
previously. 
A map $X \lra Y$ is then a weak equivalence if it induces an isomorphism of 
$H^*(A; \cO (-D)\tensor_{\cO}(\cdot ))$ for all torsion point divisors $D$.
%We thus write $\Dtp (\cOtpmod)$ for the resulting derived category.

\subsection{The equivalence.}

We are now equipped to state our second main theorem. 

\begin{thm}
\label{EAmodistpOmod}
The functors $\cMA : \EAmod \lra \tpOmod$ and $\SA: \tpOmod \lra 
\EAmod$  relating the categories of algebraic $EA$-module $\T$-spectra
and sheaves of $\cO$-modules defined in \ref{defn:MA} and \ref{defn:SA}
 induce an equivalence
$$D_{\T}(\EAmod )\simeq \Dtp (\cOtpmod)$$
of associated derived categories.
\end{thm}

\begin{remark}
Neither functor preserves infinite products, so this not an adjoint
pair or a Quillen equivalence. 
\end{remark}

We begin at the level of abelian categories. 

\begin{lemma}
There is a natural transformation of functors
$\cMA \SA \lra 1 $ which is a natural isomorphism on the line bundles
$\cO (D(V))$ for any complex representation $V$ with $V^{\T}=0$.
\end{lemma}

\begin{proof}
Suppose then that $\cY$ is a module, and let 
$$\SA^{ev}\cY =\left[ H^0(\cY \tensorcO \mcD )
\lra H^0(\cY \tensorcO \mcD (\infty tors)  )\right]$$
be the even summand of $\SA \cY$. Now 
$$\cMA (\SA^{ev}\cY)(U_{\pi})=
\HomAs (S^0 , \SA^{ev}\cY \tensor L(\pi)); $$
this contains the torsion part of $N^{ev}\cY$
and maps to 
$$\ker (c^0 \tensor \cY (\infty tors ) \lra H^0(\cY \tensorcO 
Q (\infty tors))\tensor R(\pi)), $$
which is $\cY (U_{\pi})$. This defines the map $\cMA \SA \cY \lra \cY$.

Now consider the sheaf $\cY = \cO (D)$. Combining \ref{MAsphereline}
and \ref{SAlinesphere}, 
we see that $\cMA \SA \cY \cong \cY$ in this case, and since $\cO (D)$ 
is torsion free, the natural transformation is the identity.
\end{proof}

\begin{lemma}
There is a natural transformation of functors
$\SA \cMA \lra 1 $ which is a natural isomorphism on the spheres $EA \sm S^W$
for any complex representation $W$.
\end{lemma}

\begin{proof}
We suppose $X=(N \stackrel{\beta}\lra \tf \tensor V)$ is an 
$EA$-module concentrated in even degrees  and construct a diagram
$$\begin{array}{ccc}
H^0(\mcD \tensor_{\cO} \cMA (X)) & \stackrel{\eta_n}\lra & N\\
\downarrow && \downarrow\\
H^0((\tf \tensor VA ) \tensor_{\cO}\cMA (X))& \stackrel{\eta_v}\lra & 
\tf \tensor V.
\end{array}$$
Since $(\tf \tensor VA)\tensor_{\cO}\cMA (X)$ is the constant sheaf at
$\tf \tensor V$, we take $\eta_v$ to be the identity, and it remains
to give a compatible definition for $\eta_n$. For this we use the
structure map $EA \sm X \lra X$ of the $EA$-module $X$. 

The  sheaf $\mcD \tensor_{\cO} \cMA (X)$ is  
associated to the presheaf given by a tensor product
of modules over each open set. By \ref{tNsheaf}, the presheaf 
$\tN$ is a sheaf with global sections $N$, so it suffices to 
construct a map at the presheaf level. More concretely, we need
maps
$$\mcD (U_{\pi}) \tensor_{\cO (U_{\pi})} \cMA (X) (U_{\pi})
\lra \tN (U_{\pi})=\cEi_{\pi}N$$
compatible under restriction. Now the domain is the tensor 
product of  $\mcD (U_{\pi})$ and 
$\cMA (X)(U_{\pi})= \Hom_{\As}(S^0, X \tensor L(\pi))$. 
The former can be identified with functions $f$ in $NA$ regular away 
from $\pi$ and the latter with elements $x \in \cEi_{\pi}N$ with 
$\beta (x) \in c^0 \tensor V$. We map this to $f\cdot x$ in $N$, 
and notice that this association is $\cO (U_{\pi})$ bilinear. 

Now if we take $X=EA \sm S^W$ we find $\cMA (X)=\cO (D(W))$
by \ref{MAsphereline} and $\SA \cMA (X)=X$ by \ref{SAlinesphere}.
We may check the natural transformation is an isomorphism stalkwise.
This is obvious for $W=0$, and for any other value, both $\cMA (X)_s$
and $N_s$ are free on the single element $t_s^{-w(s)}$.
\end{proof}

We may now complete the proof of \ref{EAmodistpOmod}.

\begin{proof}
We have defined the pair of functors $\cMA$ and $ \SA$
at the level of  abelian categories, and hence they preserve
actual homotopies at the level of differential graded 
categories. 
Accordingly they  induce functors at the level of derived categories 
by replacing objects with approximations using spheres or torsion point 
line bundles.  Since both functors take sphere objects to cellular objects
the  derived functor construction preserves  composites. Hence 
the functors $\cMA$ and $\SA$ on derived categories again provide an 
equivalence.
\end{proof}

\begin{cor}
\label{MAss}
For an $EA$-module $X$, there is a short exact sequence 
$$0 \lra \Sigma H^1(A; \cMA (X) ) \lra \pi^{\T}_* (X) \lra
H^0(A;\cMA (X) ) \lra 0.$$
\end{cor}

\begin{proof}
Indeed, from the equivalence of categories
the cohomology of $\cMA X$ is equal to the homotopy of 
$X$. The exact sequence for the cohomology of 
an  object $\cY$ in the derived category of sheaves is obtained
from the Adams resolution $\cY \lra \cI_0 \lra \cI_1$, with 
$\cI_j$ flabby.
\end{proof}

\section{Relation to Grojnowski's construction.}
\label{sec:Grojnowski}

The first construction of a $\T$-equivariant elliptic cohomology
was given by Grojnowski \cite{grojnowski}. It is defined for analytic
elliptic curves $A$, and takes values in $\Z /2$-graded sheaves over $A$. 
We first describe Grojnowski's construction and then show 
that it is related to the sheaf $\cMA F(X, EA)$
in the torsion point topology in the simplest possible way.

\subsection{Grojnowski's construction}

The construction works with 
an {\em analytic} elliptic curve $A$ over $\C$, presented as 
$$p : \C \lra \C /\Lambda =A$$
for a lattice $\Lambda \subseteq \C$. 
To each  {\em finite} $\T$-space $X$ it associates a sheaf
$\groj(X)$ over $A$ in the analytic topology.

An open set $U$ of $A$ is {\em small} if $p^{-1}U$ is the disjoint
union of connected components $V$ such that $p|_{V}: 
V \stackrel{\cong} \lra U$  is an isomorphism. 
The construction works with the analytic topology, because the
description needs to deal with small open sets. Accordingly we
let $\cO^{an}$ denote the sheaf of analytic functions on $A$ with 
the analytic topology.

Next, for a point $a \in A$ we write
$$X^a=
\dichotomy{X^{\T [s]} &\mbox{ if $a$ is of exact order $s$}}
{X^{\T } &\mbox{ if $a$ is of infinite order, }} $$
and we say that $a$ is {\em generic} if $X^a=X^{\T}$ and
$a$ is {\em special} otherwise. 

Finally, we say that an open cover $\{ U_a\}_{a \in A}$ of $A$
is {\em adapted to $X$} if the following five conditions are
satisfied
\bi
\item $a \in U_a$
\item each $U_a$ is small
\item if $a$ is special and $a \neq b$ then $a \not \in U_a \cap U_b$
\item if $a$ and $b$ are both special and $a \neq b$ then $U_a \cap U_b =
\emptyset$
\item if $b$ is generic, then $U_a \cap U_b$ is non-empty for at most one 
special $a$
\ei
For any finite $\T$-complex $X$, there is a cover adapted to $X$, and any 
two admit a common refinement. We say that the cover is $N$-discrete if
there is at most one point of order dividing $N$ in any $U_a$. 
For any finite $\T$-complex $X$ and any $N$, there is an
$N$-discrete  cover adapted to $X$, and any two admit a common refinement.

\begin{defn}(Grojnowski) 
Given an open cover $\{ U_a\}_a$ adapted to $X$ we define 
$\Z /2$-graded sheaves $\groj(X)_a$ over $U_a$ by 
$$\groj (X)_a(U)=H^*(E\T \times_{\T} X^a) \tensor_{\C [z]} \cO^{an}_A(U-a),$$
where $U-a$ is obtained by translating $U$ by $-a$, and
where $\cO^{an}(U-a)$ is a $\C [z]$-module since $z$ can be viewed as an 
analytic function on $U-a$ using $p$ to identify it with a neighbourhood
of $0 \in \C$.

These sheaves are compatible on intersections. Indeed, since the cover is 
adapted to $X$, we need only observe that the localization theorem gives an 
isomorphism 
$$\groj (X)_a|_U\cong
H^*(E\T \times_{\T} X^{\T}) \tensor_{\C [z]} \cO^{an}_A(U-a)$$
when $a \not \in  U$.
The cocycle condition is easily checked, so the sheaves 
 patch to give a sheaf $\groj (X)$ of $\cO^{an}$-algebras. This is 
independent of the adapted cover, since a refinement induces an isomorphism.

If $X$ has a $\T$-fixed basepoint $x_0$, the inclusion and projection 
induce a decomposition
$$\groj (X) = \grojtilde (X) \oplus \groj (x_0), $$
defining the reduced theory.
\end{defn}

\begin{remark}
(i) It is easy to adapt this to give a 2-periodic sheaf valued theory. 
Indeed, we need only replace $\cOan$ by $\Omegaan^*=\bigoplus_n \Omegaan^n$, 
and declare that $c \in H^2(B\T)$ acts as $z/dz$. We will do this without
change of notation, to allow comparison with our 2-periodic constructions.

(ii) The functor $\groj (X) $ is exact. Indeed a cofibre sequence 
$X' \lra X \lra X''$ induces a long exact sequence in Borel cohomology 
of $a$-fixed points, for each $a$. Since $z$ is not a zero-divisor
as an analytic function, $\tensor_{\C [z]} \cOan (U-a)$ preserves exactness.
Finally, exactness of  sequences of sheaves is detected stalkwise. 
\end{remark}

\subsection{The derived $\protect \cMA$ functor.}
 Grojnowski's functor preserves weak equivalences, so 
 we need to apply a homotopy invariant version of the functor $\cMA$. 
We therefore take  $\cMA F(X,EA)$, 
applying  the function spectrum functor rather than 
the Hom functor. The context makes clear that  $\cMA$ is to be interpreted as 
the total derived functor of the abelian category level functor.

We remark that this gives an exact functor. First note that a cofibre 
sequence $X' \lra X \lra X''$ of based $\T$-spaces induces a fibre sequence
$F(X', EA) \lla F(X, EA)\lla F(X'', EA)$
in the homotopy category of $EA$-modules. Applying the total derived
functor $\cMA$ we get a triangle in the derived category of 
$\Z/2$-graded sheaves. 

\subsection{Comparison}

In order to make the comparison we need to use the map 
$$j: \tp \lra \an$$
including the sets open in the torsion point topology amongst all open 
sets. Any sheaf in the analytic topology is a sheaf in the torsion point
topology by restriction and a sheaf $\cY$ in the torsion point topology 
gives a sheaf $j_*\cY$ in the analytic topology via
$$(j_*\cY)(U)=\colim_{U_{\pi} \supseteq U} \cY (U_{\pi}).$$
We also use the map  
$i: j_*\cO \lra \cO^{an}$ of sheaves of rings, giving a map 
$i_*$ converting $j_*\cO $-modules into $\cOan$-modules by taking
tensor products.

\begin{thm}
\label{Grojnowskicons}
The 2-periodic version of Grojnowksi's sheaf associated to a 
finite based $\T$-space $X$ is equivalent to
the sheaf arising from the function spectrum $F(X, EA)$: 
$$\grojtilde (X) \simeq i_*j_*\cMA (F(X, EA)). $$
\end{thm}

\begin{proof} First we construct a natural map 
$$\nu_X:  i_*j_*\cMA (F(X, EA)) \lra \grojtilde (X) $$
of $\cOan$-algebras. 

This corresponds to a map 
$$\nu'_X: j_*\cMA F(X,EA) \lra i^*\grojtilde (X)$$
of $j_*\cO$-algebras. For this we choose a cover $\{ U_a\}_{a\in A}$
adapted to $X$  and  construct a system of  maps 
$$\nu'_{X,a}:  (j_*\cMA F(X,EA))|_{U_a} \lra \grojtilde (X)_a
=H^*(E\T_+ \sm_{\T}X^a) \tensor_{\C [z]}\cOan (U_a-a)$$
compatible as $a$ varies.

Choose $g$ so that  all points of order $\geq g$ are generic, and let
$N=g!$. Now
choose an $N$-discrete cover $\{ U_a\}_{a\in A}$ adapted to $X$.
%We may choose $N$ so large that all points of order $\geq N$ are
%generic, so that for any $a$ the set  $U_a \setminus \{ a\}$ consists 
%entirely of points of order $\geq N$.
%  Choosing a refinement, of the adapted cover, 
%the number $N$ may be made arbitrarily large. Thus if $a$ is of order $s$ 
%(finite or $\infty$), the set $U_a$ contains all points of sufficiently
%large order.

\begin{lemma}
The map $ X^a \lra X$ induces an isomorphism
$$j_*\cMA F(X,EA)(U_a)\cong j_*\cMA F( X^a,EA)(U_a). $$
\end{lemma}

\begin{proof}
Write $\pi \cap U=\emptyset$ if $U$ contains no points with order
in $\pi$, so that $\pi \cap U=\emptyset$ if and only if $U \supseteq U_{\pi}$.
 We have 
$$ %\begin{array}{rcl}
j_*\cMA F(X,EA)(U_a)=\colim_{\pi \cap U_a=\emptyset}
%\subseteq [1, N] \setminus \{ s\}}
\cMA F(X,EA)(U_{\pi})
=[X, \tilde{E}(H \st |H| \cap U_a=\emptyset ) \sm EA]_{\T}^0 $$
%\end{array}$$
If $a$ is of order $s$, 
the quotient $X /X^a$ is built from cells $\T /\T [n]$ with 
$n$ special and $n \neq s$. Hence  $n \cap U_a =\emptyset$, 
and so the cell makes no contribution to the cohomology. 
\end{proof}

Now we may define the natural transformation as a composite
\begin{multline*}
j_*\cMA F(X,EA)(U_a)=j_*\cMA F( X^a,EA)(U_a)
\lra j_*\cMA F( E\T_+ \sm X^a,EA)(U_a)\\
\stackrel{\alpha}\lra H^*(E\T_+ \sm_{\T} X^a)\tensor_{\C [z]}\cOan(U_a-a).
\end{multline*}
To define $\alpha $, we use the fact  
that  $EA$ is almost ordinary (in the sense of \ref{EAHP}), so that 
$$F( E\T_+ \sm X^a,EA) \simeq F( E\T_+ \sm X^a,E\T_+ \sm EA)
\simeq F( E\T_+ \sm X^a,E\T_+ \sm HP). $$
Composing with projection $HP \lra H$, 
we may now complete the definition,  since there are maps 
$$\cMA F( E\T_+ \sm X^a,H)(U_{\pi})
\stackrel{\alpha}\lra H^*(E\T_+ \sm_{\T} X^a)\tensor_{\C [z]}\cOan(U_a-a)$$
for each $\pi$ so that $U_{\pi} \supseteq U$.

There are at least two ways to see that $\nu_X$ is an isomorphism for 
all $X$. Most directly, we can show that $\nu_X$ is an isomorphism
on stalks. Passing to limits, over neighbourhoods $U_a$ of $a$, 
we find
$$\begin{array}{rcl}
j_*\cMA F(X,EA)_a&=&\colim_{U_a}
[X, \tilde{E}(H \st |H|\cap U_{\pi}=\emptyset ) \sm EA]_{\T}^0\\
&=&[X, \tilde{E}(H \st |H|\neq s) \sm EA]_{\T}^*\\
&=&EA_{\T}^*(X^a)\tensor \cO_s
\end{array}$$
The completion theorem \ref{completionthm} shows what happens when we 
pass to $E\T_+ \sm X^a$  and then  we extend to analytic germs. 

The alternative is to use the fact that both sides are sheaf valued cohomology 
theories in $X$. It suffices to check that the natural map is an 
isomorphism for a class of $X$ sufficient to generate a thick category 
containing the suspension spectra of all finite complexes. 

By definition it suffices to deal with the homogeneous spaces
$\T /\T [k]_+$, and the cofibre sequences
$$\T /\T [k]_+ \lra S^0 \lra S^{z^k}$$
show it suffices to check that $\nu$ is an isomorphism for 
the spheres $S^V$. In this case all is well since
 since $\grojtilde (S^V)=\cOan (-D(V))$, 
and $\cMA F(S^V, EA)=\cO (-D(V))$. 
\end{proof}

In practical terms this gives a means for calculating the cohomology 
of $X$ using a spectral sequence from the sheaf cohomology of the 
Grojnowksi sheaf. 

\begin{cor}
There is a short exact sequence 
$$0 \lra \Sigma H^1(A;\grojtilde{X}) \lra EA_{\T}^*(X) \lra
H^0(A;\grojtilde{X}) \lra 0.$$
\end{cor}

\begin{proof}
This follows from \ref{MAss} and the fact that the cohomology
is unchanged by $i_*j_*$.
\end{proof}

\appendix

\section{The affine case: $\protect\T$-equivariant cohomology theories
from additive and multiplicative groups.}
\label{sec:GaGm}

The algebraic models of equivariant $K$-theory and Borel cohomology
are easily described \cite[13.1, 13.4]{s1q}. 
In this section we express the models
as  special cases of the  general functorial construction of
a cohomology theory $\EG_{\T}^*(\cdot )$
associated to a one dimensional affine group scheme $\G$ equipped
with a coordinate. 

%It is instructive to see how how the fact that 
%$\G$ is affine simplifies the process. 
%This will serve to illustrate the algebraic categories described
%in Section \ref{secModelT} and also complete the motivation
%of our construction for elliptic curves.

The additive group scheme $\Ga$ and the multiplicative group scheme
$\Gm$  are affine, and therefore the construction of associated
cohomology theories is considerably simpler than that for elliptic
curves. It turns out that the associated 2-periodic $\T$-equivariant 
theories are concentrated in even degrees and
$$(\EGa)_{\T}^0(X)=H^{ev}(\ET \times_{\T} X)$$
and
$$(\EGm)_{\T}^0(X)=K_{\T}^0(X), $$
and models for these theories were given in \cite{s1q}. We will repeat
the answer here in our present language.

 There are some features that differ from the elliptic case.
Once again, we must specify a coordinate $y$ on $\G$, which is a 
function whose vanishing defines $e$, or equivalently, a generator
of the augmentation ideal $(y)=\ker (\cO \lra k)$. However here we may use
the differential $dy$ to generate meromorphic differentials. 
%However the 
%construction will not change if  $y$ is multiplied by a unit of $k$, 
%so we write $[y]$ for the set of unit scalar multiples of $y$. 
Next, we must  choose functions defining the
points of order $s$ for each $s$. By definition $\G [n]$ is given
by the vanishing of $[n](y)$. The cyclotomic functions $\phi_s$  are defined
recursively by $[n](y)=\prod_{s|n}\phi_s$.
% Notice that, by contrast
%with the elliptic case, if $y=\phi_1$ is multiplied by a unit of $k$, 
%there is no change in $\phi_s$ for $s \geq 2$. 
Once again, $d\phi_s$ 
need not generate the K\"ahler differentials. For example if $s=3$
and we consider the multiplicative group with $y=1-z$, then 
$\phi_3=1+z+z^2$, and $d\phi_3=(1+2z)dz$. Since the zero of $1+2z$ is
not a point of finite order, the function $1+2z$ is not invertible.

\begin{thm}
\label{affineEG}
\label{isogeny}
Given a  commutative 1-dimensional affine group scheme $\G$ over a
ring containing $\Q$, and a coordinate  $y$ on $\G$ 
there is a 2-periodic cohomology theory  
$\EGT^*(\cdot)$ of type $\G$.
Since $\G$ is affine, the cohomology theory is complex oriented, 
 $\EGT^*$ is in even degrees and  $\G =\spec (\EGT^0)$.
The construction is natural for isomorphisms of $(\G ,y)$.

The construction  is also natural for quotient maps $p: \G \lra \G /\G[n]$ in the 
sense that  there is a map 
$p^*: \mathrm{infl}_{\T /\T[n]}^{\T}E(\G/\G[n]) \lra \EG$ 
of $\T$-spectra, where $\EG$ is viewed as a $\T /\T [n]$-spectrum and inflated
to a $\T$-spectrum, and the coordinate on $\G /\G[n]$ is  
 $\Pi_{a \in \G [n]} T_a y$, where $y$ is the coordinated on $\G$ and 
 $T_a$ denotes translation by $a$. 
\end{thm}

\begin{proof}
The construction was motivated in Section \ref{secFormal}. 
The idea is that all the ingredients described in 
Section \ref{sec:ModelT} are implicit in the definition of the type
(\ref{typeG}).

We will write down a rigid even object
$$M_t(E\G)=(\tf \tensor V\G  \stackrel{q}\lra T\G)$$
of the torsion category $\cA_t$ (i.e., the structure map $q$
will be surjective and $V\G$ and $T\G$ will be in even degrees).
By \ref{surjMtgivesMs} this is intrinsically formal and therefore determines
$$M_s(E\G )=(N\G \lra \tf \tensor V\G)$$
with $N\G =\ker (q)$, and the representing spectrum $E\G$.

Writing $\cO=\cO_{\G}$ for the ring of functions on $\G$,
in degree 0 we take
$$V\G_0 =\cOit$$
and
$$T\G _{0}=\cOit /\cO. $$
For other degrees we twist by $\omega$, taking
$$V_{2n}=V_0 \tensor \omega^n \mbox{ and } T_{2n}=T_0 \tensor \omega^n.$$

According to \ref{mapsq}, the map $q: \tf \tensor V\G \lra T\G$ 
may be described  compactly by giving its idempotent summands. 
We take
$$q (c^{w(s)} \tensor \alpha )_s:=
\overline{(\frac{\phi_s(y)}{dy})^{w(s)}\alpha}$$

A choice of coordinate $y$ gives a generator $dy^{\tensor n}$ of $\omega^n$, 
and  multiplication by  $dy$ gives
an isomorphism $\omega^{n} \lra \omega^{ (n+1)}$.
Now $\alpha \in V_{2n}$ can be written
$$\alpha =f\cdot  (dy)^{\tensor n}$$
for some function $f \in \cOit$ and 
$$q(c^{w(s)} \tensor f\cdot  (dy)^{\tensor n})_s:=
\overline{ \phi_s(y)^{w(s)}f}\cdot [(dy)^{\tensor (n-w(s))}].$$
Since any function $f$ only has finitely many poles, we see that this 
does map into the direct sum $T\G =\bigoplus_s e_sT\G$.

We must explain how $T\G$ is a module over $R$, and why $q$ is a
map of $R$-modules. We make $T\G $ into a module over $R$ by letting $c_s$ 
act as $\phi_s(y)/dy$ on $e_sT\G$.
Since poles are of finite order,   $T\G$ is a $\cE$-torsion module. The 
definition of the map  $q$ shows it is an $R$-map.

Finally, we must show that the homotopy groups of the resulting object
are as required in \ref{typeG}. By \ref{surjMtgivesMs}
we have $M_s(\EG)=(\beta :N \G \lra
\tf \otimes V\G )$, where $N\G =\ker (\tf \otimes V\G \lra T\G)$,
and we need to calculate
$$[S^W,\EG ]^{\T}_*=[S^w,M_s(\EG)]_* .$$
Since $q$ is epimorphic, $\beta $ is monomorphic, and $T\G$
is injective. Thus by \ref{surjMtgivesMs} we  have
the explicit injective resolution
$$0\lra M_s(\EG )
=\left(
\begin{array}{c}
N\G\\
\downarrow\\
\tf \otimes V\G
\end{array}
\right)
\lra
\left(
\begin{array}{c}
\tf \otimes V\G\\
\downarrow\\
\tf \otimes V\G
\end{array}
\right)
\lra
\left(
\begin{array}{c}
T\G \\
\downarrow\\
0
\end{array}
\right)
\lra 0. $$
Now, applying \ref{htpygroups} with $w(\T)=0$ we obtain the exact sequence
$$0 \lra \Hom_{\cA_s} (S^w, M_s(\EG)) \lra 
c^{-w} \tensor V\G_0 \stackrel{q^{-w}} \lra (\Sigma^{-w} T\G )_0 \lra 
\Ext_{\cA_s} (S^w, M_s(\EG)) \lra 0.$$
Hence 
$\Ext (S^w, M_s(\EG))=0$ since 
$q : c^{-w} \tensor V\G  \lra T\G$ is surjective. Indeed, any
torsion element $t \in (\Sigma^w T\G)_0\cong \cOit /\cO $ lifts to 
$f \in \cOit$ and hence to $1/c^W \otimes \chi (W) f$. 
It is immediate from the definition that if $w(\T)=0$, 
$$\Hom (S^w, M_s(\EG))=
\{ c^{-w} \otimes f \st f/\chi (W) \mbox{ regular }\}.$$
 By construction the divisor associated to the
function $\chi (V)$ is $D(V)$, so $f/\chi (V)$ is regular if and
only if $f \in \cO (-D(V))$ as required.

For the statement about isogenies, note that if $y$ is a coordinate on 
$\G$ then its norm
 $\Pi_{a \in \G [n]} T_a y$ is a coordinate on $\G /\G [n]$ 
(where $T_a$ denotes translation by $a$). Using these coordinates, 
we obtain equivariant spectra $E\G /\G [n]$ and $E\G$. As a first step to maps
between them, note that we have   maps $p_V^*:V(\G /\G[n]) \lra V\G$ and 
$p_T^*:T(\G /\G[n]) \lra T\G$ 
 corresponding to pullback of functions. However $p_V^*$ and $p_T^*$ do not 
give a map of $\T$-spectra $E(\G/\G[n]) \lra  \EG$; for example the non-equivariant
part of $E(\G /\G[n])$ corresponds to functions on $\G /\G [n]$ with support at
the identity, and these pull back to functions on $\G$ supported on  $\G [n]$, 
which correspond to the part of $E\G $ with isotropy contained in $\T [n]$. The
answer is to view the circle of equivariance of $E\G$ as  $\T /\T[n]$, and then to
 use the inflation functor studied in Chapters 10 and 24  of \cite{s1q} to obtain 
a $\T$-spectrum.

\end{proof}

\begin{remark} In the above proof we made use of the fact that the
Euler class $\chi (W)$ exists as a function in $\cK$. This should be
contrasted with the elliptic case,  where the Euler class is given by different
functions at different points. This corresponds to the fact that 
elliptic  cohomology  is not complex orientable, so that the
bundle specified by $W$ is not trivializable.
\end{remark}

We make the construction  explicit in four cases. Because the differentials
occur in the same way for all $s$, this has been omitted in the examples, and
the map $q$ translated to degree 0.

\begin{example} {\em (The additive group.)}
The ring of functions on $\Ga$ is $\Q [x]$, and the group structure is
defined by the coproduct $x \longmapsto 1 \otimes x +x \otimes 1$.
We choose $x$ as a coordinate about the identity, zero.
The group $\Ga [n]$ of points of order dividing $n$ is defined by
the vanishing of $\chi (z^n)=nx$, so the identity is the only element
of finite order over $\Q$-algebras. This case becomes rather
degenerate in that it only detects isotropy $1$ and $\T$.

The cohomology theory associated to $\Ga$ is  2-periodic Borel cohomology.
This is complex orientable, concentrated in even degrees  
and in each even degree is the map
$$\tf \otimes \cOit =\tf \otimes \Q[x,x^{-1}]
\lra \Q [x,x^{-1}]/\Q [x]=\cOit /\cO$$
$$s/e(V) \otimes f \longmapsto s\cdot \overline{f/\chi (V)}.$$
Here $\cO = \Q [x]$ and $\chi (z^n)=nx$. The ring
$\cOit=\Q [x,x^{-1}]$ of functions with poles only at points of finite
order is obtained by inverting the Euler class of $z$. \qqed
\end{example}

\begin{example} {\em (The multiplicative group.)}\cite[13.4.4]{s1q}
The multiplicative group is defined by $\Gm (k)=\units (k)$ 
with group structure given by the product. Accordingly, 
the ring of functions on $\Gm$ is $\cO =R(\T )=\Q [z,z^{-1}]$,
and the group structure is
defined by the coproduct $z \longmapsto z \otimes z$.
We choose $y=1-z$ as a coordinate about the identity element, 1.
The  coproduct then takes the more familiar form
$y \longmapsto 1 \otimes y +y \otimes 1-y \otimes y$.
The group $\Gm [n]$ of points of order dividing $n$ is defined by
the vanishing of $\chi (z^n)=1-z^n$.

The cohomology theory associated to $\Gm$ is  equivariant $K$-theory.
This is complex oriented, concentrated in even degrees and in each 
even degree  is the map
$$\tf \otimes \cOit\lra \cOit /\cO$$
$$s/e(V) \otimes f \longmapsto s\cdot \overline{f/\chi (V)}.$$
Here $\cO = \Q [z,z^{-1}]$ and $\chi (z^n)=1-z^n$. The ring
$\cOit$ of functions with poles only at points of finite
order is obtained by inverting all Euler classes. \qqed
\end{example}

\begin{example} {\em (The non-split one dimensional torus.)}
The ring of functions on the non-split (non-deploy\'e) torus 
$\Gnd$ is $\cO =\Q [a,b]/(a^2+b^2=1)$. Once one adjoins an element $i$ 
with $i^2=-1$,  this becomes equivalent to the multiplicative group 
(also known as the standard torus). Indeed, we may take $z=a+ib$ to see 
the equivalence. From the usual multiplication rule for complex numbers 
we see that the coproduct is given by $a\longmapsto a\tensor a -b\tensor b$ and
$b \longmapsto a \tensor b + b\tensor a$. The ideal $(1-a,b)$ of functions
vanishing at 0 is not principal, so there is no coordinate
in the previous sense.

Since there is no coordinate, a cohomology theory of type $\Gnd$
cannot be complex orientable. For example the map $S^0 \lra S^z$ induces
the inclusion $\cO \lla \cO (-(e))$ of functions vanishing at the
identity. Hence the cohomology of $S^z$ would not be a free module
of rank 1.

It is standard that $\Gnd$ can be revovered from $\Gm$ over $\Q (i)$
using an action of  $C_2$. Indeed, $C_2$  acts on 
$\cO =\Q (i)[z,z^{-1}]$ by the Galois action on $\Q (i)$ and by exchanging 
$z$ with $z^{-1}$. Thus $a=(z+z^{-1})/2$ and $b=i(z-z^{-1})/2$ are
fixed. The coordinate $y=1-z$ is not fixed, although
$1-a=z^{-1}(1-z)^2$ is fixed. Because the coordinate $y$ is not
fixed, the action of $C_2$ on $\cO$ does not extend to an action on 
$K\Q (i)$. 

We may construct a theory of type $\Gnd$ in the usual way. We let
$S$ denote the multiplicative set of functions $f$ with zeroes only 
at points of finite order, and take $\cK = S^{-1} \cO$. Now take
$V_0=\cK$ and $T_0=\bigoplus_{s}H^1_{\Gndlr{s}}(\cO)$. Here 
$H^1_{\Gndlr{s}}(\cO)$ is the local cohomology for the ideal of functions
vanishing at points of order exactly $s$. Now as before we define
$$q: \tf \otimes V_0 \lra T_0.  $$
For this we need to know that $\Gndlr{s}$ is {\em essentially} defined
by a principal ideal (generated by $\phi_s$ say), so that we may define 
$$q (c_s^{w(s)} \tensor f )_s= \overline{(\phi_s)^{w(s)}f} .$$
The point is that even though $\Gndlr{s}$ is not itself defined
by a principal ideal, it is in the appropriate local ring. 
For instance the ideal of functions vanishing at the identity is $(1-a,b)$. 
This is not a principal ideal, but at the level of local cohomology 
we have 
$$H^1_{(1-a,b)}(\cO)=H^1_{(1-a,b)}(\cO_{(1-a,b)})=H^1_{(1-a)}(\cO_{(1-a,b)}),$$
where the second equality follows since 
$$(1-a,b)=\sqrt{(1-a)} \mbox{ in } \cO_{(a-1,b)}, $$
as one sees explicitly from the equation $(1-a)(1+a)=b^2$. We therefore take
$\phi_1=1-a$ and define $\phi_s$ recursively by the equation
$$n^*\phi_1=\prod_{s| n}\phi_s.$$
\end{example}

\begin{example}{\em (Formal groups.)}
By way of completeness we also record the analogue for formal groups.
This completes the circle by establishing the universality of the
motivation described  in Section \ref{secFormal}. However, since
we must work over $\Q$, there is little difference from the
additive group above. Suppose given
a commutative one dimensional formal group $\Ghat$ over a ring $k$ containing
$\Q$, with a coordinate $y$.  We may identify the ring of functions on
$\Ghat$ with $k [[x]]$, and the group structure is
 the coproduct $x \longmapsto F(x \otimes 1, 1 \otimes x)$.
The group $\Ghat [n]$ of points of order dividing $n$ is defined by
the vanishing of $\chi (z^n)=[n](x)$ so the identity is the only element
of finite order over $\Q$-algebras. We may now make the direct analogue
of the construction in \ref{affineEG}.
This case becomes rather
degenerate in that it only detects isotropy $1$ and $\T$.

The cohomology theory associated to the formal group of a 
complex oriented 2-periodic cohomology theory $E$ is 
the 2-periodic Borel cohomology of $E$. This is 
concentrated in even degrees  and in each even degree is the map
$$\tf \otimes \cOit =\tf \otimes E^0((x))
\lra E^0((x))/E^0[[x]]=\cOit /\cO$$
$$s/e(V) \otimes f \longmapsto s\cdot \overline{f/\chi (V)}.$$
Here $\cO = E^0[[x]]$ and $\chi (z^n)=[n](x)$. The ring
$\cOit=E^0[[x]][1/x]=E^0((x))$ of functions with poles only at points of
finite order is obtained by inverting the Euler class of $z$. \qqed
\end{example}

%\input{ellTopsloops}

%\bibliography{bibweb}

\end{document}